\documentclass[10pt,a4paper]{article}
\pdfoutput=1 

\usepackage{lineno,hyperref}
\modulolinenumbers[5]
\usepackage[bottom=2cm, top = 2cm, left=2cm, right=2cm]{geometry}
\usepackage{amssymb}
\usepackage{amsmath}
\usepackage{algorithm}     
\usepackage{algpseudocode} 
\usepackage{subfig}
\usepackage{csquotes}
\usepackage{hyperref}
\usepackage{cleveref}
\usepackage{tikz}
\usetikzlibrary{arrows}
\usetikzlibrary{patterns}
\usepackage{mathtools}
\usepackage{bbm}
\usepackage{amsthm}
\usepackage{dutchcal}
\usepackage{pgfplots}
\usepackage{float}
\usepgfplotslibrary{external}

\newcommand{\email}[1]{\hspace*{\stretch{1}}\emph{\texttt{#1}}}

\makeatletter
\def\blfootnote{\xdef\@thefnmark{$\star$}\@footnotetext}
\makeatother
\newenvironment{Authors}%
  {\begin{center}\begin{bfseries}}%
  {\end{bfseries}\end{center}}
\newenvironment{Addresses}%
  {\begin{flushleft}\begin{itshape}}%
  {\end{itshape}\end{flushleft}}

 \usepackage{fancyhdr}  
  
\fancypagestyle{plain}{
	\fancyhead{}
	\fancyhead[C]{Accepted by Computer Methods in Applied Mechanics and Engineering,  November 2022  \hfill }
}
  
\newtheorem{theorem}{Theorem}[section]

\newtheorem{theorem2}{Theorem}[section]
\newtheorem{proposition}[theorem]{Proposition}

\newtheorem{remark}[theorem2]{Remark}

  \newcommand{\vertiii}[1]{{\left\vert\kern-0.25ex\left\vert\kern-0.25ex\left\vert #1 
    \right\vert\kern-0.25ex\right\vert\kern-0.25ex\right\vert}}

\begin{document}

\thispagestyle{plain}

\title{A one-shot overlapping Schwarz method for component-based model reduction: application to nonlinear elasticity.}
 \date{}
 
 \maketitle
\vspace{-50pt} 
 
\begin{Authors}
Angelo Iollo$^{1}$,
Giulia Sambataro$^{1}$,
Tommaso Taddei$^{1}$
\end{Authors}

\begin{Addresses}
$^1$
 IMB, UMR 5251, Univ. Bordeaux, 33400 Talence, France,
 Inria Bordeaux Sud-Ouest, Team MEMPHIS, 33400 Talence, France,\email{angelo.iollo@inria.fr, giulia.sambataro@inria.fr, tommaso.taddei@inria.fr,} \\
\end{Addresses}

\begin{abstract}
We propose a component-based (CB) parametric model order reduction (pMOR) formulation for parameterized nonlinear elliptic partial differential equations (PDEs) based on overlapping subdomains.
Our approach reads as a constrained optimization statement that penalizes the jump at the components' interfaces subject to the approximate satisfaction of the PDE  in each local subdomain.
Furthermore,  the approach relies on the decomposition of the local states into a port component --- associated with the solution on interior boundaries --- and a bubble component that vanishes at ports: 
since the bubble components are uniquely determined by the solution value at the corresponding port, we can
recast the constrained optimization statement into an unconstrained statement, which reads as a nonlinear least-squares problem and can be solved using the Gauss-Newton method.
We present thorough numerical investigations for a two-dimensional neo-Hookean nonlinear mechanics problem to validate our method; we further discuss the well-posedness of the mathematical formulation and the \emph{a priori} error analysis for linear coercive problems.
\end{abstract}


\section{Introduction}
\label{sec:intro}

\subsection{Component-based model order reduction for nonlinear PDEs}
Parametric model order reduction (pMOR, \cite{Haa17,HeRoSt16,QuMaNe16}) 
refers to a
class of computational techniques that aim at constructing a low-dimensional surrogate (or reduced-order)
model (ROM) for a given physical system, over a range of parameters. 
In the last few decades, 
pMOR techniques 
have received significant attention in science and engineering,   to speed up  parametric studies.
For complex, large-scale systems with many parameters, methods that combine  pMOR with  domain decomposition (DD) methods are of paramount importance to deal with high-dimensional parameterizations and changes in domain topology.
The aim of this work is to  present a general DD pMOR strategy for linear and nonlinear steady partial differential equations (PDEs).

Standard (monolithic) pMOR techniques
rely on high-fidelity (HF) solves at the training stage, which  might be unaffordable for very large-scale problems; furthermore, they rely  on the assumption  that the solution field is defined over a parameter-independent domain or over a family of diffeomorphic domains:
to address these issues, several authors have proposed component-based pMOR procedures (cf. \cite{HuKnPa13} and the review \cite{Buhr2020localized}).
During  the offline stage, a  library of \emph{archetype components} is defined, and local reduced-order bases (ROBs) as well as local ROMs are built; then, during the online stage, local components are instantiated to form the global system and the global solution is estimated by coupling local ROMs. 

CB-pMOR strategies consist of two distinct  building blocks:
(i) a rapid and reliable DD strategy for online global predictions, and
(ii) a localized training strategy exclusively based on local solves for the construction of the local approximations.
In this work, we focus exclusively on (i); we refer to 
\cite{benaceur2022port,smetana2022localized}
and   \cite[section 8.1.7]{hoang2021domain}
for recent works on localized training for nonlinear elliptic PDEs.

We propose a general component-based pMOR procedure for steady   PDEs based on overlapping subdomains, with a particular focus  on  second-order nonlinear elliptic PDEs. 
The key features of the approach are 
twofold:
(i)
a constrained  optimization statement that penalizes the jump at the components' interfaces  subject to the approximate (in a sense to be defined) satisfaction of the PDE  in each deployed (instantiated) component;
(ii)
the decomposition of the local solutions    into a \emph{port component} --- associated with the solution on interior boundaries (\emph{ports}) --- and a  \emph{bubble component}  that vanishes at \emph{ports}, to enable effective parallelization of the online solver.

\subsection{One-shot overlapping Schwarz method}

\begin{figure}[h!]
	\centering
	\begin{tikzpicture}[scale=0.95]
		\linethickness{0.3 mm}
		\linethickness{0.3 mm}
		\draw[ultra thick]  (0,0)--(5,0)--(5,3)--(0,3)--(0,0);
		\draw[ultra thick,densely dashed,blue]     (1.5,0)-- (1.5,3);
		\draw[ultra thick,densely dashed,red]      (3.5,0)-- (3.5,3);
		
		\coordinate [label={right:  {\Large {$\color{blue}\Gamma_2$}}}] (E) at (1.5, 1.5) ;
		
		\coordinate [label={left:  {\Large {$\color{red} \Gamma_1$}}}] (E) at (3.5, 1.5) ;
		\coordinate [label={left:  {\Large {$\color{red} \Omega_1$}}}] (E) at (1,2.5) ; 
		
		\coordinate [label={right:  {\Large {$\color{blue} \Omega_2$}}}] (E) at (4,2.5) ; 
		
		\coordinate [label={right:  {\Huge {$\Omega$}}}] (E) at (5,2.5) ; 
		
	\end{tikzpicture}
	\caption{
		configuration considered for 
		illustration in \cref{sec:intro}
		and for the  
		analysis of the linear coercive problem in \cref{sec:analysis}. 
	}
	\label{fig:analysis_explanation}
\end{figure}
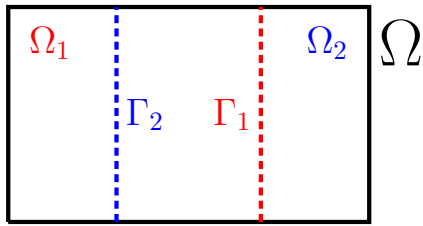

We  first introduce the formulation in the simplified case of two instantiated components  $\Omega_1,\Omega_2$ (cf. \Cref{fig:analysis_explanation}) --- to simplify notation, we do not distinguish between archetype and instantiated components; in \cref{sec:formulation}, we present the formulation in the general setting.
We denote by $\mathcal{X}_i \subset H^1(\Omega_i)$ a suitable Hilbert space in $\Omega_i$; we further define the \emph{bubble space}
$\mathcal{X}_{i,0}=\{v\in \mathcal{X}_{i}: v|_{\Gamma_i} = 0 \}$ and  the \emph{port space}
$\mathcal{U}_{i}=\{v\in \mathcal{X}_{i}: v|_{\Gamma_i} = 0 \}$, for $i=1,2$. Then, we introduce the additive or multiplicative overlapping Schwarz (OS) iterations as
\begin{equation}
	\label{eq:OS_for_dummies}
	\left\{
	\begin{array}{l}
		\displaystyle{
			{\rm find} \;\;
			u_1^{(k)} \in \mathcal{X}_1 \, :\,
			\mathcal{G}_1( u_1^{(k)} , v  ) = 0 \;\;
			\forall \, v\in \mathcal{X}_{1,0}, 
			\;\;
			u_1^{(k)} |_{\Gamma_1} = u_2^{(k-1)};
		}\\[3mm]
		\displaystyle{
			{\rm find} \;\;
			u_2^{(k)} \in \mathcal{X}_2 \, :\,
			\mathcal{G}_2( u_2^{(k)} , v  ) = 0 \;\;
			\forall \, v\in \mathcal{X}_{2,0}, 
			\;\;
			u_2^{(k)} |_{\Gamma_2} = \begin{cases*}
				u_1^{(k)}, \\
				u_1^{(k-1)}, 
			\end{cases*}
		}\\
	\end{array}
	\right.
\end{equation}
for $k=1,2,\ldots$.
Here, $u_i^{(k)}$ denotes the state estimate at the $k$-th iteration in the $i$-th subdomain, while $\mathcal{G}_1,\mathcal{G}_2$  are the variational forms associated with the PDE of interest in $\Omega_1,\Omega_2$.
\emph{Multiplicative} Schwarz iterations correspond to setting
$u_2^{(k)}|_{\Gamma_2} = 
u_1^{(k)}$  in \eqref{eq:OS_for_dummies}$_2$,
while \emph{additive} 
Schwarz iterations correspond to setting
$u_2^{(k)}|_{\Gamma_2} = 
u_1^{(k-1)}$.
Convergence of the OS iterations 
to a limit state $(u_1^{\star}, u_2^{\star})$  implies
that $\| u_1^{\star} - u_2^{\star}\|_{L^2(\Gamma_1\cup \Gamma_2)} = 0$. We thus propose to consider the formulation
\begin{equation}
	\label{eq:OS2_fordummies}
	\min_{u_1\in \mathcal{X}_1,u_2\in \mathcal{X}_2} \;
	\| u_1 - u_2  \|_{L^2(\Gamma_1\cup \Gamma_2)}
	\quad
	{\rm s.t.}  \;\;
	\mathcal{G}_i( u_i, v_i  ) = 0 \;\;
	\forall \, v_i \in \mathcal{X}_{i,0}, \;i=1,2.
\end{equation}
Clearly, the pair   $(u_1^{\star}, u_2^{\star})$ is a solution to \eqref{eq:OS2_fordummies}; in \cref{sec:analysis}, we show that, provided that the overlapping size $\delta$ is strictly positive, the solution to \eqref{eq:OS2_fordummies} is unique and depends continuously on data for linear coercive problems.
Note that, for linear problems, the solution to \eqref{eq:OS2_fordummies} can be computed directly without the need for an iterative scheme: we thus refer to our approach as  to one-shot (OS) overlapping Schwarz (OS) method and we use the abbreviation OS2\footnote{
	More rigorously, we should consider the acronym OSOS or $($OS$)^2$; however, we opted for OS2 to simplify the notation.
}. 
From this point forward, we shall use the acronym OS to refer to the standard overlapping Schwarz method.

In order to recast \eqref{eq:OS2_fordummies} into an unconstrained problem, we denote by $u_1^{\rm p},u_2^{\rm p}$ the port solutions, that is the restrictions of $u_1$ and $u_2$ to the corresponding ports; then, we introduce the extension operators
$\texttt{E}_i:\mathcal{U}_i\to \mathcal{X}_i$ and the local  \emph{port-to-bubble} solution maps 
$\texttt{F}_i:\mathcal{U}_i\to \mathcal{X}_{i,0}$ such that, given $w\in \mathcal{U}_i$, we have
$ \mathcal{G}_i(\texttt{F}_i(w) + \texttt{E}_i w , v_i  ) = 0 \;\;
\forall \, v_i \in \mathcal{X}_{i,0}$, for $i=1,2$
--- note that the \emph{port-to-bubble} field is uniquely determined by the corresponding port solution.
Then, we obtain the unconstrained OS2 statement:
\begin{equation}
	\label{eq:OS2_hybridized_fordummies}
	\min_{u_1^{\rm p}\in \mathcal{U}_1,u_2^{\rm p}\in \mathcal{U}_2} \;
	\mathfrak{f}(u_1^{\rm p}, u_2^{\rm p}):=
	\| 
	\texttt{F}_1(u_1^{\rm p}) + \texttt{E}_1 u_1^{\rm p} 
	-
	\texttt{F}_2(u_2^{\rm p}) - \texttt{E}_2 u_2^{\rm p} 
	\|_{L^2(\Gamma_1\cup \Gamma_2)}^2.
\end{equation}
The present  derivation can be viewed as a \emph{static condensation} of bubble degrees of freedom and is similar in scope to the approach in
\cite{HuKnPa13}. Following taxonomy from  the optimization literature, we might view our approach as \emph{black-box} --- as opposed to \emph{all-at-once} 
\cite[section 1.1]{herzog2010algorithms}.

Note that \eqref{eq:OS2_hybridized_fordummies} 
reads as a nonlinear least-squares problem: as in \cite{CaFaCoAm2013}, we can thus resort to the Gauss-Newton method which exploits the underlying structure of the objective function to enable rapid convergence of the CB-ROM to the optimum.

Practical implementation of a CB-pMOR approach based
on \eqref{eq:OS2_fordummies}-\eqref{eq:OS2_hybridized_fordummies} requires to address three major tasks
(i) (\emph{data compression})
the minimization statement \eqref{eq:OS2_hybridized_fordummies} is infinite-dimensional: we should thus drastically reduce the dimensionality of the port spaces $\mathcal{U}_1,\mathcal{U}_2$;
(ii) (\emph{reduction of local problems})
the local problems associated with the evaluation of the port-to-bubble maps are also infinite-dimensional: we should thus resort to standard (monolithic) MOR techniques to devise low-rank approximations of the bubble fields;
(iii) (\emph{hyper-reduction of the objective function})
evaluation of the objective function in \eqref{eq:OS2_hybridized_fordummies} requires integration over the whole curve $\Gamma_1\cup \Gamma_2$: we should thus devise a low-dimensional  quadrature rule that 
requires evaluation of the local fields in a moderate number of quadrature points.
In this work, we propose specialized MOR strategies to address these three tasks: we resort to proper orthogonal decomposition
(POD, \cite{volkwein2011model}) based on the method of snapshots \cite{Sir87} to build low-dimensional port spaces;
we rely on Galerkin ROMs 
(see, e.g., \cite{Haa17,HeRoSt16,QuMaNe16})
with hyper-reduction based on empirical quadrature/mesh sampling and weighting  \cite{farhat2015structure,yano2019lp}; 
finally, we consider two distinct approaches
to speed up the computation of the objective function: the former is based  on empirical quadrature, while the  latter relies on the empirical interpolation method (EIM, \cite{BMNP04}) .

Exploiting the static condensation of the bubble degrees of freedom, 
we can interpret the OS2 ROM as a minimum residual formulation of the port (or interface) problem associated with the underlying PDE. We discuss this interpretation for linear coercive problems in \cref{sec:analysis}.
We remark that,
similarly to \cite{lions1988schwarz},
our analysis exploits a variational interpretation of the Schwarz method.

The outline of the paper is as follows.
In \cref{sec:formulation}, we present the variational OS2 formulation for general nonlinear PDEs in arbitrary geometries.
In \cref{sec:methods}, we discuss the construction of local approximation spaces, hyper-reduction of the local models and of the objective function.
In \cref{sec:analysis},  
we discuss the well-posedness of the OS2 statement for linear coercive problems and we present an \emph{a priori} error analysis of the OS2 ROM; furthermore, 
we comment on the connection between OS and OS2 and we provide explicit convergence rates for two representative one-dimensional problems.
In \cref{sec:numerics}, we investigate performance of our method for a nonlinear elasticity problem.
\Cref{sec:conclusions} concludes the paper.

\subsection{Relation to previous works}

The aim of this work is to devise a CB-pMOR DD strategy for nonlinear PDEs: we emphasize the development of an effective solution strategy based on the Gauss-Newton method and on hyper-reduction of the objective function and of the local problems. The literature on DD for MOR and reduced-order model /full-order model (ROM-FOM) coupling is extremely vast:
CB-pMOR
strategies have  been presented in 
\cite{kaulmann2011new,HuKnPa13,hoang2021domain,iapichino2016reduced,mcbane2021component,pegolotti2021model} 
and also  recently reviewed in  \cite{Buhr2020localized};  
ROM/FOM coupling strategies have been proposed for a broad range of applications including
compressible  flows
\cite{cinquegrana2011hybrid,legresley2003dynamic,lucia2001reduced,riffaud2021dgdd}
incompressible flows
\cite{baiges2013domain,bergmann2018zonal,wicke2009modular},
and structural mechanics 
\cite{corigliano2015model,kerfriden2013partitioned,radermacher2014model} --- these methods do not distinguish between archetype and instantiated components and do not necessarily involve the training of a library of local ROMs.
Recently, several authors have proposed to couple
iterative Schwarz  DD strategies with local non-intrusive ROMs based on neural network approximations \cite{chen2021reduced,li2019d3m}.

The OS2 statement shares several features with the minimization formulation first proposed in  \cite{GEM1979} in the DD literature, for coercive linear elliptic PDEs.
OS2 is also tightly linked to the method proposed in
\cite{d2015optimization} for the coupling of local and nonlocal diffusion models (see also \cite{bochev2009optimization}): as in \cite{d2015optimization}, we interpret the OS2 statement as a control problem; while in  \cite{d2015optimization} the controls are the nonlocal volume constraint and the local boundary condition, in this work the controls are the local solutions at ports. 
We also observe that the authors of \cite{d2015optimization} do not exploit the nonlinear least-square structure of the problem and rely on a quasi-Newton scheme to approximate the solution.
We show that the choice of using the port solutions as control variables enables the definition of configuration-independent archetype components and is thus key for CB-pMOR.

Our approach is related to the Galerkin-free approach proposed in \cite{buffoni2009} and further developed in \cite{bergmann2018zonal}. In \cite{buffoni2009,bergmann2018zonal},
the authors consider a HF model in the region of interest and rely  on a low-dimensional expansion for the far-field; instead of projecting the equations in the far-field onto a low-dimensional test space, they simply rely on the objective function to compute the far-field solution coefficients (Galerkin-free). Exploiting notation introduced in the previous section, we can state  the methods in  \cite{buffoni2009, bergmann2018zonal} as:
$$
\min_{u_1\in \mathcal{X}_1, w_2 \in \mathcal{Z}_2    }
\; \|u_1 - u_2  \|_{L^2(\Omega_1 \cap \Omega_2)}
\;\;
{\rm s.t.} \;\;
\mathcal{G}_1(u_1,v_1) = 0 \; \forall \, v_1\in \mathcal{X}_{1,0},
$$
where $\mathcal{X}_1$ denotes the HF space in $\Omega_1$ and $\mathcal{Z}_2$ denotes the reduced-order space in $\Omega_2$.
The approach presented in this work is more general, more robust and also leads to more efficient online calculations, at the price of a much more involved implementation.

Our approach is linked to the minimum residual formulation in 
\cite{hoang2021domain}: the authors consider  a minimization statement  in which continuity of solution and fluxes is enforced as a constraint in the formulation, while  the global dual residual enters directly in the objective function. 
The imposition of continuity in the objective function  removes compatibility requirements at ports and allows the use of independent spaces in each archetype component;
in particular, 
the use of an overlapping partition allows us to neither explicitly enforce continuity of  the solution at ports 
nor to enforce continuity of normal fluxes. For highly-nonlinear PDEs, we found that this feature remarkably simplifies the implementation of our method and ultimately increases its flexibility.

Finally, the OS2 approach can be interpreted as an alternative to the partition-of-unity  method  (PUM,
\cite{babuvska1997partition}) considered in \cite{smetana2022localized}. Given local approximation spaces, PUM relies on the introduction of a partition of unity to define a global approximation space, and on Galerkin projection to devise the ROM for the deployed system.   PUM  has strong theoretical guarantees both in terms of  approximation and in terms of quasi-optimality properties. Similarly to OS2,  PUM requires efficient mesh interpolation to achieve online efficiency. The major difference between OS2 and PUM is that PUM relies on a global variational formulation based on a single model: on the other hand, since in OS2  local models are independent of each other, OS2 can 
be used to couple different models in different regions of the domain.


\section{Formulation}
\label{sec:formulation}
\subsection{Preliminary definitions}
\label{sec:preliminaries}

We use the superscript $(\cdot)^{\rm a}$ to indicate quantities and spaces defined for a given archetype 
component; we further denote by 
$\ell$ a generic element of the library 
$\mathcal{L}$ of archetype 
components.
We define the archetype components $\{  \Omega_{\ell}^{\rm a} \}_{\ell  \in \mathcal{L}} \subset \mathbb{R}^{d}$; we denote by 
$ {\Gamma}_{\ell}^{\rm a,dir}$ the open subset  of 
$\partial \Omega_{\ell}^{\rm a}$  where we impose Dirichlet boundary conditions, and we denote by 
$ {\Gamma}_{\ell}^{\rm a}$ the portion of 
$\partial \Omega_{\ell}^{\rm a}$ that lies inside the computational domain (``port'').
For each archetype component $\ell \in \mathcal{L}$, 
we define the local discrete 
high-fidelity (HF)
finite element (FE) space  
$\mathcal{X}_{\ell}^{\rm a} \subset  [ H_{0,    
	{\Gamma}_{\ell}^{\rm a,dir}}^1 
(   \Omega_{\ell}^{\rm a}     ) ]^D$ 
where $D$ denotes the number of state variables, 
the bubble space
$\mathcal{X}_{\ell,0}^{\rm a}  = \{ v\in  \mathcal{X}_{\ell}^{\rm a}  : v|_{{\Gamma}_{\ell}^{\rm a}}     = 0\}$, and the port space
$\mathcal{U}_{\ell}^{\rm a} = \{ 
v|_{ {\Gamma}_{\ell}^{\rm a}} : v     \in  
\mathcal{X}_{\ell}^{\rm a}  \} \subset [ H^{1/2}(   
{\Gamma}_{\ell}^{\rm a}) ]^D$.
We endow 
$\mathcal{X}_{\ell}^{\rm a}$ 
with the inner product $(\cdot, \cdot )_{\ell}$ and the induced norm
$\|  \cdot \|_{\ell} = \sqrt{ (\cdot, \cdot )_{\ell}  }$,  we define
$N_{\ell}^{\rm a} = {\rm dim} \left( 
\mathcal{X}_{\ell}^{\rm a}
\right)$, and
the extension operator
$\texttt{E}_{\ell}^{\rm a} :  \mathcal{U}_{\ell}^{\rm a}   \to \mathcal{X}_{\ell}^{\rm a}$ such that
\begin{equation}
	\label{eq:extension_operator}
	\left(  
	\texttt{E}_{\ell}^{\rm a} w, v
	\right)_{\ell} = 0 \quad \forall \, v\in 
	\mathcal{X}_{\ell,0}^{\rm a},
	\;\;
	\texttt{E}_{\ell}^{\rm a}  w \big|_{{\Gamma}_{\ell}^{\rm a}  } = w,
	\qquad \forall \, w \in  \mathcal{U}_{\ell}^{\rm a}.
\end{equation}
We define the vector of local parameters $\mu_{\ell}$ in the parameter region $\mathcal{P}_{\ell}$, which include geometric and material parameters that identify the physical model in any instantiated component of type $\ell$. We define the variational form
${\mathcal{G}}_{\ell}^{\rm a}:  
\mathcal{X}_{\ell}^{\rm a} \times \mathcal{X}_{\ell,0}^{\rm a}
\times \mathcal{P}_{\ell} \to \mathbb{R}$ such that
\begin{equation}
	\label{eq:local_model}
	\mathcal{G}_{\ell}^{\rm a} (w, v;  \mu_{\ell})
	=\sum_{k=1}^{N_{\ell}^{\rm e}  }
	\;
	\int_{\texttt{D}_{\ell,k}} \;
	\eta_{\ell}^{\rm a,e}  (w,v; \mu_{\ell}) \, dx
	\; + \;
	\int_{\partial \texttt{D}_{\ell,k}} \;
	\eta_{\ell}^{\rm a,f} (w,v; \mu_{\ell})
	\, dx
\end{equation}
where $\{ \texttt{D}_{\ell,k} \}_{k=1}^{   N_{\ell}^{\rm e}     }$ denote the elements of the FE mesh for the archetype component $\Omega_{\ell}^{\rm a}$. Furthermore, for any
$\ell\in \mathcal{L}$, 
we define the parametric mapping 
$\Phi_{\ell}^{\rm a}: 
\Omega_{\ell}^{\rm a}   \times  \mathcal{P}_{\ell}  \to \mathbb{R}^{d}$ that describes
the deformation of the archetype component
$\ell$ for the  parameter value $\mu_{\ell}\in 
\mathcal{P}_{\ell}$.

A physical system is uniquely described by a set of $N_{\rm dd}$ labels
$\{ \texttt{L}_i \}_{i=1}^{N_{\rm dd}} \subset \mathcal{L}$, and the set of parameters
${\mu}:=  (\mu_1,\ldots,\mu_{N_{\rm dd}}) \in \mathcal{P}: = \bigotimes_{i=1}^{  N_{\rm dd}   } \mathcal{P}_{\texttt{L}_i}$.
Given $ \mu \in  \mathcal{P}$,
we define 
\begin{enumerate}
	\item[(i)]
	the mappings $\{ \Phi_i \}_{i=1}^{N_{\rm dd}}$ such that $\Phi_i =\Phi_{\texttt{L}_i}^{\rm a}(\cdot; \mu_i)$ for $i=1,\ldots,N_{\rm dd}$;
	\item[(ii)]
	the instantiated overlapping partition
	$\{ \Omega_i = \Phi_i 
	(
	\Omega_{\texttt{L}_i}^{\rm a}	
	)
	\}_{i=1}^{N_{\rm dd}}$,
	the global open domain
	$\Omega \subset \mathbb{R}^d$ such that $\overline{\Omega} = \bigcup_i \overline{\Omega}_i$,  
	the ports
	$\Gamma_i = \Phi_i(   {\Gamma}_{\texttt{L}_i}^{\rm a})$ and  the Dirichlet boundaries
	$\Gamma_i^{\rm dir} = \Phi_i(  
	{\Gamma}_{\texttt{L}_i}^{\rm a,dir}	
	)$, 
	for $i=1,\ldots,N_{\rm dd}$;
	\item[(iii)]
	the deployed FE full, bubble, and port  spaces
	$\mathcal{X}_i = \{ v \circ \Phi_i^{-1} \, : \,  v \in  
	\mathcal{X}_{\texttt{L}_i}^{\rm a}  \}$, 
	$\mathcal{X}_{i,0} = \{ v\circ \Phi_i^{-1} \, : \,  v \in  
	\mathcal{X}_{\texttt{L}_i,0}^{\rm a}   \}$, and 
	$\mathcal{U}_{i} = \{ v|_{\Gamma_i}       \, : \,  v \in  
	\mathcal{X}_{i}  \}$,
	for $i=1,\ldots,N_{\rm dd}$;
	\item[(iv)]
	the extension operators
	$\texttt{E}_i: \mathcal{U}_i \to \mathcal{X}_i$ such that $\texttt{E}_i w  =  \texttt{E}_{\texttt{L}_i}^{\rm a} \left( w  \circ \Phi_i \right) \circ \Phi_i^{-1}$ for 
	$i=1,\ldots,N_{\rm dd}$;
	\item[(v)]
	the deployed variational forms
	$ \mathcal{G}_i:  
	\mathcal{X}_i \times \mathcal{X}_{i,0}\to \mathbb{R}$ such that
	\begin{equation}
		\label{eq:deployed_local_model}
		\mathcal{G}_i(w, v)
		\;   =  \;
		{\mathcal{G}}_{\texttt{L}_i}^{\rm a} (w \circ \Phi_i, v \circ \Phi_i;  \mu_i ).
	\end{equation}
\end{enumerate}
Given $i=1,\ldots,N_{\rm dd}$, 
we further define  the set of neighboring elements ${\rm Neigh}_i = \{j: \Omega_j \cap \Omega_i \neq \emptyset, j\neq i \}$, and  the partition of $\Gamma_i$ $\{ \Gamma_{i,j} = \Gamma_i \cap   \Omega_j : j\in {\rm Neigh}_i \}$ --- note that $\Gamma_{i,j}  \neq \Gamma_{j,i}$.  

Given the archetype mesh 
$\mathcal{T}_{\ell}^{\rm a}=
\left(\{x_{\ell,j}^{\rm a,v} \}_{j=1}^{N_{\ell}^{\rm v} }, 
\texttt{T}_{\ell}
\right)$, with nodes $
\{x_{\ell,j}^{\rm a,v} \}_{j=1}^{N_{\ell}^{\rm v} }
$, connectivity matrix $\texttt{T}_{\ell}$ and elements
$\{ \texttt{D}_{k,\ell}  \}_{k=1}^{N_{\ell}^{\rm e}}$, we denote by $u$ a generic element of $\mathcal{X}_{\ell}$ and we denote by $\mathbf{u}\in \mathbb{R}^{D N_{\ell}^{\rm v}}$ the corresponding FE vector associated with the Lagrangian basis of $\mathcal{T}_{\ell}^{\rm a}$, for all    $\ell  \in \mathcal{L}$. 
Following \cite{taddei2021discretize}, we pursue a discretize-then-map treatment of parameterized geometries: 
given the mesh
$\mathcal{T}_{\texttt{L}_i}^{\rm a}$, we state the local variational problems in the deformed mesh 
$\Phi_i \left( \mathcal{T}_{\texttt{L}_i}^{\rm a}  \right) = 
\left(\{  \Phi_i\left( x_{j, \texttt{L}_i}^{\rm a,v} \right) \}_{j=1}^{N_{\texttt{L}_i}^{\rm v} }, 
\texttt{T}_{\texttt{L}_i}
\right)$. In \cref{sec:hyper_port2bubble}, we discuss the hyper-reduced formulation of the local problems.
Note that if 
$\left( \mathcal{T}_{\ell}^{\rm a}, \mathbf{u} \right)$ is associated with the element $u\in \mathcal{X}_{\ell}$, then $\left(\Phi_i( \mathcal{T}_{\ell}^{\rm a} ), \mathbf{u}\right)$ approximates $u \circ \Phi^{-1}$.

\subsection{Model problem}
\label{sec:model_problem}
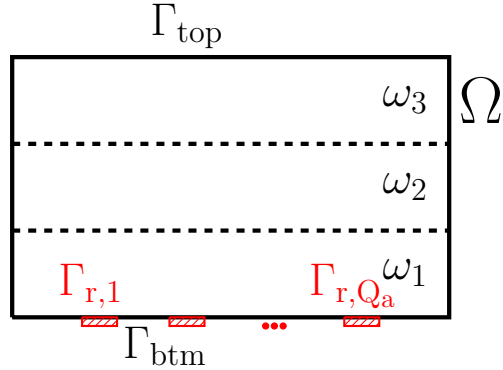
\begin{figure}[h!]
	\centering
	\begin{tikzpicture}[scale=1.15]
		\linethickness{0.3 mm}
		\linethickness{0.3 mm}
		
		\draw[ultra thick]  (0,0)--(5,0)--(5,3)--(0,3)--(0,0);
		
		\draw[ultra thick,dashed]  (0,1)--(5,1);
		\draw[ultra thick,dashed]  (0,2)--(5,2);

		\draw[thick,red,pattern=north east lines,pattern color=red]    (0.8,0)-- (1.2,0)-- (1.2,-0.1)-- (0.8,-0.1)-- (0.8,0);
		
		\draw[thick,red,pattern=north east lines,pattern color=red]    (1.8,0)-- (2.2,0)-- (2.2,-0.1)-- (1.8,-0.1)-- (1.8,0);
		
		\filldraw[color=red] (3,-0.1)    circle (1pt);
		\filldraw[color=red] (3.1,-0.1) circle (1pt);
		\filldraw[color=red] (2.9,-0.1) circle (1pt);
		\draw[thick,red,pattern=north east lines,pattern color=red]    (3.8,0)-- (4.2,0)-- (4.2,-0.1)-- (3.8,-0.1)-- (3.8,0);
		
		\coordinate [label={above:  {\LARGE {$\omega_1$}}}] (E) at (4.5, 0.25) ;
		
		\coordinate [label={above:  {\LARGE {$\Gamma_{\rm top}$}}}] (E) at (2, 3) ;
		
		\coordinate [label={below:  {\LARGE {$\Gamma_{\rm btm}$}}}] (E) at (1.75, 0) ;
		
		\coordinate [label={above:  {\LARGE {$\omega_2$}}}] (E) at (4.5, 1.25) ;
		
		\coordinate [label={above:  {\LARGE {$\omega_3$}}}] (E) at (4.5, 2.25) ;

		\coordinate [label={right:  {\Huge {$\Omega$}}}] (E) at (5, 2.5) ;
		
		\coordinate [label={above:  {\LARGE {${\color{red}\Gamma_{\rm r,1}}$}}}] (E) at (0.9, 0) ;
		
		\coordinate [label={above:  {\LARGE {$
					{\color{red} \Gamma_{\rm r,Q_{\rm a}}}$}}}] (E) at (3.9, 0) ;
	\end{tikzpicture}
	\caption{global system.
		$\Gamma_{\rm top}$ and 
		$\Gamma_{\rm r,1},\ldots, \Gamma_{\rm r,Q_{\rm a}}$ are associated with the stress conditions;
		the regions 
		$\{  \Gamma_{{\rm r},q} \}_q$ are of equal size $\ell_{\rm r}>0$, and the distance between consecutive regions is constant and equal to $\mathcal{d}>\ell_{\rm r}$.
	}
	\label{fig:picture_component_based_MOR}
\end{figure}

We illustrate the many elements of the formulation for the two-dimensional (plane stress) nonlinear (neo-Hookean) elasticity problem 
considered in the numerical experiments.
The problem shares the same geometric configuration with the problem studied in \cite{iollo2022adaptive} for radioactive management applications.
We consider the constitutive law for the first Piola Kirchhoff stress tensor
\begin{subequations}
	\begin{equation}
		\label{eq:Piola_kirchhoff}
		P(F(u)) =  
		\lambda_2 \left( F(u) - F(u)^{-T} \right)
		\, + \, 
		\lambda_1 \log \left( {\rm det} (  F(u)  ) \right)   \,  F(u)^{-T} .
	\end{equation}
	Here, $F(u) = \mathbbm{1} + \nabla u$ is the deformation gradient associated with the displacement $u$, $\lambda_1, \lambda_2$ are the Lam{\'e} constants given by 
	\begin{equation}
		\label{eq:lame_constants}
		\lambda_1 = \frac{E\nu}{1-\nu^2},
		\qquad
		\lambda_2 = \frac{E }{2(1 + \nu)},
	\end{equation}
	where $E$ is the Young's modulus, and $\nu$ is the Poisson's ratio.
	We consider the domain $\Omega=(0,1)^2$ depicted in Figure \ref{fig:picture_component_based_MOR}; we set $\nu=0.3$ and we consider $E=E_k$ in $\omega_k$ for $k=1,2,3$.
	We prescribe normal homogeneous Dirichlet conditions on  the left and right boundaries; homogeneous Dirichlet conditions on the bottom boundary $\Gamma_{\rm btm}$ and the Neumann  conditions:
	\begin{equation}
		\label{eq:neumann_conditions}
		P(F(u))  \mathbf{n} \big|_{ \Gamma_{\rm top}   }
		=
		g_{\rm top} :=
		\left[
		\begin{array}{l}
			0 \\
			-4 x_1(1-x_1) \\
		\end{array}
		\right],
		\quad
		P(F(u))  \mathbf{n} \big|_{ \Gamma_{\rm r,q}   }
		=
		g_{\rm r} :=
		-s
		\left[
		\begin{array}{l}
			0 \\
			1 \\
		\end{array}
		\right],
		\quad
		q=1,\ldots,Q_{\rm a}
	\end{equation}
	with $s>0$. 
\end{subequations}

The system of equations below summarizes the problem: we seek the solution $u:\Omega \to \mathbb{R}^2$ to the system 
\begin{equation}
	\label{eq:PDE_model}
	\left\{
	\begin{array}{ll}
		-\nabla \cdot P(F(u)) = 0 & {\rm in} \; \Omega \\[3mm]
		u \cdot \mathbf{n}= 0 & {\rm on} \; \{ 0 , 1 \} \times (0,1) \\[3mm]
		P(F(u))   \mathbf{n} = g_{\rm r}
		& \displaystyle{{\rm on} \; \Gamma_{\rm r} }\\[3mm]
		P(F(u))  \mathbf{n} = g_{\rm top}
		& \displaystyle{{\rm on} \;   \Gamma_{\rm top}}\\[3mm]
		u = 0
		& \displaystyle{{\rm on} \;   \Gamma_{\rm btm} = (0,1) \times \{0\} \setminus  \Gamma_{\rm r} }
		\\ 
	\end{array}
	\right.
\end{equation}
where $\Gamma_{\rm r} = \bigcup_{q=1}^{Q_{\rm a}} \Gamma_{\rm r,q}$. Our goal is to estimate the solution to \eqref{eq:PDE_model} for any choice of the  Young's moduli 
$(E_1,E_2,E_3)$ associated with the  regions $\omega_1,\omega_2,\omega_3$  in $ [25,30]\times [10,20]\times [10,20]$,
any value of $s\in [0.4,1]$  in \eqref{eq:neumann_conditions}, and
any $Q_{\rm a} \in \{2,\ldots,7\}$. Note that variations of $Q_{\rm a}$ induce topological changes that prevent the application  of standard monolithic techniques.

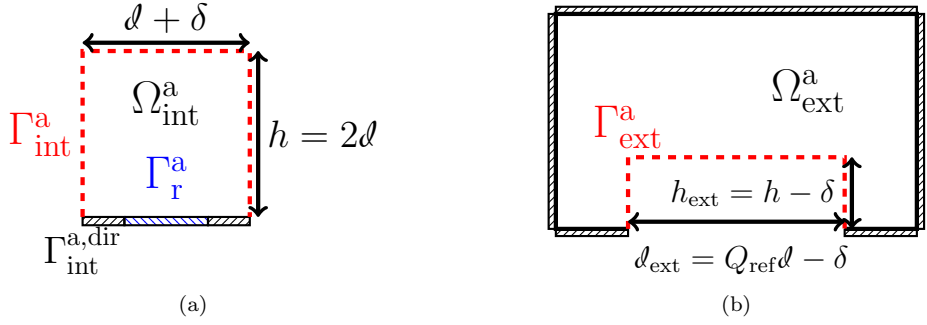
\begin{figure}[h!]
	\centering
	\subfloat[ ]{
		\begin{tikzpicture}[scale=1.1]
			\linethickness{0.3 mm}
			\linethickness{0.3 mm}
			
			\coordinate [label={above:  {\LARGE {${\Omega}_{\rm int}^{\rm a}$}}}] (E) at (1, 1) ;
			
			\coordinate [label={above:  {\LARGE {$\color{blue}{ {\Gamma}_{\rm r}^{\rm a} }$}}}] (E) at (1, 0.1) ;
			
			\draw[thick,pattern=north east lines,pattern color=black]    (0,0)-- (0.5,0)-- (0.5,-0.1)-- (0,-0.1)-- (0,0);
			
			\draw[thick,pattern=north west lines,pattern color=blue]    (0.5,0)-- (1.5,0)-- (1.5,-0.1)-- (0.5,-0.1)-- (0.5,0);
			
			\draw[thick,pattern=north east lines,pattern color=black]    (1.5,0)-- (2,0)-- (2,-0.1)-- (1.5,-0.1)-- (1.5,0);
			
			\draw[ultra thick,dashed,red]  (0,0)--(0,2) --(2,2)--(2,0);
			
			\coordinate [label={left:  {\LARGE {$\color{red}{ \Gamma_{\rm int}^{\rm a}}$}}}] (E) at (0, 1) ;
			
			\draw[<->,ultra thick]  (0,2.1)--(2,2.1);
			\coordinate [label={above:  {\Large {$\mathcal{d}+\delta$}}}] (E) at (1, 2.1) ;
			
			\draw[<->,ultra thick]  (2.1,0)--(2.1,2);
			\coordinate [label={right:  {\Large {$h = 2 \mathcal{d}$}}}] (E) at (2.1, 1) ;
			
			\coordinate [label={below:  {\Large {$\color{black}{ {\Gamma}^{\rm a,dir}_{\rm int}}$}}}] (E) at (0, 0) ;
		\end{tikzpicture}
	}
	~~~~~~~~~~~~~~~
	\subfloat[ ]{
		\begin{tikzpicture}[scale=0.95]
			\linethickness{0.3 mm}
			\linethickness{0.3 mm}
			
			\draw[ultra thick] (1,0)-- (0,0)--(0,3)-- (5,3)--(5,0)--(4,0);
			
			\draw[thick,pattern=north east lines,pattern color=black]    (0,0)-- (1,0)-- (1,-0.1)-- (0,-0.1)-- (0,0);
			
			\draw[thick,pattern=north east lines,pattern color=black]    (4,0)-- (5,0)-- (5,-0.1)-- (4,-0.1)-- (4,0);
			
			\draw[thick,pattern=north east lines,pattern color=black]    (0,0)-- (0,3)-- (-0.1,3)-- (-0.1,0)-- (0,0);
			
			\draw[thick,pattern=north east lines,pattern color=black]    (5,0)-- (5,3)-- (5.1,3)-- (5.1,0)-- (5,0);
			
			\draw[thick,pattern=north east lines,pattern color=black]    (0,3)-- (5,3)-- (5,3.1)-- (0,3.1)-- (0,3);
			
			\draw[ultra thick,dashed,red]   (1,0)-- (1,1)-- (4,1)-- (4,0);
			
			\coordinate [label={above:  {\LARGE {$\Omega_{\rm ext}^{\rm a}$}}}] (E) at (3.5, 1.5) ;
			
			\coordinate [label={above:  {\LARGE {$\color{red}{\Gamma_{\rm ext}^{\rm a}}$}}}] (E) at (1, 1) ;
			
			\draw[<->,ultra thick]  (1,0.1)--(4,0.1);
			\coordinate [label={below:  {\large {
						$\mathcal{d}_{\rm ext} = Q_{\rm ref} \mathcal{d} - \delta$}}}] (E) at (2.5,- 0.1) ;
			
			\draw[<->,ultra thick]  (4.1,0)--(4.1,1);
			\coordinate [label={left:  {\large {
						$h_{\rm ext} =h - \delta$}}}] (E) at (4,0.5) ;
			
		\end{tikzpicture}
	}
	
	\caption{geometrical configuration. Archetype components. ($Q_{\rm ref}=5$).}
	\label{fig:picture_components}
\end{figure}

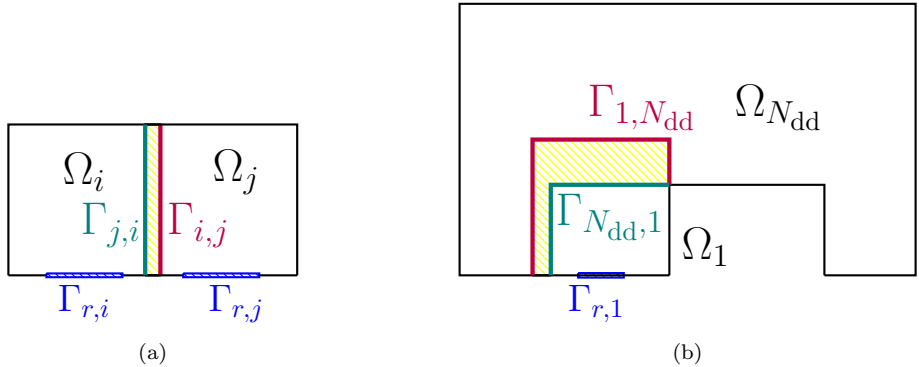
\begin{figure}
	[h!]
	\centering
	\subfloat[ ]{
		\begin{tikzpicture}[scale=1]
			\linethickness{0.3 mm}
			\linethickness{0.3 mm}
			
			\draw[thick,gray,pattern=north west lines,pattern color=yellow]    (1.8,2)--(2,2)--(2,0)--(1.8,0);

			\draw[thick,pattern=north east lines,pattern color=black]    (0,0)-- (0.5,0);
			\draw[thick,pattern=north east lines,pattern color=black]    (1.5,0)-- (2,0);
			\draw[thick,blue,pattern=north west lines,pattern color=blue]    (0.5,0.025)-- (1.5,0.025)--(1.5,-0.025)--(0.5,-0.025)--(0.5,0.025);
			\draw[thick,black]  (0,0)--(0,2) --(2,2)--(2,0);
			\coordinate [label={below:  {\Large {$\color{blue}{\Gamma_{r,i}}$}}}] (E) at (1, 0) ;
			
			\coordinate [label={above:  {\LARGE {${\Omega}_{i}$}}}] (E) at (1, 1) ;
			
			\draw[ultra thick,purple]  (2,2)--(2,0);
			\coordinate [label={above:  {\LARGE {$\color{purple}{\Gamma_{i,j}}$}}}] (E) at (2.5, 0.25) ;
			
			\draw[thick,pattern=north east lines,pattern color=black]    (1.8,0)-- (2.3,0);
			\draw[thick,blue,pattern=north west lines,pattern color=blue]    (2.3,0.025)-- (3.3,0.025)--(3.3,-0.025)--(2.3,-0.025)--(2.3,0.025);
			\draw[thick,pattern=north east lines,pattern color=black]    (3.3,0)-- (3.8,0);
			\draw[thick,black]  (1.8,0)--(1.8,2) --(3.8,2)--(3.8,0);
			\coordinate [label={below:  {\Large {$\color{blue}{\Gamma_{r,j}}$}}}] (E) at (3, 0) ;
			
			\coordinate [label={above:  {\LARGE {${\Omega}_{j}$}}}] (E) at (3, 0.95) ;
			
			\draw[ultra thick,teal]  (1.8,2)--(1.8,0);
			\coordinate [label={above:  {\LARGE {$\color{teal}{\Gamma_{j,i}}$}}}] (E) at (1.35, 0.25) ;
			
		\end{tikzpicture}
	}
	~~~~~~~~~~~~~~~
	\subfloat[ ]{
		\begin{tikzpicture}[scale=1.2]
			\linethickness{0.3 mm}
			\linethickness{0.3 mm}
			
			\draw[thick,gray,pattern=north west lines,pattern color=yellow]   (0.8,0)--(0.8,1.5)--(2.3,1.5)--(2.3,1)--(1,1)--(1,0)--(0.8,0);

			\draw[thick,black] (1,0)-- (0,0)--(0,3)-- (5,3)--(5,0)--(4,0);
			\draw[thick,black]   (1,0)-- (1,1)-- (4,1)-- (4,0);
			\coordinate [label={above:  {\LARGE {$\Omega_{N_{\rm dd}}$}}}] (E) at (3.5, 1.5) ;
			
			\coordinate [label={above:  {\LARGE {$\Omega_{1}$}}}] (E) at (2.7, 0) ;
			\draw[thick,black]  (0.8,0)--(0.8,1.5) --(2.3,1.5)--(2.3,0);
			\draw[thick,blue,pattern=north west lines,pattern color=blue]    (1.3,0.025)-- (1.8,0.025)--(1.8,-0.025)--(1.3,-0.025)--(1.3,0.025);
			\draw[thick,black] (0.8,0)--(1.3,0);
			\draw[thick,black] (1.3,0)--(2.3,0);
			\coordinate [label={below:  {\Large {$\color{blue}{\Gamma_{r,1}}$}}}] (E) at (1.5, 0) ;
			
			\draw[ultra thick,purple]  (0.8,0)--(0.8,1.5)--(2.3,1.5)--(2.3,1);
			\coordinate [label={above:  {\LARGE {$\color{purple}{\Gamma_{1,N_{\rm dd}}}$}}}] (E) at (2,1.5) ;
			
			\draw[ultra thick,teal]  (1,0)--(1,1)--(2.3,1);
			\coordinate [label={below:  {\LARGE {$\color{teal}{\Gamma_{N_{\rm dd},1}}$}}}] (E) at (1.65, 1) ;
		\end{tikzpicture}
	}
	\caption{geometrical configuration. Examples of deployed components. (a): $i,j=1, \ldots, Q_{\rm a}$, (b): $i=1,\,  j=N_{\rm dd}=Q_{\rm a}+1$. The overlap area is marked in yellow.}
	\label{fig:picture_deployed_components}
\end{figure}

We introduce the library of components 
${\Omega}_{\rm int}^{\rm a}$ and ${\Omega}_{\rm ext}^{\rm a}$ depicted in   \Cref{fig:picture_components}; in \Cref{fig:picture_deployed_components} we show examples of instantiated components and we identify the corresponding ports.
We denote by $\delta>0$ the size of the overlap.
The mapping ${\Phi}_{\rm int}^{\rm a}$ associated with the internal component is a simple horizontal shift, while the mapping ${\Phi}_{\rm ext}^{\rm a}$ associated with the external component consists in a piecewise-linear map in the horizontal direction and the identity map in the vertical direction.
The internal component is uniquely described by the vector of parameters $\mu_{\rm int} = [E_1,s, x_{\rm shift}]$ where $x_{\rm shift}$ denotes the magnitude of the horizontal shift;
the external component is described by the vector of parameters
$\mu_{\rm ext} = [E_1,E_2, E_3, \mathcal{d}_{\rm ext}]$ with 
$\mathcal{d}_{\rm ext} = Q_{\rm a} \mathcal{d} - \delta$. Note that the external archetype component (cf. Figure~\ref{fig:picture_components}) corresponds to the choice $Q_{\rm a}=Q_{\rm ref}$ with $Q_{\rm ref}=5$. 
We then introduce the variational forms:
\begin{subequations}
	\label{eq:local_models_neohook}
	\begin{equation}
		\left\{
		\begin{array}{l}
			\displaystyle{
				{\mathcal{G}}_{\rm int}^{\rm a} (w,v; \mu_{\rm int})
				=
				\int_{\Omega_{\rm int}^{\rm a}}
				\eta_{\rm int}^{\rm a,e} (w,v; \mu_{\rm int}) \, dx
				\,+\,
				\int_{ \Gamma_{\rm r}^{\rm a}}
				\eta_{\rm int}^{\rm a,f} (w,v; \mu_{\rm int}) \, dx,
			}\\[3mm]
			\displaystyle{
				{\mathcal{G}}_{\rm ext}^{\rm a} (w,v; \mu_{\rm ext})
				=
				\int_{\Omega_{\rm ext}^{\rm a}}
				\eta_{\rm ext}^{\rm a,e} (w,v; \mu_{\rm int}) \, dx
				\,+\,
				\int_{ \Gamma_{\rm top}^{\rm a}}
				\eta_{\rm ext}^{\rm a,f} (w,v; \mu_{\rm ext}) \, dx.
			}
			\\
		\end{array}
		\right.
	\end{equation}
	Explicit expressions of $\eta_{\ell}^{\rm a,e}$ and $\eta_{\ell}^{\rm a,f}$ can be obtained  by resorting  to change-of-variable formulas: given the mapping $\Phi$, we denote by $\nabla_{\Phi} = \nabla \Phi^{-T} \nabla$ the corresponding ``mapped'' gradient and we define 
	$\nabla_{\rm s, \Phi} =\frac{1}{2} \left(
	\nabla_{\Phi} + \nabla_{\Phi}^{T}
	\right)$ and $F_{\Phi}  = \mathbbm{1} + \nabla_{\Phi}$. Then, we have
	(we omit dependence on the parameter to shorten notation)
	\begin{equation}
		\begin{array}{l}
			\displaystyle{
				\eta_{\rm int}^{\rm a,e} (w,v)
				=
				\eta_{\rm ext}^{\rm a,e} (w,v)
				=
				P( F_{\Phi}(w)    ) : \nabla_{\rm s, \Phi} v \, {\rm det} (\nabla \Phi),
			}\\[3mm]
			\displaystyle{
				\eta_{\rm int}^{\rm a,f} (w,v)
				=
				v\cdot \left( g_{\rm r} \circ \Phi  \right)\| \nabla \Phi \widehat{\mathbf{t}}  \|_2 ,
				\quad
				\eta_{\rm ext}^{\rm a,f} (w,v)
				=
				v \cdot  \left(  g_{\rm top} \circ \Phi \right) \| \nabla \Phi \widehat{\mathbf{t}}  \|_2,
			}
			\\
		\end{array}
	\end{equation}
	where  
	$\widehat{\mathbf{t}} $ denotes the tangent vector to the surface.
\end{subequations}

\subsection{Hybridized statement}
\label{sec:hybridization}

\subsubsection{High-dimensional formulation}
We generalize below the OS2 statement introduced in \cref{sec:intro}.
Given the set of parameters
${\mu} =  (\mu_1,\ldots,\mu_{N_{\rm dd}}) \in \mathcal{P}  = \bigotimes_{i=1}^{  N_{\rm dd}   } \mathcal{P}_{\texttt{L}_i}$, we 
propose the CB full-order model: find
$u^{\rm hf} = (u_1^{\rm hf},\ldots,u_{N_{\rm dd}}^{\rm hf}) \in 
{\mathcal{X}} := 
\bigotimes_{i=1}^{  N_{\rm dd}   } \mathcal{X}_i$ to minimize
\begin{equation}
	\label{eq:OS2_original}
	\min_{u \in  {\mathcal{X}}}
	\;
	\frac{1}{2}  \sum_{i=1}^{N_{\rm dd}} 
	\sum_{j \in {\rm Neigh}_i} \;
	\|u_i - u_j \|_{L^2(\Gamma_{i,j})}^2
	\quad
	{\rm s.t.} \;\;
	\mathcal{G}_i(u_i,v_i) = 0 \;\; \forall \, v_i \in \mathcal{X}_{i,0}, \; i=1,\ldots,N_{\rm dd}.
\end{equation}
Note that \eqref{eq:OS2_original} reduces to \eqref{eq:OS2_fordummies} for the case of two overlapping components.

To derive the hybridized formulation, we define the port-to-bubble maps $\texttt{F}_i: \mathcal{U}_i \to \mathcal{X}_{i,0}$ such that, given $w\in  \mathcal{U}_i$, 
\begin{subequations}
	\label{eq:hybridizedOS2}
	\begin{equation} 
		\label{eq:port2bubble}
		\mathcal{G}_i\left(  \texttt{F}_i(w) + \texttt{E}_i w, v \right) \, = \, 0 \quad
		\forall \, v\in \mathcal{X}_{i,0}. 
	\end{equation}
	Note that \eqref{eq:port2bubble} corresponds to the FE solution to a localized PDE problem with datum $w$ on $\Gamma_i$. Then, we rewrite \eqref{eq:OS2_original} as the unconstrained least-square problem:
	find $u^{\rm hf,p}  = (u_1^{\rm hf,p},\ldots,u_{N_{\rm dd}}^{\rm hf,p}) \in 
	{\mathcal{U}} := 
	\bigotimes_{i=1}^{  N_{\rm dd}   } \mathcal{U}_i$ to minimize
	\begin{equation}
		\label{eq:hybridizedOS2_b}
		\min_{u^{\rm p}  \in {\mathcal{U}}}
		\;
		\frac{1}{2}  \sum_{i=1}^{N_{\rm dd}} 
		\sum_{j \in {\rm Neigh}_i} \;
		\|   u_i^{\rm p}    -  \texttt{E}_j  u_j^{\rm p}  - \texttt{F}_j(u_j^{\rm p} ) 
		\|_{L^2(\Gamma_{i,j})}^2.
	\end{equation}
	Minimization problem \eqref{eq:hybridizedOS2_b} reads as a nonlinear least-square problem; in the following we devise a low-dimensional 
	reduced-order 
	approximation 
	of  \eqref{eq:hybridizedOS2_b}
	based on Galerkin projection of the port-to-bubble maps.
\end{subequations}

\subsubsection{Reduced-order formulation}
For all $\ell \in \mathcal{L}$, 
we introduce the low-dimensional archetype bubble and port spaces 
$\mathcal{Z}_{\ell}^{\rm a,b} \subset \mathcal{X}_{\ell,0}$,
$\mathcal{Z}_{\ell}^{\rm a,p} \subset 
\mathcal{U}_{\ell}$ and the extended port spaces
$\mathcal{W}_{\ell}^{\rm a,p} = \{  
\texttt{E}_{\ell} \zeta : \zeta \in 
\mathcal{Z}_{\ell}^{\rm a,p} 
\} \subset 
\mathcal{X}_{\ell}$; we denote by $n$ and $m$ the dimensions of the bubble and port spaces, respectively; for simplicity, we assume that the dimension of the spaces is the same for all archetype components. 
We also define the archetype ROBs 
$Z_{\ell}^{\rm a,b}:\mathbb{R}^n \to \mathcal{Z}_{\ell}^{\rm a,b}$ and
$W_{\ell}^{\rm a,b}:\mathbb{R}^m \to \mathcal{W}_{\ell}^{\rm a,b}$. 
Given the deployed system, we introduce the instantiated (or deployed) bubble and port spaces
$\mathcal{Z}_i^{\rm b} = \{ \zeta\circ \Phi_i^{-1}: \zeta\in  \mathcal{Z}_{\texttt{L}_i}^{\rm a,b}
\}$ and
$\mathcal{W}_i^{\rm p} = \{ \zeta\circ \Phi_i^{-1}: \zeta\in  \mathcal{W}_{\texttt{L}_i}^{\rm a,p} \}$ with 
ROBs 
$Z_i^{\rm b}=[\zeta_{i,1}^{\rm b},\ldots,  \zeta_{i,n}^{\rm b}]: \mathbb{R}^n \to \mathcal{Z}_i^{\rm b}$ 
and
$W_i^{\rm b}=[\psi_{i,1}^{\rm p},\ldots,  \psi_{i,m}^{\rm p}]: \mathbb{R}^m \to \mathcal{W}_i^{\rm p}$, respectively. 
Then, we define the ansatz:
\begin{equation}
	\label{eq:ansatz_local}
	\widehat{u}_i ( \widehat{\boldsymbol{\alpha}}_i ,   
	\widehat{\boldsymbol{\beta}}_i )= Z_i^{\rm b} \widehat{\boldsymbol{\alpha}}_i +
	W_i^{\rm p} \widehat{\boldsymbol{\beta}}_i,
	\quad
	i=1,\ldots,N_{\rm dd}.
\end{equation}
We observe that $ \widehat{u}_i^{\rm b} = Z_i^{\rm b} \widehat{\boldsymbol{\alpha}}_i$ should approximate the bubble field
$u|_{\Omega_i} - \texttt{E}_i ( u|_{\Gamma_i} )$, while
$ \widehat{u}_i^{\rm p} = W_i^{\rm p} \widehat{\boldsymbol{\beta}}_i$ is an approximation of the (extended) port field $\texttt{E}_i ( u|_{\Gamma_i} )$: we refer to 
$ \widehat{u}_i^{\rm b},\widehat{u}_i^{\rm p}$ as to the bubble and port estimates of the solution field in the $i$-th component.

\begin{subequations}
	\label{eq:port2bubble_all}
	To obtain the low-dimensional formulation, 
	we introduce the local residuals\footnote{The superscript $^{\rm hf}$ encodes the fact that the local residuals are computed using the HF mesh.}
	(cf. \eqref{eq:local_model}, \eqref{eq:deployed_local_model} and \eqref{eq:local_models_neohook})
	\begin{equation}
		\label{eq:local_problems}
		\widehat{\mathbf{R}}_i^{\rm hf} : \mathbb{R}^n\times \mathbb{R}^m \to \mathbb{R}^n
		\quad
		{\rm s.t.} \;\;
		\left(
		\widehat{\mathbf{R}}_i^{\rm hf} \left( {\boldsymbol{\alpha}}_i , {\boldsymbol{\beta}}_i \right)
		\right)_j
		=
		\mathcal{G}_i \left( \widehat{u}_i( {\boldsymbol{\alpha}}_i,   {\boldsymbol{\beta}}_i )  \,, \, 
		\zeta_{i,j}^{\rm b} \right),
		\;\;
		i=1,\ldots,N_{\rm dd}, \;\;
		j=1,\ldots,n,
	\end{equation}
	and the approximate port-to-bubble maps
	$\widehat{\texttt{F}}_i^{\rm hf}:\mathbb{R}^m \to \mathbb{R}^n$ such that
	$
	\widehat{\mathbf{R}}_i^{\rm hf} 
	\left(
	\widehat{\texttt{F}}_i^{\rm hf}
	\left(
	{\boldsymbol{\beta}}_i\right) , {\boldsymbol{\beta}}_i \right)
	= \mathbf{0}.
	$
	Computation of the port-to-bubble maps $\{ \widehat{\texttt{F}}_i^{\rm hf}\}_i$ is expensive due to the need to integrate over the whole computational mesh. We thus replace the residuals $\{ \widehat{\mathbf{R}}_i^{\rm hf}\}_i$ with the empirical quadrature (EQ) approximations   $\{ \widehat{\mathbf{R}}_i^{\rm eq}\}_i$ and we define the hyper-reduced  port-to-bubble maps
	$\widehat{\texttt{F}}_i^{\rm eq}:\mathbb{R}^m \to \mathbb{R}^n$ such that
	\begin{equation}
		\label{eq:port2bubble_reduced}
		\widehat{\mathbf{R}}_i^{\rm eq} 
		\left(
		\widehat{\texttt{F}}_i^{\rm eq}
		\left(
		{\boldsymbol{\beta}}_i\right) , {\boldsymbol{\beta}}_i \right)
		= \mathbf{0}.
	\end{equation}
	We discuss in \cref{sec:hyper_port2bubble} the hyper-reduction strategy employed to construct the approximate residuals $\widehat{\mathbf{R}}_i^{\rm eq} $; here, we observe that
	the gradient of the port-to-bubble map can be obtained by differentiating \eqref{eq:port2bubble_reduced}:
	\begin{equation}
		\label{eq:port2bubble_reduced_gradient}
		\nabla 
		\texttt{F}_i^{\rm eq}
		\left(
		\boldsymbol{\beta}_i\right)
		= \, -\,
		\left(
		\partial_{ \boldsymbol{\alpha}_i   }
		\widehat{\mathbf{R}}_i^{\rm eq} \right)^{-1}
		\partial_{ \boldsymbol{\beta}_i   }
		\widehat{\mathbf{R}}_i^{\rm eq}
		\Big|_{( \boldsymbol{\alpha}_i, \boldsymbol{\beta}_i     ) \, = \, (  \widehat{\texttt{F}}_i^{\rm eq}
			\left( {\boldsymbol{\beta}}_i\right),  {\boldsymbol{\beta}}_i    )  } 
	\end{equation}
\end{subequations}
We remark that the existence and well-posedness of  the port-to-bubble maps
\eqref{eq:port2bubble_reduced}
is conditioned to the existence of solutions to  the nonlinear systems of equations 
$\widehat{\mathbf{R}}_i^{\rm hf} = \mathbf{0}$ 
and to the fact that
$\partial_{\boldsymbol{\alpha}_i} \widehat{\mathbf{R}}_i^{\rm hf}$ is non-singular at the optimum. It thus depends on the particular problem of interest, and might also depend on the overlapping partition considered and on the reduced-order approximation spaces.

We now focus on the objective function. We observe that
$$
\begin{array}{l}
	\displaystyle{
		\frac{1}{2}
		\sum_{i=1}^{N_{\rm dd}}
		\sum_{j \in {\rm Neigh}_i} 
		\int_{\Gamma_{i,j}} \| \widehat{u}_i(x) - \widehat{u}_j(x)   \|_2^2 \, dx
	} 
	\\[3mm]
	=
	\displaystyle{
		\frac{1}{2}
		\sum_{i=1}^{N_{\rm dd}}
		\, \int_{\Gamma_{i}} 
		\left( 
		\sum_{j \in {\rm Neigh}_i: x\in \Omega_j} 
		\| \widehat{u}_i(x) - \widehat{u}_j(x)   \|_2^2 
		\right)
		\, dx
	} 
	\\[3mm]
	=
	\displaystyle{
		\frac{1}{2}
		\sum_{i=1}^{N_{\rm dd}}
		\, \int_{ {\Gamma}_{\texttt{L}_i}^{\rm a}} 
		\left( 
		\sum_{j \in {\rm Neigh}_i: \Phi_i(\widehat{x}) \in \Omega_j} 
		\| \widehat{u}_i(  \Phi_i(\widehat{x}) )- \widehat{u}_j(\Phi_i(\widehat{x}))   \|_2^2 
		\right)
		J_i^{\rm bnd}(  \widehat{x}  )
		\, d \widehat{x}
	} 
	\\
\end{array}
$$
where $J_i^{\rm bnd}  = \|
{\rm det} ( \nabla \Phi_i  )
\nabla \Phi_i^{-T}  {\mathbf{n}}_{\texttt{L}_i}^{\rm a} \|_2$ and 
${\mathbf{n}}_{\ell}^{\rm a}$ is the outward normal to ${\Gamma}_{\ell}^{\rm a}$. Note that in the last identity we used the Narson formula; furthermore, to shorten notation, we omitted  dependence of $\widehat{u}_i,\widehat{u}_j$ on bubble and port coefficients (cf. \eqref{eq:ansatz_local}).
We introduce the HF quadrature rules 
$ \{ ( x_{\ell,q}^{\rm p}, \rho_{\ell,q}^{\rm p}  ) \}_{q=1}^{N_{\ell}^{\rm p}}$
on the archetype ports $\Gamma_{\ell}^{\rm a}$ for $\ell\in \mathcal{L}$; then, we have
\begin{subequations}
	\label{eq:tedious_formula}
	\begin{equation}
		\label{eq:tedious_formula_a} 
		\frac{1}{2}
		\sum_{i=1}^{N_{\rm dd}}
		\sum_{j \in {\rm Neigh}_i} 
		\int_{\Gamma_{i,j}} \| \widehat{u}_i( \boldsymbol{\alpha}_i, \boldsymbol{\beta}_i) - \widehat{u}_j( \boldsymbol{\alpha}_j, \boldsymbol{\beta}_j)   \|_2^2 \, dx
		\, \approx \,
		\frac{1}{2} 
		\sum_{i=1}^{N_{\rm dd}}
		\boldsymbol{\rho}_{\texttt{L}_i}^{\rm p} \cdot
		\boldsymbol{\eta}_{i}^{\rm p}
		(\boldsymbol{\alpha},\boldsymbol{\beta}   ) 
	\end{equation}			
	where $\boldsymbol{\alpha} = [\boldsymbol{\alpha}_1,\ldots,\boldsymbol{\alpha}_{N_{\rm dd}}]\in \mathbb{R}^N$ with $N:=n N_{\rm dd}$,  
	$\boldsymbol{\beta} = [\boldsymbol{\beta}_1,\ldots,\boldsymbol{\beta}_{N_{\rm dd}}]\in \mathbb{R}^M$, and
	\begin{equation}
		\label{eq:tedious_formula_b}	
		\boldsymbol{\eta}_{i}^{\rm p}
		(\boldsymbol{\alpha},\boldsymbol{\beta}   ) 
		= \left[
		\begin{array}{l}
			\eta_i^{\rm p}\left( \Phi_i(  x_{\ell,1}^{\rm p} );     
			\boldsymbol{\alpha},\boldsymbol{\beta} \right) \\
			\vdots \\
			\eta_i^{\rm p}\left( \Phi_i(  x_{\ell,N_{\ell}^{\rm p}}^{\rm p} );     
			\boldsymbol{\alpha},\boldsymbol{\beta} \right) \\
		\end{array}
		\right]
	\end{equation}
	with 
	\begin{equation}
		\label{eq:tedious_formula_c}	
		\eta_i^{\rm p}\left( x; 
		\boldsymbol{\alpha},\boldsymbol{\beta} \right)
		=
		\left(
		\sum_{j \in {\rm Neigh}_i: x \in \Omega_j} 
		\| \widehat{u}_i(x; \boldsymbol{\alpha}_i,\boldsymbol{\beta}_i ) \, - \,
		\widehat{u}_j(x; \boldsymbol{\alpha}_j,\boldsymbol{\beta}_j ) 
		\|_2^2 
		\right)
		J_i^{\rm bnd}(  \Phi_i^{-1}(x) ),
	\end{equation}
	for $i=1,\ldots,N_{\rm dd}$.
\end{subequations}

Evaluation of \eqref{eq:tedious_formula} is expensive due to   the need to integrate over the port boundaries
$\bigcup_{i=1}^{N_{\rm dd}}$  
$\bigcup_{j\in {\rm Neigh}_i}  \Gamma_{i,j}$: we should thus replace the HF quadrature vectors
$\{  \boldsymbol{\rho}_{\ell}^{\rm p}  \}_{\ell\in \mathcal{L}}$ with sparse EQ vectors
$\{  \boldsymbol{\rho}_{\ell}^{\rm p,eq}  \}_{\ell\in \mathcal{L}}$. In conclusion, we obtain the discrete OS2 formulation:  find
$\widehat{\boldsymbol{\beta}}  = [ \widehat{\boldsymbol{\beta}}_1  ,\ldots, \widehat{\boldsymbol{\beta}}_{N_{\rm dd}}  ] \in \mathbb{R}^{M}$ such that
\begin{subequations}
	\label{eq:OS2_eq}
	\begin{equation}
		\label{eq:OS2_eq_a}
		\widehat{\boldsymbol{\beta}} 
		\in {\rm arg}
		\min_{\boldsymbol{\beta}\in \mathbb{R}^M}
		\mathfrak{f}^{\rm eq} (  {\boldsymbol{\beta}}   )
		\,	= \,
		\mathfrak{F} \left(
		\widehat{\texttt{F}}^{\rm eq}( \boldsymbol{\beta}  ), \,
		\boldsymbol{\beta}, 
		\{    \boldsymbol{\rho}_{\ell}^{\rm p,eq}  \}_{\ell \in \mathcal{L}}
		\right)
	\end{equation}
	where 
	$\widehat{\texttt{F}}^{\rm eq}: \mathbb{R}^{M} \to \mathbb{R}^{N}$  is 
	the full port-to-bubble map  
	such that
	$\widehat{\texttt{F}}^{\rm eq}({\boldsymbol{\beta}} )
	=\left[
	\widehat{\texttt{F}}_1^{\rm eq}(\boldsymbol{\beta}_1 ),\ldots,
	\widehat{\texttt{F}}_{N_{\rm dd}}^{\rm eq}(\boldsymbol{\beta}_{N_{\rm dd}} )
	\right]^T$, and 
	\begin{equation}
		\label{eq:OS2_eq_b}
		\mathfrak{F} \left(
		\boldsymbol{\alpha} , \,
		\boldsymbol{\beta}, 
		\{    \boldsymbol{\rho}_{\ell}^{\rm p,eq}  \}_{\ell \in \mathcal{L}}
		\right)
		=
		\frac{1}{2}
		\;
		\sum_{i=1}^{N_{\rm dd}}
		\boldsymbol{\rho}_{\texttt{L}_i}^{\rm p,eq} \cdot
		\boldsymbol{\eta}_{i}^{\rm p}
		(\boldsymbol{\alpha},\boldsymbol{\beta}   ) .
	\end{equation}
	If we denote by $Q$ the total number of quadrature points with repetitions times the number of state variables $D$,
	\begin{equation}
		\label{eq:definitionQ}
		Q : =
		D \left(
		\sum_{i=1}^{N_{\rm dd}}
		\sum_{q=1}^{N_{\texttt{L}_i}^ {\rm p}}
		{\rm card} 
		\left\{
		j: \Phi_i(   x_{\texttt{L}_i,q}^{\rm p}   ) \in \Omega_j
		\right\}
		H(   \rho_{\texttt{L}_i,q}^{\rm p,eq}    )
		\right),
		\quad
		{\rm with} \; 
		H(x) = \left\{
		\begin{array}{ll}
			1 & {\rm if} \, x> 0 \\
			0 & {\rm otherwise}  \\
		\end{array}
		\right.
	\end{equation}
	we find that there exist $\mathbf{P}\in \mathbb{R}^{Q \times N}$ and 
	$\mathbf{Q}\in \mathbb{R}^{ Q \times M}$ such that
	\begin{equation}
		\label{eq:algebraic_objective_fun}
		\mathfrak{f}^{\rm eq}
		\left(  {\boldsymbol{\beta}}
		\right)
		=
		\frac{1}{2}
		\| \mathbf{r}^{\rm eq} \left(   \boldsymbol{\beta} \right)
		\|_2^2,
		\quad
		{\rm where} \;\; 
		\mathbf{r}^{\rm eq} \left(   {\boldsymbol{\beta}}
		\right)
		=
		\mathbf{P} \;  \widehat{\texttt{F}}^{\rm eq}({\boldsymbol{\beta}} )
		+ \mathbf{Q} {\boldsymbol{\beta}}.
	\end{equation}
\end{subequations}

\subsection{Discussion}

The remarks below provide a number of comments on the OS2 statement introduced in the previous section.

\begin{remark}
	\textbf{Algebraic representation of the local ROBs.}
	Exploiting  notation introduced at the end of section  \ref{sec:preliminaries}, the archetype bubble ROB $Z_{\ell}^{\rm a,b}:\mathbb{R}^{n} \rightarrow \mathcal{Z}_{\ell}^{\rm a,b}$ admits the algebraic representation 
	$Z_{\ell}^{\rm a,b}: 
	\boldsymbol{\alpha}\in \mathbb{R}^n 
	\mapsto  
	\big(
	\mathcal{T}_{\ell}^{\rm a}, 
	\mathbf{Z}_{\ell}^{\rm b}  \boldsymbol{\alpha}
	\big)$
	for some $\mathbf{Z}_{\ell}^{\rm b} \in \mathbb{R}^{N_{\ell}^{\rm a} \times n}$, while the deployed operators can be stated as $Z_i^{\rm b}:\boldsymbol{\alpha}\in \mathbb{R}^{n}\mapsto 
	\big(
	\Phi_i(\mathcal{T}_{\texttt{L}_i}^{\rm a} ), 
	\mathbf{Z}_{\texttt{L}_i}^{\rm b}  
	\boldsymbol{\alpha}\big)$, for $i=1,\ldots, N_{\rm{dd}}$. 
	Note that by virtue of the correspondence between archetype and deployed spaces, we do not have to explicitly instantiate --- and then store --- the bubble  ROBs for each configuration.
	The same applies for the port bases.  
\end{remark}

\begin{remark}
	\label{remark:dirichlet}
	\textbf{Extension to non-homogeneous Dirichlet  conditions.}
	The  OS2 formulation can readily deal with non-homogeneous Dirichlet boundary conditions. Towards this end, 
	for $i=1,\ldots,N_{\rm dd}$,
	given the Dirichlet datum $g_{i}^{\rm dir} : \Gamma_i^{\rm dir} \to \mathbb{R}^D$, 
	we introduce the lift $u_i^{\rm dir}$ such that
	$u_i^{\rm dir}|_{\Gamma_i^{\rm dir}   } = g_{i}^{\rm dir}$, and the ansatz
	$$
	\widehat{u}_i( \widehat{\boldsymbol{\alpha}}_i ,   
	\widehat{\boldsymbol{\beta}}_i )
	\,  =  \, u_{i}^{\rm dir} + Z_i^{\rm b} \widehat{\boldsymbol{\alpha}}_i +
	W_i^{\rm p} \widehat{\boldsymbol{\beta}}_i,
	\quad
	i=1,\ldots,N_{\rm dd}.
	$$
	Here,  $ \widehat{u}_i^{\rm b} = Z_i^{\rm b} \widehat{\boldsymbol{\alpha}}_i$ should approximate the bubble field $u|_{\Omega_i} - \texttt{E}_i ( (u-u_{i}^{\rm dir})|_{\Gamma_i} ) - u_i^{\rm dir}$, while $ \widehat{u}_i^{\rm p} = W_i^{\rm p} \widehat{\boldsymbol{\beta}}_i$ is an approximation of the (extended) port field $\texttt{E}_i ( (u-u_{i}^{\rm dir})|_{\Gamma_i} )$. Then, we can proceed as before to derive the reduced port-to-bubble maps and the low-dimensional OS2 formulation. We refer to 
	\cite{gunzburger2007reduced} for a thorough discussion on the imposition of Dirichlet boundary conditions in Galerkin ROMs.
\end{remark}


\begin{remark}
	\textbf{Computation of the matrices $\mathbf{P}, \mathbf{Q} $.}
	The matrices $\mathbf{P}, \mathbf{Q} $ depend on the configuration of interest but are independent of the port coefficients
	${\boldsymbol{\beta}}$: they can thus be defined after having instantiated the system and before solving the optimization problem. 
	Since the port quadrature points $\{ 
	\Phi_i(x_{\texttt{L}_i, q}^{\rm p})
	\}_{i,q}$ are configuration-dependent, we should resort to mesh interpolation to assemble the matrices $\mathbf{P}$ and $\mathbf{Q}$.
	In this work, we rely on structured meshes in the archetype components that enable logarithmic-in-$N_{\ell}^{\rm v}$ FE interpolations. 
\end{remark}

\begin{remark}
	\textbf{Hyper-reduction.}
	As required in CB-pMOR, hyper-reduction should be
	defined at the component level and is then translated  to the deployed system  using the mappings $\{\Phi_i\}_i$. From an algorithmic standpoint, an archetype component 
	$\ell \in \mathcal{L}$ 	
	should be interpreted as a complex data structure that comprises 
	(i) bubble and port ROBs; 
	(ii) the approximate residual $\widehat{\mathbf{R}}_{\ell}^{\rm eq}$ that enables effective computations of port-to-bubble maps;
	(iii)  the port quadrature rule 
	$\boldsymbol{\rho}_{\ell}^{ \rm p,eq}$
	associated with the approximate objective function \eqref{eq:OS2_eq_a}; and
	(iv)  a (structured) mesh structure for which efficient 
	(i.e., logarithmic-in-$N_{\ell}^{\rm v}$) interpolation procedures are available 
	for the computation of the matrices $\mathbf{P}, \mathbf{Q}$.
\end{remark}

\subsection{Solution to the OS2 minimization problem}
\label{sec:minimization} 
In view of the description of the numerical solution  to \eqref{eq:OS2_eq}, we observe that the Jacobian  of the global port-to-bubble map
$\widehat{\texttt{F}}^{\rm eq}:\mathbb{R}^M \to \mathbb{R}^N$ is block-diagonal (cf. \eqref{eq:port2bubble_reduced_gradient}):
\begin{subequations}
	\label{eq:objective_function_gradient}
	\begin{equation}
		\widehat{\mathbf{J}}_{\texttt{F}}^{\rm eq} 
		( {\boldsymbol{\beta}}  )
		=
		{\rm diag} 
		\left[
		\widehat{\mathbf{J}}_{\texttt{F}_1}^{\rm eq}
		(\boldsymbol{\beta}_1),
		\ldots,
		\widehat{\mathbf{J}}_{\texttt{F}_{N_{\rm dd}}}^{\rm eq} 
		( \boldsymbol{\beta}_{N_{\rm dd}} )
		\right],
		\;
		\widehat{\mathbf{J}}_{\texttt{F}_i}^{\rm eq} 
		( \boldsymbol{\beta}_i  ) : =
		\, -\,
		\left(
		\partial_{ \boldsymbol{\alpha}_i   }
		\widehat{\mathbf{R}}_i^{\rm eq} \right)^{-1}
		\partial_{ \boldsymbol{\beta}_i   }
		\widehat{\mathbf{R}}_i^{\rm eq}
		\Big|_{( \boldsymbol{\alpha}_i, \boldsymbol{\beta}_i ) \, = \, (  \widehat{\texttt{F}}_i^{\rm eq}
			\left(
			{\boldsymbol{\beta}}_i \right),  {\boldsymbol{\beta}}_i)}. 
	\end{equation}
	Then, we observe that
	\begin{equation}
		\nabla 
		\mathbf{r}^{\rm eq} \, = \,
		\mathbf{P}   \widehat{\mathbf{J}}_{\texttt{F}}^{\rm eq} 
		+ \mathbf{Q},
		\quad
		\nabla \,
		\mathfrak{f}^{\rm eq} \, = \,
		\left(
		\mathbf{P}   \widehat{\mathbf{J}}_{\texttt{F}}^{\rm eq} 
		+ \mathbf{Q} \right)^T \, \mathbf{r}^{\rm eq}.
		\label{eq:eq:algebraic_grad_objective_fun}
	\end{equation}
\end{subequations}
If $N_{\rm dd} \gg m$ (as in the cases considered  e.g. in \cite{HuKnPa13,mcbane2021component}), the Jacobian
$\widehat{\mathbf{J}}_{\texttt{F}}^{\rm eq}$ is highly sparse; note that explicit assembly of the local Jacobians requires to solve $m$ linear systems of size $n$, while matrix-vector  multiplications 
$\widehat{\mathbf{J}}_{\texttt{F}}^{\rm eq} \mathbf{v}$ and 
$\mathbf{v}^T \widehat{\mathbf{J}}_{\texttt{F}}^{\rm eq}$ require $N_{\rm dd}$ $n\times m$ matrix-vector multiplications and 
$N_{\rm dd}$ linear solves of size $m$.

The nonlinear least-square problem \eqref{eq:OS2_eq} can be solved using (i) steepest-descent or quasi-Newton methods, or (ii) Gauss-Newton or Levenberg-Marquandt algorithms, \cite{nocedal2006numerical}.
\begin{itemize}
	\item[(i)]
	Steepest-descent or quasi Newton methods
	only require the explicit calculation  of the objective function 
	$\mathfrak{f}^{\rm eq}$ and its gradient  
	$\nabla \, 
	\mathfrak{f}^{\rm eq}$, which can be computed without explicitly forming $\widehat{\mathbf{J}}_{\texttt{F}}^{\rm eq} $. However, these methods do not exploit the underlying least-square  structure of the optimization problem and might thus exhibit slower convergence  and/or  might be more prone to divergent behaviors.
	\item[(ii)]
	The Gauss-Newton method (GNM) reads as 
	$$
	\widehat{\boldsymbol{\beta}}^{(k+1)}
	=
	\widehat{\boldsymbol{\beta}}^{(k)}
	\, -\,
	\left(
	\nabla 
	\mathbf{r}^{\rm eq}
	\left(
	\widehat{\boldsymbol{\beta}}^{(k)}
	\right)
	\right)^{\dagger}
	\mathbf{r}^{\rm eq}
	\left(
	\widehat{\boldsymbol{\beta}}^{(k)}
	\right)
	$$
	where $( \cdot )^{\dagger}$ denotes the Moore-Penrose pseudo-inverse. The Levenberg-Marquandt algorithm (LMA) is a generalization of GNM that is typically more robust for poor choices of the initial condition. Note that GNM/LMA are the methods of choice for least-squares problems; however, they require the assembly of $\widehat{\mathbf{J}}_{\texttt{F}}^{\rm eq} $ at each iteration.
\end{itemize}

\begin{algorithm}[H]                      
	\caption{Solution to \eqref{eq:OS2_eq} through the Gauss-Newton method.}     
	\label{alg:objective}     
	
	\small
	\begin{flushleft}
		\emph{Inputs:}  
		$\boldsymbol{\alpha}^{(0)} = [\boldsymbol{\alpha}_1^{(0)} , \ldots,\boldsymbol{\alpha}_{N_{\rm dd}}^{(0)} ]$,
		$\boldsymbol{\beta}^{(0)}  = [\boldsymbol{\beta}_1^{(0)} , \ldots,\boldsymbol{\beta}_{N_{\rm dd}}^{(0)}]$ initial conditions (cf. Eq. \eqref{eq:IC_GNM}),
		$tol>0, \texttt{maxit}$.
		\smallskip
		
		\emph{Outputs:} 
		$\widehat{\boldsymbol{\beta}}$ port coefficients,
		$\widehat{\boldsymbol{\alpha}}= \widehat{\texttt{F}}^{\rm eq}(  \widehat{\boldsymbol{\beta}}  )$  bubble coefficients.
	\end{flushleft}                      
	
	\normalsize 
	
	\begin{algorithmic}[1]
		\State
		Compute the matrices
		$\mathbf{P},\mathbf{Q}$ in  \eqref{eq:algebraic_objective_fun}.
		\smallskip

		\State Set $\widehat{\boldsymbol{\beta}}^{(0)} = \boldsymbol{\beta}^{(0)}$ and 
		$\widehat{\boldsymbol{\alpha}}=\boldsymbol{\alpha}^{(0)}$.
		\medskip
		
		\For {$k=1, \ldots, \texttt{maxit}$ }
		
		\For {$i=1, \ldots, N_{\rm dd}$ }
		
		\State
		Compute $\boldsymbol{\alpha}_i$ s.t. 
		$\widehat{\mathbf{R}}_i^{\rm eq}( \boldsymbol{\alpha}_i, {\boldsymbol{\beta}}_i^{(k)}) = \mathbf{0}$ using Newton's method with initial condition 
		$\widehat{\boldsymbol{\alpha}}_i$.
		\vspace{3pt}
		
		\State
		Compute 
		$\widehat{\mathbf{J}}_{\texttt{F}_i}^{\rm eq} 
		(\boldsymbol{\beta}_i^{(k)})$
		(cf. \eqref{eq:objective_function_gradient}).
		\vspace{3pt}
		
		\EndFor
		\vspace{3pt}
		
		\State
		Update  $\widehat{\boldsymbol{\alpha}} =[\boldsymbol{\alpha}_1,\ldots,\boldsymbol{\alpha}_{N_{\rm dd}}]$.
		\vspace{3pt}
		
		\State
		Compute $\mathbf{r}^{\rm eq, (k)} = \mathbf{P} \widehat{\boldsymbol{\alpha}} + \mathbf{Q} \widehat{\boldsymbol{\beta}}_i^{(k)}$ and
		$\nabla
		\mathbf{r}^{\rm eq,(k)} \, = \,
		\mathbf{P}   \widehat{\mathbf{J}}_{\texttt{F}}^{\rm eq} 
		+ \mathbf{Q}$.
		\vspace{3pt}
		
		\State
		Compute 
		$\widehat{\boldsymbol{\beta}}^{(k+1)}
		=
		\widehat{\boldsymbol{\beta}}^{(k)}
		\, -\,
		\left(
		\nabla 
		\mathbf{r}^{\rm eq,(k)}
		\right)^{\dagger}
		\mathbf{r}^{\rm eq,(k)}
		$
		\vspace{3pt}
		
		\If {$\| \widehat{\boldsymbol{\beta}}^{(k+1)} - \widehat{\boldsymbol{\beta}}^{(k)}   \|_2 < tol \| \widehat{\boldsymbol{\beta}}^{(k)} \|_2$}, \texttt{BREAK}
		\label{alg:line:exit}
		\EndIf
		
		\EndFor
		
		\State
		\label{alg:line:optimality}
		Return 
		$\widehat{\boldsymbol{\beta}} = \widehat{\boldsymbol{\beta}}^{(k+1)}$  and 
		$\widehat{\boldsymbol{\alpha}} = \widehat{\texttt{F}}^{\rm eq}(  \widehat{\boldsymbol{\beta}} )$. 
		
	\end{algorithmic}
\end{algorithm}

Algorithm \ref{alg:objective} summarizes the overall procedure as implemented in our code, which relies on GNM to solve \eqref{eq:OS2_eq}; we envision that our approach can cope with LMA with only minor changes: we omit the details.
Note that we update at each iteration the estimates of the bubble coefficients: this is important to speed up the solution to the local Newton problems.
In addition, the algorithm requires to provide an initial guess for port and bubble coefficients; we discuss the choice of the initial condition in section \ref{sec:methods} (cf. Eq. \eqref{eq:IC_GNM}).

As explained in \cite{madsen2004nlls}, for nonlinear least-squares problems of the form \eqref{eq:algebraic_objective_fun}, Gauss-Newton's method shows   quadratic convergence if $\mathbf{r}^{\rm eq}\big(\widehat{\boldsymbol{\beta}}  \big)=\mathbf{0}$  and a super-linear convergence if $\|
\mathbf{r}^{\rm eq}\big(\widehat{\boldsymbol{\beta}}  \big)
\|_2$ is small.
In the numerical results, we also investigate  performance of a quasi-Newton method  --- the limited-memory 
BFGS method  \cite{nocedal2006numerical}. Note that the implementation of the latter  follows a similar  procedure 
as in \Cref{alg:objective} with only minor changes: we omit the details.

\begin{remark}
	We remark that the internal loop at  lines 4-7 in \Cref{alg:objective} and the construction of the matrices $\mathbf{P},\mathbf{Q}$ are embarrassingly parallelizable. 
\end{remark}

\section{Methodology}
\label{sec:methods}
\subsection{Data compression}
\label{sec:data_compression}
In this work, we resort to global solves to construct the archetype ROBs $\{  (Z_{\ell}^{\rm a,b},   W_{\ell}^{\rm a,b}) \}_{\ell \in \mathcal{L}}$,
$Z_{\ell}^{\rm a,b} = [\zeta_{\ell,1}^{\rm a,b},\ldots, \zeta_{\ell,n}^{\rm a,b}  ]$,
$W_{\ell}^{\rm a,p} = [\psi_{\ell,1}^{\rm a,p},\ldots, \psi_{\ell,n}^{\rm a,p}]$. 
We generate $n_{\rm train}$ global configurations $\{  \mu^{(k)}  \}_{k=1}^{n_{\rm train}}$  and we denote by  $\{ \left( \Omega_i^{(k)}, \texttt{L}_i^{(k)} \right) \}_{i,k}$ the corresponding labeled partitions;
we estimate the global solutions 
$\{  u^{(k)} \}_{k=1}^{n_{\rm train}}$
using a standard FE solver and we assemble the datasets
\begin{subequations}
	\begin{equation}
		\label{eq:data_compression_a}
		\mathcal{D}_{\ell} =
		\left\{
		u^{(k)}  |_{\Omega_i^{(k)}} \circ \Phi_i^{(k)}  \, : \,
		\texttt{L}_i^{(k)} = \ell, \;
		k=1,\ldots,n_{\rm train} 
		\right\} \subset \mathcal{X}_{\ell}^{\rm a},
		\quad
		\ell \in \mathcal{L};
	\end{equation}
	we further define the bubble and port datasets
	\begin{equation}
		\label{eq:data_compression_b}
		\mathcal{D}_{\ell}^{\rm b}
		:=
		\left\{
		w - \texttt{E}_{\ell}^{\rm a} (w|_{\Gamma_{\ell}^{\rm a}})  \, : \,
		w \in \mathcal{D}_{\ell}
		\right\},
		\quad
		\mathcal{D}_{\ell}^{\rm p}
		:=
		\left\{
		\texttt{E}_{\ell}^{\rm a} (w|_{\Gamma_{\ell}^{\rm a}})  \, : \,
		w \in \mathcal{D}_{\ell}
		\right\};
	\end{equation}
\end{subequations}
finally, we apply proper orthogonal decomposition (POD, \cite{volkwein2011model}) based on the method of snapshots \cite{Sir87} with inner product 
$(\cdot, \cdot)_{\ell}$, 
to obtain the local approximation spaces. 
\Cref{alg:data_compression} summarizes the computational procedure.

In view of the application of the empirical quadrature procedures described in sections \ref{sec:hyper_port2bubble} and \ref{sec:EQ_bnd}, for all
$\ell \in \mathcal{L}$
we further  compute the projected coefficients
$\{ \boldsymbol{\alpha}_{\ell,j}   \}_{j=1}^{n_{\rm train,\ell}}$,
$\{ \boldsymbol{\beta}_{\ell,j}   \}_{j=1}^{n_{\rm train,\ell}}$ 
\begin{equation}
	\label{eq:projected_coefficients}
	\left( \boldsymbol{\alpha}_{\ell,j}  \right)_i
	=
	\left(u_{\ell,j}^{\rm b}, \zeta_{\ell,i}^{\rm a,b} \right)_{\ell},
	\quad
	\left( \boldsymbol{\beta}_{\ell,j} \right)_q
	=
	\left(u_{\ell,j}^{\rm p}, \psi_{\ell,q}^{\rm a,p} \right)_{\ell},
	\qquad
	\ell \in \mathcal{L},
\end{equation}
for $i=1,\ldots,n$,  $q=1,\ldots,m$, 
$j=1,\ldots,n_{\rm train, \ell}$, 
where $u_{\ell,j}^{\rm b}$ (resp., $u_{\ell,j}^{\rm p}$) denotes the $j$-th  bubble (resp., port) solution in the dataset 
$\mathcal{D}_{\ell}^{\rm b}$ (resp., $\mathcal{D}_{\ell}^{\rm p}$).

\begin{algorithm}[H]
	\caption{Data compression based on global solves}
	\label{alg:data_compression}
	\small
	\begin{flushleft}
		\emph{Inputs:}  training parameters $\{ \mu^{(k)}\}_{k=1}^{n_{\rm train}}$; $m,n$ ROB dimensions.  \\
		\smallskip
		\emph{Outputs:} 
		$\{  (Z_{\ell}^{\rm a,b},   W_{\ell}^{\rm a,b}) \}_{\ell \in \mathcal{L}}$ ROBs;
		$\{ \boldsymbol{\alpha}_{\ell}^{(k)}    \}_{k=1}^{n_{\rm train,\ell}}$,
		$\{ \boldsymbol{\beta}_{\ell}^{(k)}    \}_{k=1}^{n_{\rm train,\ell}}$ local optimal coefficients.
	\end{flushleft}

	\normalsize 
	\begin{algorithmic}[1]
		\State
		Initialize $\mathcal{D}_{\ell}^{\rm b} =\mathcal{D}_{\ell}^{\rm p} =  \emptyset$ for $\ell \in \mathcal{L}$.
		\smallskip

		\For {$k=1, \ldots, n_{\rm train}$}
		\State 
		Estimate the global solution $u_{\mu}$  to  \eqref{eq:PDE_model} using a global FE method.
		\smallskip
		
		\State
		Update the datasets  $\mathcal{D}_{\ell}^{\rm b}$ and 
		$\mathcal{D}_{\ell}^{\rm p}$ using \eqref{eq:data_compression_b}.
		\EndFor
		\medskip	
		
		\State 
		Perform POD to obtain the ROBs
		$Z_{\ell}^{\rm a,b} = [\zeta_{\ell,1}^{\rm a,b},\ldots, \zeta_{\ell,n}^{\rm a,b}  ]$ and
		$W_{\ell}^{\rm a,p} = [\psi_{\ell,1}^{\rm a,b},\ldots, \psi_{\ell,n}^{\rm a,b}]$
		\smallskip
		
		\State 
		Define the optimal coefficients 
		$\{ \boldsymbol{\alpha}_{\ell}^{(k)}    \}_{k=1}^{n_{\rm train,\ell}}$,
		$\{ \boldsymbol{\beta}_{\ell}^{(k)}    \}_{k=1}^{n_{\rm train,\ell}}$  using \eqref{eq:projected_coefficients}.
	\end{algorithmic}
\end{algorithm}

We remark that the proposed approach
---
which was previously considered in \cite{pegolotti2021model}
---
might be  highly inefficient since it requires global solves that are often unfeasible in the framework of CB-pMOR. 
We envision to further extend the localized training approach in \cite{smetana2022localized} to address this issue.
For practical applications, we envision that global solves should be performed using a standard FE solver and then resorting to FE interpolation routines to extract the local solutions: this procedure inevitably introduces an error at the scale of the FE mesh size between full-order and reduced-order models. 
Even if this error might be negligible for applications, it   hinders the interpretations  of the  numerical investigations. 
To avoid this issue, in the numerical experiments, we  rely on  the HF model  \eqref{eq:hybridizedOS2} to generate the dataset of local solutions.

\subsection{Hyper-reduction of port-to-bubble problems}
\label{sec:hyper_port2bubble}

We here rely on element-wise EQ, that is we replace the residuals \eqref{eq:local_model} in \eqref{eq:local_problems} with the weighted residual associated  with the variational form
\begin{equation} 
	\label{eq:weighted_residual}
	\mathcal{G}_{\ell}^{\rm a, eq} (w, v;  \mu_{\ell} )
	=\sum_{k=1}^{N_{\ell}^{\rm e}}
	\;
	\rho_{\ell,k}^{\rm eq} 
	\left(
	\int_{\texttt{D}_{\ell,k}} \;
	\eta_{\ell}^{\rm a,e}  (w,v; \mu_{\ell}) \, dx
	\; + \;
	\int_{\partial \texttt{D}_{\ell,k}} \;
	\eta_{\ell}^{\rm a,f} (w,v; \mu_{\ell})
	\, dx
	\right),
\end{equation}
where $\boldsymbol{\rho}_{\ell}^{\rm eq} = 
[\rho_{\ell,1}^{\rm eq}, \ldots, \rho_{\ell, N_{\ell}^{\rm e}}^{\rm eq}]^T$ is a sparse vector of non-negative weights. 

This hyper-reduction   approach,
which has been considered in a number of previous works including \cite{iollo2022adaptive},
is   discussed for completeness in \ref{sec:app_hyper_port2bubble}.
We anticipate that the algorithm takes as input the projected coefficients
\eqref{eq:projected_coefficients}
generated by \Cref{alg:data_compression}
and the associated local parameters,
$\{ (\boldsymbol{\alpha}_{\ell}^{(j)} , 
\boldsymbol{\beta}_{\ell}^{(j)}   ,
\mu_{\ell}^{(j)}  \}_{j=1}^{n_{\rm train,\ell}}$.

We remark that, as discussed in \cite{du2022efficient}, the use of elementwise-
(as opposed to pointwise-)
reduced quadrature formulations leads to significantly less efficient ROMs, particularly for high-order FE discretizations. On the other hand, elementwise reduced quadrature formulations are significantly easier to implement and can easily cope with geometry deformations \cite{taddei2021discretize}. We refer to 
\cite{farhat2021computational,yano2019discontinuous,yano2021model} for a thorough introduction to state-of-the-art hyper-reduction techniques.

\subsection{Hyper-reduction of the objective function}
\label{sec:hyper_bnd}

Exploiting \eqref{eq:OS2_eq_b},
it is easy to verify that --- we here stress dependence on the parameter value $\mu$ --- 
$$
\mathfrak{F} \left(
\boldsymbol{\alpha} , \,
\boldsymbol{\beta}, 
\{    \boldsymbol{\rho}_{\ell}^{\rm p,eq}  \}_{\ell \in \mathcal{L}}, \mu
\right)
=
\frac{1}{2}
\;
\sum_{\ell \in \mathcal{L}} 
\left(
\sum_{i: \texttt{L}_i=\ell}
\boldsymbol{\eta}_{i}^{\rm p}
(\boldsymbol{\alpha},\boldsymbol{\beta}  , \mu )
\right)
\cdot 
\boldsymbol{\rho}_{\ell}^{\rm p,eq}  
\, = \,
\frac{1}{2}
\sum_{\ell \in \mathcal{L}} 
\sum_{j=1}^{N_{\rm dd, \ell}}
\left(
\mathbf{G}_{\ell}^{\rm p}
(\boldsymbol{\alpha},\boldsymbol{\beta} , \mu  )
\boldsymbol{\rho}_{\ell}^{\rm p,eq}  
\right)_j,
$$
where $N_{\rm dd, \ell}$ is the number of components of type $\ell$ and $\{\mathbf{G}_{\ell}^{\rm p}\}_{\ell}$ are suitable matrices;  to provide a concrete example, for the model problem of section \ref{sec:model_problem}, we have
$$
\mathbf{G}_{\rm int}^{\rm p}(\boldsymbol{\alpha},\boldsymbol{\beta} , \mu  )
=
\left[
\begin{array}{c}
	\left(  \boldsymbol{\eta}_1^{\rm p}
	(\boldsymbol{\alpha},\boldsymbol{\beta} , \mu   ) \right)^T
	\\
	\vdots \\
	\left(  \boldsymbol{\eta}_{Q_{\rm a}}^{\rm p}
	(\boldsymbol{\alpha},\boldsymbol{\beta} , \mu   ) \right)^T
	\\
\end{array}
\right],
\quad
\mathbf{G}_{\rm ext}^{\rm p}(\boldsymbol{\alpha},\boldsymbol{\beta} , \mu  )
=
\left(  \boldsymbol{\eta}_{Q_{\rm a}+1}^{\rm p}
(\boldsymbol{\alpha},\boldsymbol{\beta} , \mu   ) \right)^T.
$$
In order to speed up the evaluation of $\mathfrak{f}^{\rm eq}$, it is necessary to build a sparse quadrature rule 
$\{ \boldsymbol{\rho}_{\ell}^{\rm p,eq}  \}_{\ell \in \mathcal{L}}$. In the remainder of this section, we propose two different strategies to address this task: the former relies on the solution to  a suitable sparse representation problem and is tightly linked to the EQ procedure employed for hyper-reduction of the port-to-bubble maps; the latter relies on a variant of the empirical interpolation method (EIM, \cite{BMNP04}) for vector-valued functions.

\subsubsection{Empirical quadrature method}
\label{sec:EQ_bnd}

We denote by $(\boldsymbol{\alpha}^{(k)}, \boldsymbol{\beta}^{(k)} )$ the projected bubble and port coefficients associated with the $k$-th configuration $\mu^{(k)}$ and Eq. \eqref{eq:projected_coefficients}; we further denote by 
$(\boldsymbol{\alpha}_0^{(k)}, \boldsymbol{\beta}_0^{(k)} )$ the bubble and port coefficients associated with the sample means,
\begin{subequations}
	\label{eq:IC_GNM}
	\begin{equation}
		\boldsymbol{\alpha}_0^{(k)} = \left[
		\begin{array}{c}
			\boldsymbol{\alpha}_{0,1}^{(k)} \\
			\vdots \\
			\boldsymbol{\alpha}_{0,N_{\rm dd,(k)}}^{(k)}  \\
		\end{array}
		\right],
		\;\; 
		\boldsymbol{\beta}_0^{(k)} = \left[
		\begin{array}{c}
			\boldsymbol{\beta}_{0,1}^{(k)} \\
			\vdots \\
			\boldsymbol{\beta}_{0,N_{\rm dd,(k)}}^{(k)}  \\
		\end{array}
		\right], 
		\quad
		k=1,\ldots,n_{\rm train},
	\end{equation}
	where  $ \boldsymbol{\alpha}_{0,i}^{(k)} =  \boldsymbol{\alpha}_{\texttt{L}_i^{(k)}}^{\rm avg}$ and
	$\boldsymbol{\beta}_{0,i}^{(k)} =  \boldsymbol{\beta}_{\texttt{L}_i^{(k)}}^{\rm avg}$, with 
	\begin{equation}
		\boldsymbol{\alpha}_{\ell}^{\rm avg} :=
		\frac{1}{n_{\rm train,\ell}} \sum_{j=1}^{  n_{\rm train,\ell}  }  \boldsymbol{\alpha}_{\ell,j},
		\quad
		\boldsymbol{\beta}_{\ell}^{\rm avg} :=
		\frac{1}{n_{\rm train,\ell}} \sum_{j=1}^{  n_{\rm train,\ell}  }  \boldsymbol{\beta}_{\ell,j},
		\qquad \forall \, \ell \in \mathcal{L}.
	\end{equation}
	We anticipate that \eqref{eq:IC_GNM} is used in the numerical results to initialize the Gauss-Newton's algorithm. 
\end{subequations}

Given the random samples $s^{(k)} \overset{\rm iid}{\sim} {\rm Uniform}(0,1)$, we define the matrices
\begin{subequations}
	\label{eq:NNLS_bnd}
	\begin{equation}
		\label{eq:NNLS_bnd_a}
		\mathbf{C}_{\ell}
		=
		\left[
		\begin{array}{l}
			\mathbf{G}_{\ell}^{\rm p}
			(\widetilde{\boldsymbol{\alpha}}^{(1)},\widetilde{\boldsymbol{\beta}}^{(1)} , \mu^{(1)}  )
			\\[2mm]
			\vdots \\
			\mathbf{G}_{\ell}^{\rm p}
			(\widetilde{\boldsymbol{\alpha}}^{(n_{\rm train})},\widetilde{\boldsymbol{\beta}}^{(n_{\rm train})} , \mu^{(n_{\rm train})}  )
			\\[2mm]
			\mathbf{1}^T
		\end{array}
		\right],
		\quad
		\forall \, \ell \in \mathcal{L},
	\end{equation}
	where $\mathbf{1}$ is the vector with entries all equal to one, and $\widetilde{\boldsymbol{\alpha}}^{(k)} $ and 
	$\widetilde{\boldsymbol{\beta}}^{(k)}$ are random convex  interpolations between the projected bubble and port coefficients $(\boldsymbol{\alpha}^{(k)}, \boldsymbol{\beta}^{(k)} )$  and 
	the initial conditions  for the GNM 
	$(\boldsymbol{\alpha}_0^{(k)}, \boldsymbol{\beta}_0^{(k)} )$,
	\begin{equation}
		\label{eq:NNLS_bnd_b}
		\widetilde{\boldsymbol{\alpha}}^{(k)} =
		(1 - s^{(k)}) {\boldsymbol{\alpha}}^{(k)} 
		+
		s^{(k)} {\boldsymbol{\alpha}}_0^{(k)} ,
		\quad
		\widetilde{\boldsymbol{\beta}}^{(k)} =
		(1 - s^{(k)}) {\boldsymbol{\beta}}^{(k)} 
		+
		s^{(k)} {\boldsymbol{\beta}}_0^{(k)} ,
		\quad
		k=1,\ldots,n_{\rm train}.
	\end{equation}
	The first $n_{\rm train}$  blocks of $\mathbf{C}_{\ell}$ are associated to the ``manifold accuracy constraints'', while the last row is associated to the ``constant accuracy constraint'' \cite{yano2019lp}.
	Then, we compute the empirical weights 
	$\{  \boldsymbol{\rho}_{\ell}^{\rm p,eq} \}_{\ell\in \mathcal{L}}$ by approximately solving the non-negative least-square problem
	\begin{equation}
		\label{eq:NNLS_bnd_c}
		\min_{ \boldsymbol{\rho} \in \mathbb{R}^{N_{\ell}^{\rm p}}   } \| \mathbf{C}_{\ell} \left(  \boldsymbol{\rho} - 
		\boldsymbol{\rho}_{\ell}^{\rm p} \right)   \|_2, 
		\quad
		{\rm s.t.} \;\; \boldsymbol{\rho} \geq 0
	\end{equation}
	up to a tolerance $tol_{\rm eq}^{\rm obj}$ using the Matlab function \texttt{lsqnonneg}, which implements the iterative procedure proposed in \cite{lawson1974solving}.
\end{subequations}

The choice of the port and bubble coefficients
$\{(\widetilde{\boldsymbol{\alpha}}^{(k)},\widetilde{\boldsymbol{\beta}}^{(k)} )\}_k$
for the ``accuracy constraints''
in \eqref{eq:NNLS_bnd_a}
is justified by the fact that the objective function should be accurate for all 
port and bubble coefficients considered during the GNM iterations; this choice  is found to empirically improve the conditioning of the non-negative least-square problem and ultimately improve performance --- compared to the choice $\widetilde{\boldsymbol{\alpha}}^{(k)} =
{\boldsymbol{\alpha}}^{(k)} $, $\widetilde{\boldsymbol{\beta}}^{(k)} =
{\boldsymbol{\beta}}^{(k)}$.
The constant function accuracy constraint, which was first proposed in \cite{yano2019lp} for hyper-reduction of monolithic ROMs, is important to bound the $\ell^1$ norm of the empirical weights; we have indeed
\begin{equation}
	\label{eq:importance_constant_accuracy}
	\|  \boldsymbol{\rho}_{\ell}^{\rm p,eq} \|_1
	\leq
	| \mathbf{1} \cdot (\boldsymbol{\rho}_{\ell}^{\rm p,eq} -  \boldsymbol{\rho}_{\ell}^{\rm p} ) | \, + \,
	\|  \boldsymbol{\rho}_{\ell}^{\rm p} \|_1
	\leq
	\| \mathbf{C}_{\ell} \left(  \boldsymbol{\rho}_{\ell}^{\rm p,eq} -  \boldsymbol{\rho}_{\ell}^{\rm p} \right)  \|_2
	+
	\|  \boldsymbol{\rho}_{\ell}^{\rm p} \|_1,
	\quad
	\forall \, \ell \in \mathcal{L}.
\end{equation}
We also observe that, even if hyper-reduction is ultimately performed at the local level,  for each archetype component, the EQ procedure requires global solves to  define the matrices $\{  \mathbf{C}_{\ell}  \}_{\ell}$.

\subsubsection{Empirical interpolation method}
\label{sec:EIM_bnd}

The objective function $\mathfrak{F}$ is designed to penalize the jump of the solution at the components' interface. Since the jumps are dictated by the behavior of the port modes $\{ \psi_{\ell,i}^{\rm a,p}  \}_{i=1}^m$ on the ports $\Gamma_{\ell}$, we propose to replace the integral in \eqref{eq:tedious_formula} with the discrete sum
\begin{equation}
	\label{eq:EIM_approximation}
	\frac{1}{2}
	\sum_{i=1}^{N_{\rm dd}}
	\sum_{j \in {\rm Neigh}_i} 
	\int_{\Gamma_{i,j}} \| \widehat{u}_i( \boldsymbol{\alpha}_i, \boldsymbol{\beta}_i) - \widehat{u}_j( \boldsymbol{\alpha}_j, \boldsymbol{\beta}_j)   \|_2^2 \, dx
	\, \approx \,
	\frac{1}{2}
	\sum_{q\in \texttt{I}_{\ell}^{\rm p,eq}}
	\left(
	\boldsymbol{\eta}_i^{\rm p}( \boldsymbol{\alpha}, \boldsymbol{\beta}, \mu  )
	\right)_q,
\end{equation}
where $\texttt{I}_{\ell}^{\rm p,eq} \subset \{1, \ldots,N_{\ell}^{\rm p} \}$ are chosen so that we can adequately recover any element of $\mathcal{Z}_{\ell}^{\rm a,p}$ based on the information at the points
$\{ x_{\ell,j}^{\rm p}  \}_{j\in \texttt{I}_{\ell}^{\rm p,eq}}$. Note that the approximation \eqref{eq:EIM_approximation} is an inconsistent approximation of the $L^2$ integral \eqref{eq:tedious_formula}; however, we  expect --- and we verify numerically ---  that the minimization of the right-hand side of \eqref{eq:EIM_approximation}
should control the jump at elements' interfaces and ultimately ensure accurate performance. 

We here rely on a variant of EIM to select the quadrature indices 
$\texttt{I}_{\ell}^{\rm p,eq}$. EIM was first proposed in \cite{BMNP04} to identify accurate interpolation points for arbitrary sets of scalar functions. In this work, we resort to the extension of EIM to vector-valued fields considered in \cite{taddei2019offline}. We refer to the MOR literature for other variants of EIM for vector-valued fields; in particular, we observe that the present algorithm returns exactly $m$ quadrature points: we refer to 
\cite[Algorithm 2]{maday2015parameterized} and to \cite{maday2021overcolloc}
for extensions of EIM that resort to over-collocation to improve performance. 

Algorithm \ref{alg:EIM} reviews the computational procedure: 
note that,
for each $\ell\in \mathcal{L}$, the algorithm takes as input the port functions 
$\{ \psi_{\ell,i}^{\rm a,p}  \}_{i=1}^m$ and returns the indices 
$\texttt{I}_{\ell}^{\rm p,eq}$.
Given the set of indices $\texttt{I}_{\ell}^{\rm p,eq}$ and the space
$\mathcal{Z}_{\ell}^{\rm a,p}$, we denote by $\mathcal{I}_{\ell,m}$ the approximation least-square operator
$$
\mathcal{I}_{\ell,m}(v):=
\mathcal{I}
\left( v;  \texttt{I}_{\ell}^{\rm p,eq},   \mathcal{Z}_{\ell}^{\rm a,p} \right)
= {\rm arg} \min_{\psi \in    \mathcal{Z}_{\ell}^{\rm a,p}} \;
\sum_{j\in   \texttt{I}_{\ell}^{\rm p,eq} }
\|  v(   x_{\ell,j}^{\rm p}  )   - \psi( x_{\ell,j}^{\rm p})   \|_2^2,
\quad
\forall \, v\in C(\Gamma_{\ell}; \mathbb{R}^D), \;\ell \in \mathcal{L}.
$$
Note  that for $D>1$  
$ \mathcal{I}_{\ell,m}$ is not an interpolation operator.

\begin{algorithm}[H]
	\caption{Empirical Interpolation Method for vector-valued fields}
	\label{alg:EIM}
	\hspace*{\algorithmicindent} \textbf{Input}:
	$\{ \psi_{\ell,i}^{\rm a,p}  \}_{i=1}^m$, $\ell \in \mathcal{L}$\\[2mm]
	\hspace*{\algorithmicindent} \textbf{Output}:
	$\texttt{I}_{\ell}^{\rm p,eq} = \{\texttt{i}_{\ell,1}^{\star}, \ldots,
	\texttt{i}_{\ell,m}^{\star} \}$ \\[-2mm]
	
	\begin{algorithmic}
		\State 
		Set $\displaystyle{\texttt{i}_{\ell,1}^{\star}:=\arg \max_{j  \in \{1,\ldots,N_{\ell}^{\rm p} \} } \|  \psi_{\ell,1}^{\rm a,p}(  x_{\ell,j}^{\rm p}   )  \|_2}$,
		and  define
		$ \mathcal{I}_{\ell,1} :=
		\mathcal{I}
		\left( \cdot; 
		\{ \texttt{i}_{\ell,1}^{\star} \} 
		,  {\rm span}\{ \psi_{\ell,1}^{\rm a,p}   \} \right)$
		\medskip
		
		\For{$m'=2, \ldots, m$}
		
		\State 
		Compute
		$r_{m'} \, = \, 
		\psi_{\ell,m'}^{\rm a,p}  -
		\mathcal{I}_{\ell, m'-1}  \left( \psi_{\ell,m'}^{\rm a,p} \right)$
		\medskip
		
		\State Set $
		\displaystyle{\texttt{i}_{\ell,m'}^{\star}:=\arg \max_{
				j  \in \{1,\ldots,N_{\ell}^{\rm p} \} 
			} \|  
			r_{m'}
			(  x_{\ell,j}^{\rm p}   ) 
			\|_2
		}$
		\medskip

		\State
		Update
		$ \mathcal{I}_{\ell,m'} :=
		\mathcal{I}
		\left( \cdot; 
		\{ \texttt{i}_{\ell,j}^{\star} \}_{j=1}^{m'} 
		,  {\rm span}\{ \psi_{\ell,j}^{\rm a,p}   \} _{j=1}^{m'} \right)$.
		\EndFor
		
	\end{algorithmic}
\end{algorithm}

\section{Analysis and interpretation  for linear coercive problems}
\label{sec:analysis}

We analyze the OS2 statement for linear coercive problems. To simplify the presentation, we consider the case with two subdomains   depicted in \Cref{fig:analysis_explanation}. We denote by $\left(\mathcal{X}, \|  \cdot \|_{\Omega} \right)$ the global ambient space such that
$H_0^1(\Omega) \subset \mathcal{X} \subset H^1(\Omega)$; given the ports $\Gamma_1, \Gamma_2$ (cf. \Cref{fig:analysis_explanation}), we define the bubble and port spaces:
$$
\mathcal{X}_{i,0} :=\left\{ v\in \mathcal{X}_i \, : \, v|_{\Gamma_i} = 0  \right\},
\quad
\mathcal{U}_i :=\left\{ v|_{\Gamma_i} \, : \, v  \in \mathcal{X}_i  \right\},
\quad
i=1,2.
$$
We introduce  the bilinear form $a:\mathcal{X} \times \mathcal{X} \to \mathbb{R}$  with continuity constant $\gamma$ and coercivity constant $\alpha>0$, and we introduce the linear functional  $f\in \mathcal{X}'$. Then, we introduce the model problem:
\begin{equation}
	\label{eq:lin_coercive}
	{\rm find} \; u^{\star} \in \mathcal{X} \, : \,
	a(u^{\star},v) = f(v) \quad \forall \; v\in \mathcal{X}.
\end{equation}
In   \cref{sec:port_formulation}, we derive the  port formulation of the problem \eqref{eq:lin_coercive};
in   \cref{sec:analysis_port_pb}, we present two important results for the port problem; in   \cref{sec:OS2_results} we exploit the results of the previous section to derive an \emph{a priori} bound for the OS2 statement;  in  \cref{sec:OS2_alternative_variational}, 
we comment on an alternative variational interpretation of the OS2 statement;
finally, in   \cref{sec:OSvsOS2_rates}, we derive explicit estimates for two representative model problems.

Given $v\in  \mathcal{X}_{i,0}$, we denote by $v^{\rm ext} \in \mathcal{X}$ the trivial extension of $v$ to $\Omega$ that is zero in $\Omega\setminus \Omega_i$.
We assume that $a$ and $f$ are associated to a differential (elliptic) problem; in particular, we assume that 
\begin{equation}
	\label{eq:technical_assumption}
	a(u,  v^{\rm ext}) = a\left(
	u\big|_{\Omega_i}, v
	\right),
	\quad
	\forall \, v\in \mathcal{X}_{i,0}.
\end{equation}
Note that by construction we have $f( v^{\rm ext}) = f\left( v \right)$ for all $v\in  \mathcal{X}_{i,0}$.

\subsection{Port formulation}
\label{sec:port_formulation}
We define the tensor-product space $\mathcal{U} = \mathcal{U}_1\times \mathcal{U}_2$ endowed with the inner product $\langle  w,v  \rangle = \sum_{i=1,2} (w_i,v_i)_{H^{1/2}(\Gamma_i)}$ and the induced norm
$\vertiii{\cdot} = \sqrt{\langle  \cdot , \cdot  \rangle}$. We introduce the local solution operators $T_i: \mathcal{U}_i \to \mathcal{X}_i$ and 
$G_i: \mathcal{X}' \to \mathcal{X}_{i,0}$ such that:
\begin{equation}
	\label{eq:analysis_Ti0}
	\left(T_i \lambda \right) \big|_{\Gamma_i} = \lambda,
	\qquad
	a(  T_i \lambda , v      ) = 0 \quad
	\forall \, v\in \mathcal{X}_{i,0};
\end{equation}
\begin{equation}
	\label{eq:analysis_Gi0}
	\left(G_i f \right)  \big|_{\Gamma_i} = 0,
	\qquad
	a(  G_i f   , v      ) = f(v) \quad
	\forall \, v\in \mathcal{X}_{i,0}.
\end{equation}
Since the elements of $\mathcal{X}_{i,0}$ can be trivially extended to zero in $\Omega \setminus  \Omega_i$, 
we have that the form $a$ is continuous and coercive in $\mathcal{X}_{i,0}$ with continuity and  coercivity constants bounded from above and below by $\gamma$ and $\alpha$, due to the fact that $\mathcal{X}_{i,0} \subset \mathcal{X}$.

Therefore, $T_i$ and $G_i$ are well-defined linear bounded operators.
By comparing the previous definitions with 
\eqref{eq:port2bubble}, we note
that
the affine operators
$F_i  := T_i  - E_i  + G_i f  $ correspond to the port-to-bubble maps that are exploited to derive the hybridized formulation in \cref{sec:formulation}:
we have
$ u^{\star}|_{\Omega_i}=
F_i \lambda_i^{\star} +E_i \lambda_i^{\star}  =T_i \lambda_i^{\star} +G_i f $, where $\lambda_i^{\star}\in \mathcal{U}_i$ is equal to $u^{\star}|_{\Gamma_i}.$

Given the trace operators
$\chi_{\Gamma_1} : \mathcal{X}_2 \to \mathcal{U}_1$, 
$\chi_{\Gamma_2} : \mathcal{X}_1 \to \mathcal{U}_2$, we introduce the operators
$T:\mathcal{U} \to \mathcal{U}$ and 
$G:\mathcal{X}' \to \mathcal{U}$ such that
\begin{subequations}
	\label{eq:port_problem}
	\begin{equation}
		\label{eq:port_problem_a}
		T\lambda 
		=
		\left[
		\begin{array}{l}
			\chi_{\Gamma_1}  T_2 \lambda_2 \\
			\chi_{\Gamma_2}  T_1  \lambda_1 \\
		\end{array}
		\right],
		\quad
		G f 
		=
		\left[
		\begin{array}{l}
			\chi_{\Gamma_1} G_2 f \\
			\chi_{\Gamma_2} G_1 f \\
		\end{array}
		\right],
		\quad
		\forall \, \lambda \in \mathcal{U}, \quad
		f \in \mathcal{X}'.
	\end{equation}
	Finally, we introduce the port problem:
	find $\lambda^{\star} \in \mathcal{U}$ such that
	\begin{equation}
		\label{eq:port_problem_b}
		a_{\rm p} (\lambda^{\star}  ,v ) = f_{\rm p}(v) \quad
		\forall \, v\in \mathcal{U},
		\;\;
		{\rm where} \;\;
		a_{\rm p} (\lambda ,v )
		:=
		\langle \lambda - T \lambda, v \rangle,
		\;\;
		f_{\rm p} ( v )
		:=
		\langle G f, v \rangle.
	\end{equation}
\end{subequations}

\begin{remark}
	\textbf{Connection with OS methods.}
	We can rewrite standard additive and multiplicative OS iterations using the operators introduced in \eqref{eq:port_problem}. In more detail, 
	multiplicative OS iterations can  be written as (see, e.g., \cite[Chapter 1]{quarteroni1999domain})
	$$
	\left[
	\begin{array}{ll}
		Id & 0 \\ 
		-\chi_{\Gamma_2} T_1 & Id  \\
	\end{array}
	\right]
	\lambda^{(k+1)}
	\, = \,
	\left[
	\begin{array}{ll}
		0 &  \chi_{\Gamma_1} T_2 \\
		0 &  0 \\
	\end{array}
	\right]
	\lambda^{(k)}
	\, + \,
	G f,
	\qquad
	k=1,2,\ldots, 
	$$
	while additive OS iterations can  be written as 
	$$
	\left[
	\begin{array}{ll}
		Id & 0 \\ 0 & Id  \\
	\end{array}
	\right]
	\lambda^{(k+1)}
	\, = \,
	\left[
	\begin{array}{ll}
		0 &  \chi_{\Gamma_1} T_2 \\
		\chi_{\Gamma_2} T_1 &  0 \\
	\end{array}
	\right]
	\lambda^{(k)}
	\, + \,
	G f,
	\qquad
	k=1,2,\ldots.
	$$
	These identities imply that any fixed point of the OS iterations satisfies \eqref{eq:port_problem_b} and thus OS and OS2 converge to the same limit as $k\rightarrow +\infty.$ As discussed in the introduction,  this connection between OS and OS2 formulations is valid for both linear and nonlinear problems; however, the analysis is strictly restricted to the linear case.
\end{remark}

\subsection{Analysis of the port problem}
\label{sec:analysis_port_pb}

\Cref{th:link_port2full} clarifies the relationship between the variational statement  \eqref{eq:lin_coercive} and the port problem \eqref{eq:port_problem}; on the other hand, 
\Cref{th:stability_port}  is key for the analysis of the OS2 ROM.
Proofs are postponed to \ref{sec:tedious_proofs}.
The results rely on the introduction of a 
partition-of-unity (PoU, \cite{babuvska1997partition}) $\{  \phi_i \}_{i=1}^{2} \subset {\rm Lip}(\Omega; \mathbb{R})$ associated with   $\{ \Omega_i  \}_{i=1}^{2}$ such that
$$
\sum_{i=1}^{2} \phi_i(x) = 1, 
\quad
\left\{
\begin{array}{ll}
	0\leq \phi_i(x)\leq 1  & \forall \, x\in \Omega, \\[2mm]
	\phi_i(x) = 0 & \forall \, x\notin \Omega_i, \\
\end{array}
\right.
\;\;
i=1,2.
$$  

\begin{proposition}
	\label{th:link_port2full}
	Let $u^{\star}$ be the solution to \eqref{eq:lin_coercive}. Then, $\lambda^{\star}  = \left( u^{\star}|_{\Gamma_1}, u^{\star}|_{\Gamma_2}   
	\right)$ solves \eqref{eq:port_problem_b}.
	Conversely, if $\lambda^{\star}$ is a solution 
	to \eqref{eq:port_problem_b}, then $
	u^{\star} = \sum_{i=1}^2 \left( T_i \lambda_i^{\star} + G_i f  \right) \phi_i$
	solves \eqref{eq:lin_coercive}.
\end{proposition}

\begin{proposition}
	\label{th:stability_port}
	Let the operator $T$ in \eqref{eq:port_problem_a} be compact. Then, the form $a_{\rm p}:\mathcal{U} \times \mathcal{U} \to \mathbb{R}$ defined in \eqref{eq:port_problem_b} is inf-sup stable and continuous, that is
	\begin{equation}
		\label{eq:stability_port}
		\alpha_{\rm p} = \inf_{w\in \mathcal{U}} 
		\sup_{v\in \mathcal{U}}  \frac{a_{\rm p}(w,v)}{\vertiii{w} \vertiii{v}}>0,
		\quad
		\gamma_{\rm p} = \sup_{w\in \mathcal{U}} 
		\sup_{v\in \mathcal{U}}  \frac{a_{\rm p}(w,v)}{\vertiii{w} \vertiii{v}}<\infty.
	\end{equation}
\end{proposition}

The proof of the compactness of the operator $T$ depends on the underlying PDE. For several problems, including the Laplace equation, the advection-diffusion-reaction equation,
the Stokes equations, and the Helmholtz's equation, we can prove compactness of the operator $T$ using Caccioppoli's inequalities: we refer to 
\cite[Appendix C]{Tad17} and also 
\cite{SmePat16} for further details.
We further observe that 
\Cref{th:stability_port} does not provide an explicit
relationship among 
the stability constant $\alpha_{\rm p}$ in  \Cref{th:stability_port}, the PDE of interest  and  the 
size of the overlap.
We envision that the derivation of explicit bounds for the stability constant $\alpha_{\rm p}$ in terms of the PDE of interest and the size of the overlap
will shed light on the underlying properties of the OS2 formulation and might also lead to new algorithmic developments.
We note that there is a vast body of works that address the derivation of sharp estimates for the convergence of  overlapping  Schwarz methods (see, e.g., \cite{chan1991geometry,ciaramella2017analysis}): the derivation of analogous results for this setting is beyond the scope of the present paper.

As discussed in \ref{sec:tedious_proofs}, proofs of Propositions \ref{th:link_port2full} and \ref{th:stability_port} rely on the fact that, if we introduce the spaces $\mathcal{X}_{1,2} = \{v|_{\Omega_1\cap \Omega_2} : v\in\mathcal{X} \}$ and 
$\mathcal{X}_{1,2}^0 = \{v \in \mathcal{X}_{1,2}: 
v|_{\Gamma_1\cup \Gamma_2}  = 0  \}$, the problem of finding
$u\in \mathcal{X}_{1,2}$ such that
$$
a(u,v) = 0 \quad \forall \, v\in \mathcal{X}_{1,2}^0, \;\;
u|_{\Gamma_1} = \lambda_1, \;\;
u|_{\Gamma_2} = \lambda_2, 
$$
admits a unique solution for any $(\lambda_1,\lambda_2)\in \mathcal{U}$. This result is trivial for coercive problems, but it is significantly less trivial --- and requires additional assumptions --- for inf-sup stable problems and is not addressed in this work. 
On the other hand, we envision that the analysis for nonlinear PDEs requires more sophisticated tools and is   beyond the scope of this work.

\subsection{Analysis of the OS2 statement}
\label{sec:OS2_results}
We consider the following  OS2 formulation for the linear problem \eqref{eq:lin_coercive}:
\begin{equation}
	\label{eq:OS2_lincoercive_mod}
	{\rm find} \; \widehat{\lambda} = {\rm arg} \min_{\lambda \in \mathcal{Z}^{\rm p} } \vertiii{ \lambda - \widehat{T} \lambda - \widehat{G} f   }.
\end{equation} 
Note that \eqref{eq:OS2_lincoercive_mod} corresponds to the OS2 statement \eqref{eq:hybridizedOS2} with the important difference that we replace the $L^2 $ norm with the $H^{1/2}$ norm  $\vertiii{\cdot}$.
In particular, 
in our work,
the space $\mathcal{Z}^{\rm p}$ in \eqref{eq:OS2_lincoercive_mod} 
is given by the tensor product of the local port spaces, 
$\mathcal{Z}^{\rm p} = \mathcal{Z}_1^{\rm p}  \times \mathcal{Z}_2^{\rm p}$,
and $\widehat{T},\widehat{G}$ are associated to the approximate local solution operators that are obtained by Galerkin projection.
We further introduce the OS2 formulation with perfect local operators:
\begin{equation}
	\label{eq:OS2_lincoercive_mod_perfect}
	{\rm find} \; \widetilde{\lambda} = {\rm arg} \min_{\lambda \in \mathcal{Z}^{\rm p} } \vertiii{ \lambda - T \lambda -  {G} f   }.
\end{equation} 

We observe that \eqref{eq:OS2_lincoercive_mod_perfect} corresponds to the minimum residual formulation of the port problem \eqref{eq:port_problem};  we have indeed
$$
\vertiii{ \lambda - T \lambda - {G} f   }
=
\sup_{v\in \mathcal{U} } 
\frac{\langle   \lambda - T \lambda - {G} f, v \rangle}{\vertiii{v}} 
=
\sup_{v\in \mathcal{U} } 
\frac{ a_{\rm p}( \lambda,  v ) - f_{\rm p}(v)}{\vertiii{v}}. 
$$
Recalling the result in \cite{xu2003some}, we thus have
\begin{equation}
	\label{eq:quasi_optimality_perfect_OS2}
	\vertiii{ \lambda^{\star} -\widetilde{\lambda}   }
	\leq
	\frac{\gamma_{\rm p}}{\alpha_{\rm p}}
	\inf_{\lambda\in \mathcal{Z}^{\rm p}} \
	\vertiii{ \lambda^{\star} - {\lambda}   },
\end{equation}
which proves the quasi-optimality of the OS2 statement with perfect local operators \eqref{eq:OS2_lincoercive_mod_perfect}.

To estimate the error $ \vertiii{ \widehat{\lambda} -\widetilde{\lambda}   }$, we resort to a perturbation analysis.
We denote by $\widehat{\alpha}_{\rm p}$ and $\widehat{\gamma}_{\rm p}$  the stability and continuity constants associated with the problem \eqref{eq:OS2_lincoercive_mod}: it is possible to resort to a perturbation analysis to estimate these constants; since the argument is completely standard, we omit the details.
We define the quantities $\varepsilon_{\rm T}$ and  $\varepsilon_{\rm G}$ as follows:
\begin{equation}
	\label{eq:local_errors_source}
	\varepsilon_{\rm T} : = 
	\sup_{\psi\in \mathcal{Z}^{\rm p}}
	\frac{\vertiii{ (T-\widehat{T})\psi   }}{\vertiii{ \psi   }},
	\quad
	\varepsilon_{\rm G} : = 
	\vertiii{(G-\widehat{G})f}.
\end{equation}
Then, it is possible to show that  
\begin{equation}
	\label{eq:perturbation_bound}
	\vertiii{ \widetilde{\lambda}-\widehat{\lambda}    }
	\leq
	\frac{1}{\alpha_{\rm p}^2}
	\left(
	M \left( \gamma_{\rm p} + \widehat{\gamma}_{\rm p}  \right)
	\frac{\vertiii{\widehat{G} f}}{\widehat{\alpha}_{\rm p}}
	\varepsilon_{\rm T}
	\,+\,
	\sqrt{M}
	\left(  \widehat{\gamma}_{\rm p} \varepsilon_{\rm G} 
	+  \vertiii{\widehat{G} f} \varepsilon_{\rm T}
	\right)
	\right).
\end{equation}
We postpone the proof of \eqref{eq:perturbation_bound} to 
\ref{sec:tedious_proofs}.

By combining \eqref{eq:perturbation_bound} with \eqref{eq:quasi_optimality_perfect_OS2}, we obtain the following result.
We observe that \eqref{eq:error_bound_OS2} is the sum of two terms: the first term is associated with the approximation properties of the port space, while the second term is directly linked to the accuracy of the local solution operators.

\begin{proposition}
	\label{th:error_bound_OS2}
	Let $\gamma_{\rm p}, \alpha_{\rm p}$ be the continuity and stability constants of the form $a_{\rm p}$ and let 
	$\widehat{\gamma}_{\rm p}, \widehat{\alpha}_{\rm p}$ be the continuity and stability constants of the form $\widehat{a}_{\rm p}(\lambda,v) = \langle \lambda  - \widehat{T} \lambda \, , \,  v \rangle$. Given the $M$-dimensional space $\mathcal{Z}^{\rm p} \subset \mathcal{U}$, we have
	\begin{equation}
		\label{eq:error_bound_OS2}
		\vertiii{  {\lambda}^{\star}-\widehat{\lambda}    }
		\leq
		\frac{1}{\alpha_{\rm p}}
		\left(
		\gamma_{\rm p}
		\inf_{\lambda\in \mathcal{Z}^{\rm p}} \
		\vertiii{ \lambda^{\star} - {\lambda}   } \, + \,
		\frac{1}{\alpha_{\rm p}}  \left(
		M \left( \gamma_{\rm p} + \widehat{\gamma}_{\rm p}  \right)
		\frac{\vertiii{\widehat{G} f}}{\widehat{\alpha}_{\rm p}}
		\varepsilon_{\rm T}
		\,+\,
		\sqrt{M}
		\left(  \widehat{\gamma}_{\rm p} \varepsilon_{\rm G} 
		+  \vertiii{\widehat{G} f} \varepsilon_{\rm T}
		\right)
		\right)
		\right).
	\end{equation} 
\end{proposition}

\subsection{Alternative variational interpretation of the OS2 statement}
\label{sec:OS2_alternative_variational}

Following \cite{lions1988schwarz}, we might also consider the alternative variational framework of the OS limit formulation
(see also \cite[Chapter 1.5.2]{quarteroni1999domain}):
find $(u_1^{\rm b}, u_1^{\rm p}, u_2^{\rm b}, u_2^{\rm p}) \in  \bigotimes_{i=1}^2  \mathcal{X}_{i,0} \times \mathcal{U}_i$ such that
\begin{equation} 
	\label{eq:Lions_variational_formulation}
	\left\{
	\begin{array}{ll}
		a\left(   u_i^{\rm b} + E_i u_i^{\rm p}, \, v_i \right)
		= f(v_i) & \forall \, v_i\in \mathcal{X}_{i,0}, \;\;i=1,2 ; \\[3mm]
		\left(
		u_1^{\rm p}
		-
		\chi_{\Gamma_1} \left(
		u_2^{\rm b}
		+ 
		E_2 u_2^{\rm p}
		\right), \; \psi_1
		\right)_{H^{1/2}(\Gamma_1)}
		+
		\left(
		u_2^{\rm p}
		-
		\chi_{\Gamma_2} \left(
		u_1^{\rm b}
		+
		E_1 u_1^{\rm p}
		\right), \; \psi_2
		\right)_{H^{1/2}(\Gamma_2)}
		= 0
		&
		\forall \, \psi=(\psi_1,\psi_2) \in \mathcal{U}; 
		\\
	\end{array}
	\right.
\end{equation}
where $E_1,E_2$ are the extension operators, $u_1^{\rm b},u_2^{\rm b}$ are the bubble solutions and  $u_1^{\rm p},u_2^{\rm p}$ are the port solutions.
Given the reduced spaces
$\mathcal{Z}_i^{\rm b} \subset \mathcal{X}_{i,0}$ and
$\mathcal{Z}_i^{\rm p} \subset \mathcal{U}_{i}$,
and the  approximate port-to-bubble maps
$\widehat{F}_i = \widehat{T}_i + \widehat{G}_i f - E_i$,
for $i=1,2$,  
the reduced-order OS2 formulation can be stated as follows:
find $(\widehat{u}_1^{\rm b}, \widehat{u}_1^{\rm p}, \widehat{u}_2^{\rm b}, \widehat{u}_2^{\rm p}) \in  \bigotimes_{i=1}^2  \mathcal{Z}_i^{\rm b} \times \mathcal{Z}_i^{\rm p} $ such that
\begin{subequations}
	\label{eq:Lions_variational_formulation_OS2}
	\begin{equation}
		\left\{
		\begin{array}{ll}
			a\left(   \widehat{u}_i^{\rm b} + E_i \widehat{u}_i^{\rm p}, \, v_i \right)
			= f(v_i) & \forall \, v_i\in \mathcal{Z}_{i}^{\rm b}, \;\;i=1,2 ; \\[3mm]
			\left(
			\widehat{u}_1^{\rm p}
			-
			\chi_{\Gamma_1} \left(
			\widehat{u}_2^{\rm b}
			+
			E_2  \widehat{u}_2^{\rm p}
			\right), \; \psi_1
			\right)_{H^{1/2}(\Gamma_1)}
			+
			\left(
			\widehat{u}_2^{\rm p}
			-
			\chi_{\Gamma_2} \left(
			\widehat{u}_1^{\rm b}
			+
			E_1  \widehat{u}_1^{\rm p}
			\right), \; \psi_2
			\right)_{H^{1/2}(\Gamma_2)}
			= 0
			&
			\forall \, \psi=(\psi_1,\psi_2) \in \widetilde{\mathcal{Z}}^{\rm p}; 
			\\
		\end{array}
		\right.
	\end{equation}
	where
	$\widetilde{\mathcal{Z}}^{\rm p} \subset \mathcal{U}$ is the $M$-dimensional space given by
	\begin{equation}
		\widetilde{\mathcal{Z}}^{\rm p}
		=\left\{
		\left(
		\zeta_1^{\rm p} -   \chi_{\Gamma_1}  \widehat{T}_2(\zeta_2^{\rm p}), \;
		\zeta_2^{\rm p} -   \chi_{\Gamma_2}  \widehat{T}_{1}(\zeta_{1}^{\rm p})
		\right)
		\,: \,
		\zeta_i^{\rm p}  \in \mathcal{Z}_i^{\rm p}, i=1,2
		\right\},
	\end{equation}
	and $\widehat{T}_i \zeta$ satisfies
	$\widehat{T}_i \zeta = u^{\rm b}(\zeta) + E_i \zeta$ with  $u_i^{\rm b}(\zeta) \in \mathcal{Z}_i^{\rm b}$ and 
	$a(u_i^{\rm b}(\zeta) + E_i \zeta, v) = 0$ for all $v\in \mathcal{Z}_i^{\rm b}$.
\end{subequations}

The proof of 
\eqref{eq:Lions_variational_formulation_OS2} is
straightforward, and it is provided for completeness 
in \ref{sec:tedious_proofs}.
Note that the OS2 statement reads as a Petrov-Galerkin projection of \eqref{eq:Lions_variational_formulation} for a suitable choice of the test space 
$\widetilde{\mathcal{Z}}^{\rm p}$.
We envision that \eqref{eq:Lions_variational_formulation_OS2} could be exploited to devise an alternative error analysis for the OS2 statement. We do not address this issue in the present work.

\subsection{Explicit convergence rates for two one-dimensional model problems}
\label{sec:OSvsOS2_rates}
Given $\Omega=(-1,1)$ and the partition $\Omega_1 = (-1,\delta)$,
$\Omega_2 = (-\delta,1)$, we study the convergence of (multiplicative) OS and OS2 for the problems
\begin{subequations}
	\begin{equation}
		\label{eq:1Dpb_a}
		\left\{
		\begin{array}{ll}
			u'' = 2  &  {\rm in} \;  \Omega, \\
			u(-1)=u(1)=1; &    \\
		\end{array}
		\right.
	\end{equation}
	and
	\begin{equation}
		\label{eq:1Dpb_b}
		\left\{
		\begin{array}{ll}
			-u'' +  \gamma u' = 0  &  {\rm in} \;  \Omega, \\
			u(-1)=0, \;\; u(1)=1; &    \\
		\end{array}
		\right.
	\end{equation}
	in the limit $|\delta | \ll 1$. 
	
	The analysis can be readily extended to the additive OS method. 
	
\end{subequations}
For OS2, we resort to the gradient descent method with optimal choice of the step size, and to the Gauss-Newton method (OS2-GN) --- the choice of 
the gradient descent method is intended to simplify calculations (compared to quasi-Newton methods).
The motivation of this analysis is twofold: first, we show that the use of gradient-based methods --- as opposed to Gauss-Newton --- is increasingly sub-optimal as $\delta \rightarrow 0$; second, we provide explicit estimates for the constants $\alpha_{\rm p}$ and $\gamma_{\rm p}$ of Proposition \ref{th:stability_port} for two representative model problems.

We denote by $\widehat{u}_i$ the approximation of the solution in $\Omega_i$ for $i=1,2$; we define
$\beta_1 = \widehat{u}_1(\delta)$ and
$\beta_2 = \widehat{u}_2(-\delta)$. 
We can show that OS and OS2 iterations can be written as 
$$
\boldsymbol{\beta}^{(k)}
=
\mathbf{P}_{\delta}^{\rm os} \boldsymbol{\beta}^{(k-1)}
+
\mathbf{F}_{\delta}^{\rm os},
\quad
\boldsymbol{\beta}^{(k)}
=
\mathbf{P}_{\delta}^{\rm os2} \boldsymbol{\beta}^{(k-1)}
+
\mathbf{F}_{\delta}^{\rm os2},
$$
for $k=1,2,\ldots$ and suitable choices of $\big( \mathbf{P}_{\delta}^{\rm os}, \mathbf{F}_{\delta}^{\rm os} \big)$ and
$\big( \mathbf{P}_{\delta}^{\rm os2}, \mathbf{F}_{\delta}^{\rm os2} \big)$. On the other hand, since the problems are linear, OS2-GN reduces to a direct method and can be stated as 
$$
\mathbf{A}_{\delta} 
\boldsymbol{\beta} 
=
\mathbf{F}_{\delta} 
$$
for suitable choices of $\big( \mathbf{A}_{\delta}, \mathbf{F}_{\delta} \big)$.

In \ref{sec:tedious_proofs}, we show that the spectral radii $\rho_{\delta}^{\rm os}$ and $\rho_{\delta}^{\rm os2}$ of the transition matrices 
$\mathbf{P}_{\delta}^{\rm os}$ and 
$\mathbf{P}_{\delta}^{\rm os2}$ satisfy
\begin{subequations}
	\label{eq:tedious_proof}
	\begin{equation}
		\label{eq:tedious_proof_a}
		\begin{array}{lll}
			\displaystyle{
				\rho_{\delta}^{\rm os} \sim 1 - 4 \delta},
			&
			\displaystyle{
				\rho_{\delta}^{\rm os2} \sim 1 - 4 \delta^2
			}
			&
			{\rm for} \;\; \eqref{eq:1Dpb_a}; \\[2mm]
			\displaystyle{
				\rho_{\delta}^{\rm os} \sim 1 - 2
				\frac{e^{\gamma}+1}{e^{\gamma}-1}
				\gamma
				\delta},
			&
			\displaystyle{
				\rho_{\delta}^{\rm os2} \sim 1 -  
				\frac{e^{\gamma}+2}{8(e^{\gamma}-1)}
				\gamma^2
				\delta^2 
			}
			&
			{\rm for} \;\; \eqref{eq:1Dpb_b}; \\[2mm]
		\end{array}
	\end{equation}
	while the condition number of the linear system associated to OS2-GN satisfies
	\begin{equation}
		\label{eq:tedious_proof_b}
		\begin{array}{ll}
			\displaystyle{
				{\rm cond} \left( 
				\mathbf{A}_{\delta} \right)
				= \frac{1}{\delta} },
			&
			{\rm for} \;\; \eqref{eq:1Dpb_a}; \\[2mm]
			\displaystyle{
				{\rm cond} \left( 
				\mathbf{A}_{\delta}  \right)
				\sim   
				\frac{ 4 (e^{\gamma}-1)   }{4 (e^{\gamma}+2) \gamma}
				\delta^{-1} 
			}
			&
			{\rm for} \;\; \eqref{eq:1Dpb_b}; \\[2mm]
		\end{array}
	\end{equation}
	and the constants $\alpha_{\rm p}$ and $\gamma_{\rm p}$ defined in \Cref{th:stability_port} satisfy
	\begin{equation}
		\label{eq:tedious_proof_c}
		\begin{array}{ll}
			\displaystyle{
				\alpha_{\rm p}
				= \frac{2\delta}{1+\delta},
				\;\;
				\gamma_{\rm p}
				= \frac{2}{1+\delta}
			},
			&
			{\rm for} \;\; \eqref{eq:1Dpb_a}; \\[4mm]
			\displaystyle{
				\alpha_{\rm p}  
				\sim
				\frac{4 (e^{\gamma}+2) \gamma\delta}{ 2 (e^{\gamma}-1)   },
				\quad
				\gamma_{\rm p}  
				\sim 2,
			}
			&
			{\rm for} \;\; \eqref{eq:1Dpb_b}; \\[2mm]
		\end{array}
	\end{equation}
\end{subequations}
As expected, OS, OS2 and OS2-GN become increasingly ill-conditioned as $\delta$ decreases to zero and do not converge for $\delta=0$; however, we observe that for small values of $\delta$ OS exhibits significantly faster convergence rates than OS2 based on the  gradient-descent method: this observation further strengthens the importance of exploiting the least-square structure of the OS2 statement.

\section{Numerical results}
\label{sec:numerics}
\label{sec:numerics}

\subsection{Assessment metrics and training parameters}

We train the CB-ROM based on $n_{\rm train}=70$ global parameters $\Xi_{\rm train}=\{   \mu^{(k)} \}_{k=1}^{n_{\rm train}}$ such that
$$
(E_1^{(k)}, E_2^{(k)},E_3^{(k)},s^{(k)}) \overset{\rm iid}{\sim} {\rm Uniform} \left(
[25,30]\times [10,20]^2 \times [0.4,1]
\right),
\quad
Q_{\rm a}^{(k)}
\overset{\rm iid}{\sim} {\rm Uniform} \left(
\{ 2,\ldots, 7 \}
\right);
$$
on the other hand, we assess performance based on $n_{\rm test}=20$ out-of-sample global parameters 
$\Xi_{\rm test}=\{   \widetilde{\mu}^{(j)} \}_{j=1}^{n_{\rm test}}$ generated using the same distribution. In view of the assessment, we also define the  PoU $\{  \phi_i \}_{i=1}^{N_{\rm dd}} \subset {\rm Lip}(\Omega; \mathbb{R})$ associated with the partition $\{ \Omega_i  \}_{i=1}^{N_{\rm dd}}$ such that
$$
\sum_{i=1}^{N_{\rm dd}} \phi_i(x) = 1, 
\quad
\left\{
\begin{array}{ll}
	0\leq \phi_i(x)\leq 1  & \forall \, x\in \Omega, \\[2mm]
	\phi_i(x) = 0 & \forall \, x\notin \Omega_i, \\
\end{array}
\right.
\;\;
i=1,\ldots,N_{\rm dd}.
$$
Given $u\in \mathcal{X} :=\bigotimes_{i=1}^{N_{\rm dd}} \mathcal{X}_i$, we define the PoU operator
\begin{equation}
	\label{eq:pum_operator}
	\texttt{P}_{\rm pu}[u]:=\displaystyle{\sum_{i=1}^{N_{\rm dd}}}\phi_i \, u_i \in H^1(\Omega).
\end{equation}
Note that we omit the dependence of $\{ \phi_i \}_i$ and also $N_{\rm dd}$ on the parameter to shorten notation.
Finally, we define the out-of-sample average and maximum prediction errors
\begin{subequations}
	\begin{equation}
		\label{eq:Eavg}
		E_{\rm avg}
		:=
		\frac{1}{n_{\rm test}} \; 
		\sum_{\mu \in \Xi_{\rm test}}
		\;
		\frac{\|  \texttt{P}_{\rm pu}[u_{\mu}^{\rm hf}]  - \texttt{P}_{\rm pu}[\widehat{u}_{\mu}]     \|_{H^1(\Omega)}}{
			\|  \texttt{P}_{\rm pu}[u_{\mu}^{\rm hf}] \|_{H^1(\Omega)}
		},
	\end{equation}
	\begin{equation}
		\label{eq:Emax}	
		E_{\rm max}
		:=
		\max_{\mu \in \Xi_{\rm test}}
		\;
		\frac{\|  \texttt{P}_{\rm pu}[u_{\mu}^{\rm hf}]  - \texttt{P}_{\rm pu}[\widehat{u}_{\mu}]     \|_{H^1(\Omega)}}{
			\|  \texttt{P}_{\rm pu}[u_{\mu}^{\rm hf}] \|_{H^1(\Omega)}
		}.
	\end{equation}
	
\end{subequations}
As mentioned in section \ref{sec:methods}, we here resort to the HF CB solver to generate HF data for training and test, to simplify interpretation of the numerical results. In several figures, we compare the prediction error \eqref{eq:Eavg} with the 
error associated with the  mapped $H^1(\Omega_{\texttt{L}_i}^{\rm a})$ projection of $u_{\mu}^{\rm hf} \circ \Phi_i$, for $i=1,\ldots,N_{\rm dd}$,
\begin{equation}
	\label{eq:Eavg_opt}
	E_{\rm avg}^{\rm opt}
	:=
	\frac{1}{n_{\rm test}} \; 
	\sum_{\mu \in \Xi_{\rm test}}
	\;
	\frac{\|  \texttt{P}_{\rm pu}[u_{\mu}^{\rm hf}]  - \texttt{P}_{\rm pu}[\widehat{u}_{\mu}^{\rm opt}]     \|_{H^1(\Omega)}}{
		\|  \texttt{P}_{\rm pu}[u_{\mu}^{\rm hf}] \|_{H^1(\Omega)}
	},
	\quad
	{\rm with} \;
	\left( \widehat{u}_{\mu}^{\rm opt} \right)_i
	= \left( \Pi_{  \mathcal{Z}_{\texttt{L}_i}^{\rm a,b} \cup     \mathcal{W}_{\texttt{L}_i}^{\rm a,p}     } u_{\mu}^{\rm hf} \circ \Phi_i \right)\circ \Phi_i^{-1},
\end{equation}
for $i=1,\ldots,N_{\rm dd}.$
Note that \eqref{eq:Eavg_opt} is not  optimal --- that is, it is not the relative $H^1(\Omega)$ projection error associated with the instantiated spaces --- but it can be shown to be quasi-optimal exploiting 
\cite[Theorem 1]{babuvska1997partition}. We omit the details.

We resort to a P2 FE discretization with 
$N_{\rm int}^{\rm e}=1120$ and $N_{\rm ext}^{\rm e}=3960$ elements, and 
$N_{\rm int}^{\rm p}=272$ and $N_{\rm ext}^{\rm p}=200$  port quadrature points. 
We emphasize that the HF component-based discretization is constructed to ensure that the local grids match exactly for $Q_{\rm a}=Q_{\rm ref}$ (cf. Figure \ref{fig:picture_components}); however, we remark that internal and external meshes do not lead to a global conforming discretization for any other value of $Q_{\rm a}$. 

All simulations are performed in Matlab 2020b on a commodity laptop.
The implementation of the method does not resort to any parallelization of offline and online solves.

\subsection{Reduced-order model with HF quadrature}

We show the performance of the OS2 ROM without hyper-reduction. 
First, we show the behavior of the percentage of retained energy of the POD eigenvalues $\{ \lambda_i \}_i$ of the Gramian matrix associated with  the snapshot set.
To facilitate visualization, we show the average in-sample error
$E_{\rm n} = 1 - \frac{\sum_{i=1}^n  \lambda_i  }{\sum_{j=1}^{n_{\rm train}}  \lambda_j}$ for several values of $n$, for port and bubble components, and for the two 
archetype components. We observe  that the POD eigenvalues decay extremely rapidly, for both components.

\begin{figure}[H]
	\centering
	\subfloat[$\Omega_{\rm int}^{\rm a}$]{
		\begin{tikzpicture}
			\begin{loglogaxis}[
				xmode=linear,
				ymode=log,
				grid=both,
				minor grid style={gray!25},
				major grid style={gray!25},
				title = {},
				xlabel={$n$},
				ylabel={avg in-sample error},
				ymin= 0.000000001,
				ymax=0.1,
				width=0.43\textwidth,
				xtick={2,4,6,8,10,12,14,16},
				]
				]
				\addplot[ 
				black,
				mark=o,
				mark options={ color = black},
				]  %
				table{dat/pod/port_int.dat};
				\addlegendentry{\small{port}};
				\addplot[ 
				black,
				mark=square,
				mark options={ color = red},
				]  %
				table{dat/pod/bubble_int.dat};
				\addlegendentry{\small{bubble}};
			\end{loglogaxis}
		\end{tikzpicture}
	}
	\subfloat[$\Omega_{\rm ext}^{\rm a}$]{
		\begin{tikzpicture}
			\begin{loglogaxis}[
				xmode=linear,
				ymode=log,
				grid=both,
				minor grid style={gray!25},
				major grid style={gray!25},
				title = {},
				xlabel={$n$},
				ylabel={avg in-sample error},
				ymin= 0.000000001,
				ymax=0.1,
				width=0.43\textwidth,
				xtick={2,4,6,8,10,12,14,16},
				]
				]
				\addplot[ 
				black,
				mark=o,
				mark options={ color = black},
				]  %
				table{dat/pod/port_ext.dat};
				\addlegendentry{\small{port}};
				\addplot[ 
				black,
				mark=square,
				mark options={color = red},
				]  %
				table{dat/pod/bubble_ext.dat};
				\addlegendentry{\small{bubble}};
			\end{loglogaxis}
		\end{tikzpicture}
	}
	\caption{behavior of the average squared in-sample error $E_{\rm n} = 1 - \frac{\sum_{i=1}^n  \lambda_i  }{\sum_{j=1}^{n_{\rm train}}  \lambda_j}$ for several values of $n$, for port and bubble components, and for the two 
		archetype components.}
	\label{fig:pod_eigen}
\end{figure}
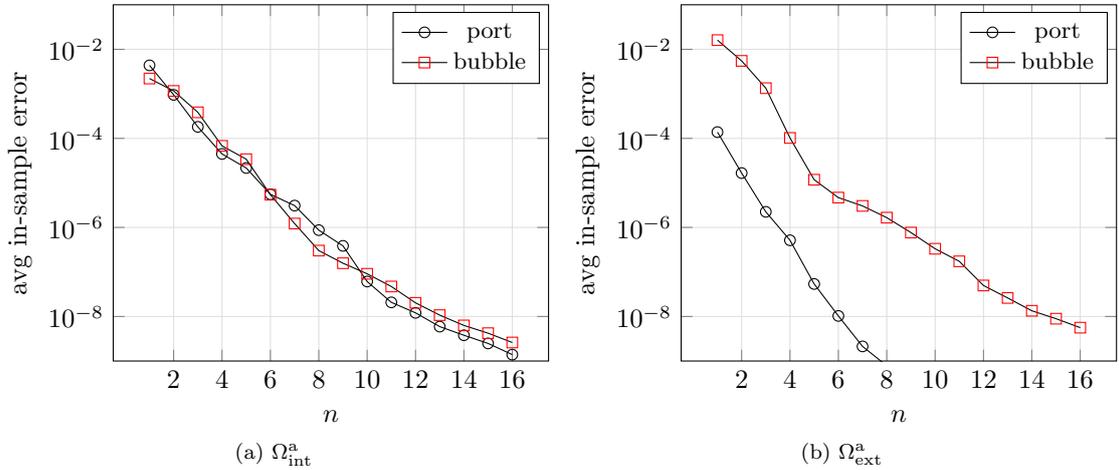

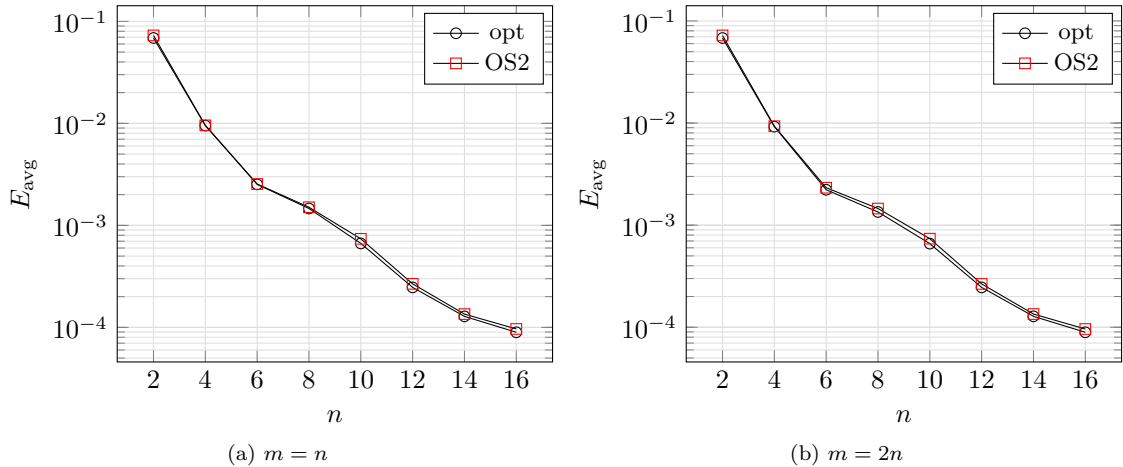
\begin{figure}[H]
	\centering
	\subfloat[$m=n$]{
		\begin{tikzpicture}
			\begin{loglogaxis}[
				xmode=linear,
				ymode=log,
				grid=both,
				minor grid style={gray!25},
				major grid style={gray!25},
				title = {},
				xlabel={$n$},
				ylabel={$E_{\rm avg}$},
				width=0.43\textwidth,
				xtick={2,4,6,8,10,12,14,16},
				]
				]
				\addplot[ 
				black,
				mark=o,
				mark options={ color = black},
				]  %
				table{dat/proj.dat};
				\addlegendentry{\small{opt}};
				\addplot[ 
				black,
				mark=square,
				mark options={ color = red},
				]  %
				table{dat/Gal.dat};
				\addlegendentry{\small{OS2}};
			\end{loglogaxis}
		\end{tikzpicture}
	}
	\subfloat[$m=2n$]{
		\begin{tikzpicture}
			\begin{loglogaxis}[
				xmode=linear,
				ymode=log,
				grid=both,
				minor grid style={gray!25},
				major grid style={gray!25},
				title = {},
				xlabel={$n$},
				ylabel={$E_{\rm avg}$},
				width=0.43\textwidth,
				xtick={2,4,6,8,10,12,14,16},
				]
				]
				\addplot[ 
				black,
				mark=o,
				mark options={ color = black},
				]  %
				table{dat/proj_2n.dat};
				\addlegendentry{\small{opt}};
				\addplot[ 
				black,
				mark=square,
				mark options={ color =red},
				]  %
				table{dat/Gal_2n.dat};
				\addlegendentry{\small{OS2}};
			\end{loglogaxis}
		\end{tikzpicture}
	}
	\caption{out-of-sample performance of OS2 ROM without hyper-reduction for several values of $m$, with $n=m$ and $n=2m$; comparison with optimal (``opt'')
		average error   $E_{\rm avg}^{\rm opt}$ \eqref{eq:Eavg_opt}}.
	\label{fig:err_nohyper}
\end{figure}

In \Cref{fig:err_nohyper},  we compare the average error 
$E_{\rm avg}$ \eqref{eq:Eavg} associated with the OS2 ROM  for several values of $m$ and $n=m$ and $n=2m$, with the average error 
$E_{\rm avg}^{\rm opt}$ \eqref{eq:Eavg_opt} obtained through projection. We  observe that the OS2 ROM achieves near-optimal performance for all choices of the port and bubble ROBs.
We also observe that doubling the number of port modes $m$ by keeping the same number of bubble modes $n$ does not lead to relevant differences in terms of both projection and OS2 prediction error. In the remainder, we 
set $m=n$.

\begin{figure}[h!]
	\centering
	\subfloat[geometry configuration and slice ]{
		\includegraphics[width=0.43\textwidth]
		{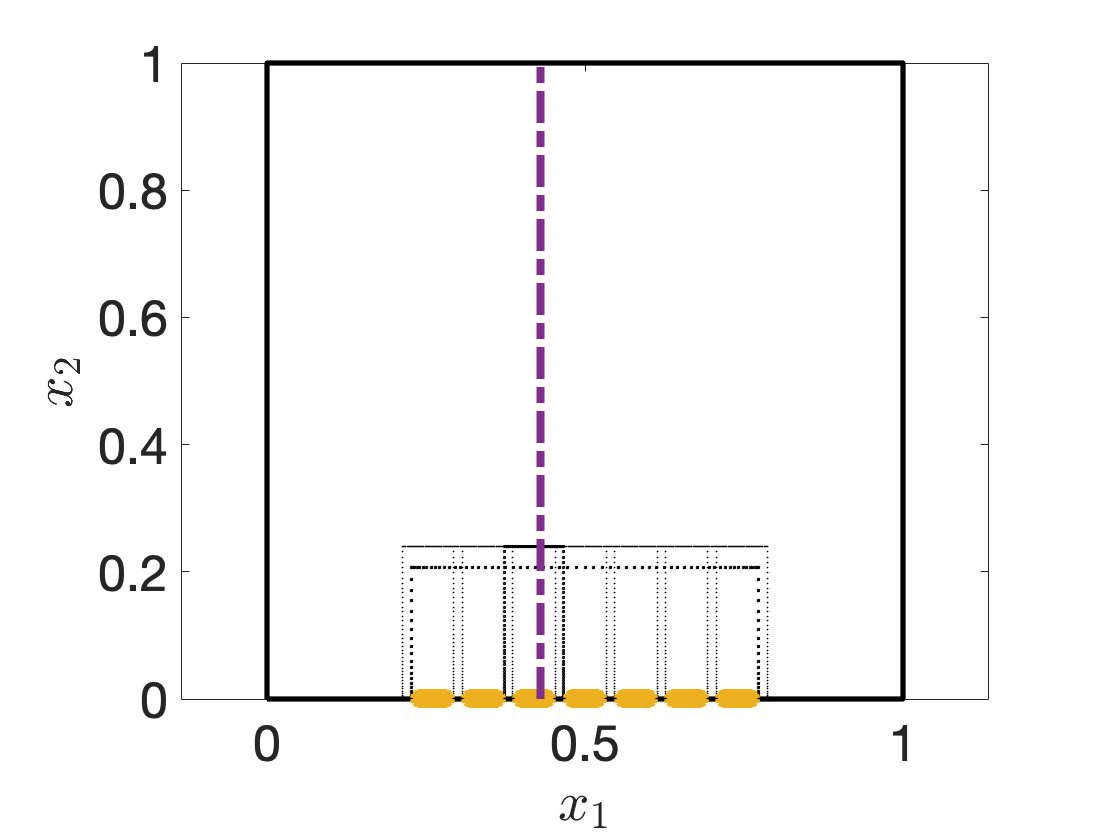}
	}
	\quad
	\subfloat[horizontal displacement ]{
		\includegraphics[width=0.42\textwidth]
		{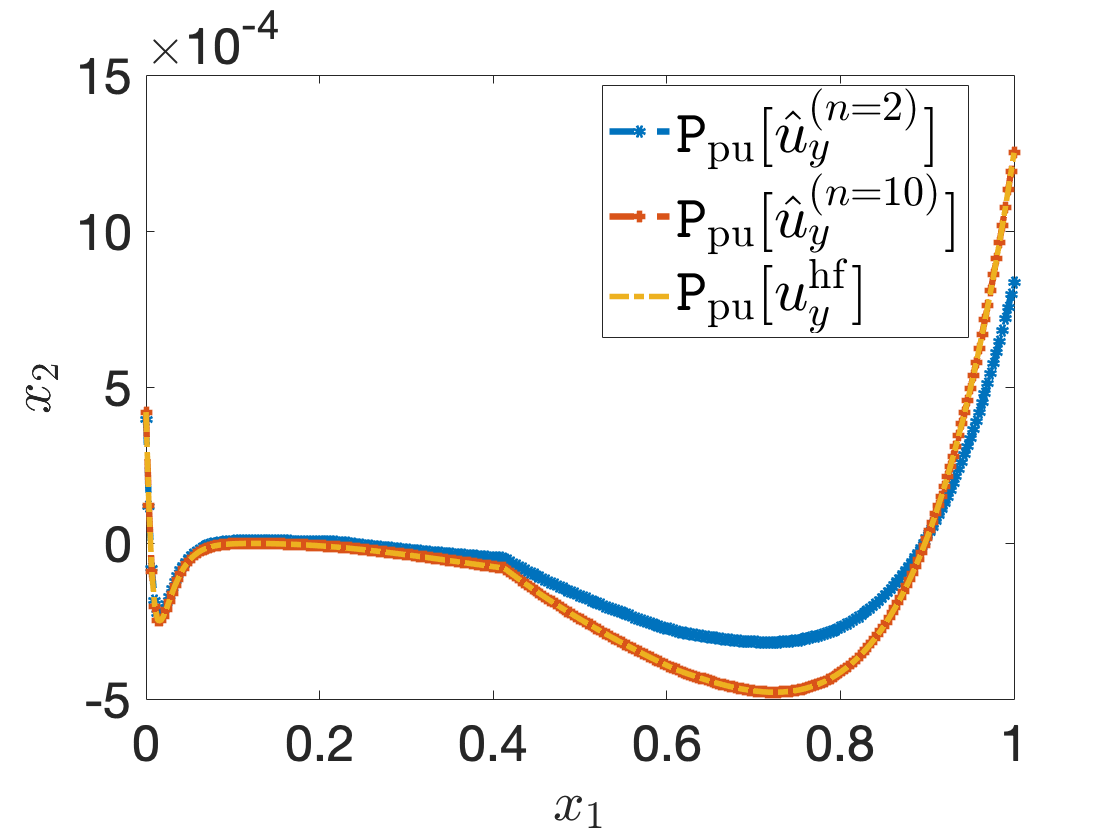}
	}
	\subfloat[ vertical displacement]{
		\includegraphics[width=0.42\textwidth]
		{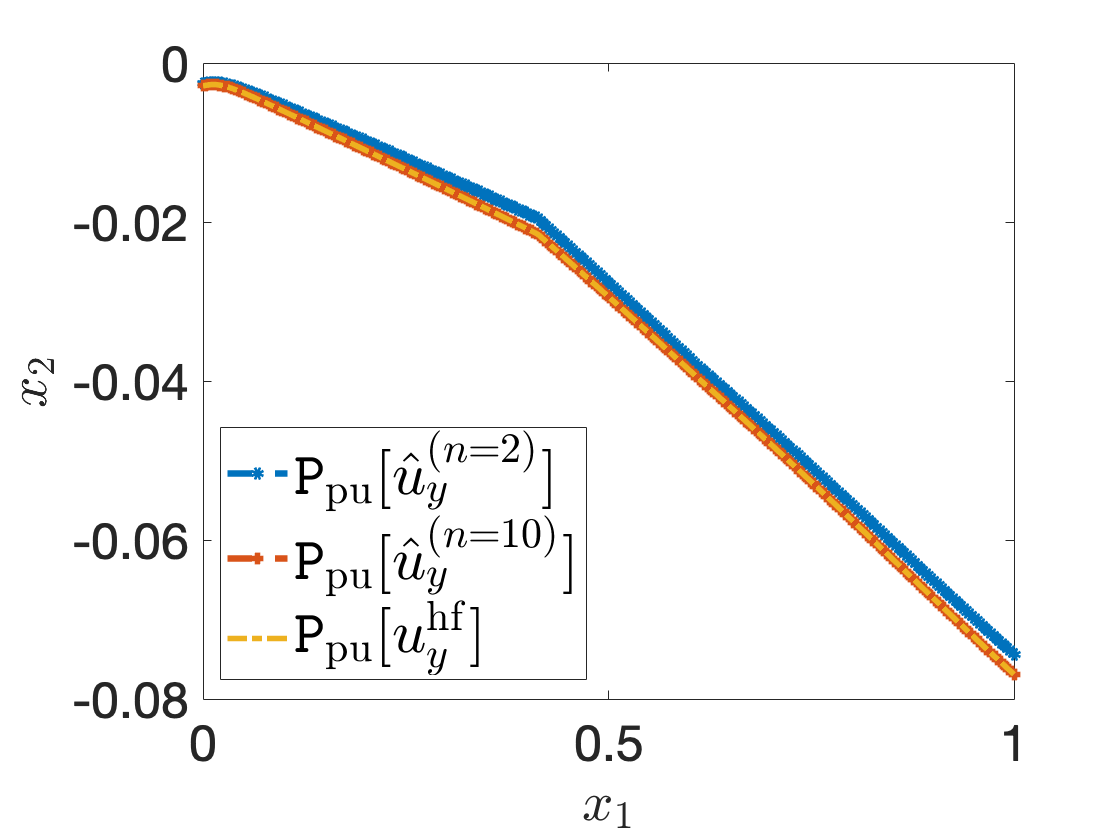}
	}
	
	\caption{
		visualization of the horizontal and vertical displacement components for a vertical slice.
	}
	\label{fig:vis_slices}
\end{figure}

\Cref{fig:vis_slices} shows the behavior of the solution over a vertical slice of the domain for a test configuration with $Q_{\rm a}=7$; boundaries of the $Q_{\rm a}$ internal subdomains associated with repositories and the external subdomain are marked as black dots in \ref{fig:vis_slices} (a); the vertical slice, drawn as a purple dashed line, corresponds to points $(x,y)$ such that $x=\bar{x}=0.43$, $0\leq y \leq 1$.
Points of the slice belong to either the instantiated component $\Omega_3$ or $\Omega_8$ (or both). We apply the partition of unity operator \eqref{eq:pum_operator} to generate globally-defined solutions. We compute therefore approximate solutions $\texttt{P}_{\rm pu}[\hat{u}_{\star}^{(n=2)}]$ $\texttt{P}_{\rm pu}[\hat{u}_{\star}^{(n=10)}]$ corresponding to two choices of the ROB size $n=m=2$ and $n=m=10$ and for subscript $\star$ corresponding to $x$ and $y$ components; we also compare the reduced solutions with the HF globally defined solutions  $\texttt{P}_{\rm pu}[u^{\rm hf}_{\star}]$. We observe that  the choice $n=m=2$ enables qualitatively accurate approximations of the vertical displacement (cf. \ref{fig:vis_slices}(c)), but extremely inaccurate approximations of the horizontal  displacement (cf. \Cref{fig:vis_slices}(b)), while the choice $n=m=10$ leads to accurate predictions for both horizontal and vertical displacements.

\subsection{Hyper-reduction of the port-to-bubble maps}

\begin{figure}[H]
	\centering
	\subfloat[]{
		\begin{tikzpicture}
			\begin{loglogaxis}[
				xmode=linear,
				ymode=log,
				grid=both,
				minor grid style={gray!25},
				major grid style={gray!25},
				title = {},
				xlabel={$m$},
				ylabel={$E_{\rm avg}$},
				width=.35\textwidth,
				xtick={2,4,6,8,10,12,14,16},
				legend style={at={(1,0.9)},anchor=east},
				ymin=5*1e-5,
				ymax=0.25,
				]
				]
				\addplot[ 
				red,
				mark=triangle*,
				mark options={ color = red},
				]  %
				table{dat/max_errors/onlyEQ/avg_errHFQ.dat};
				\addlegendentry{\tiny{HFQ}};
				\addplot[ 
				black,
				mark=square,
				mark options={ color = black},
				]  %
				table{dat/max_errors/onlyEQ/avg_err_tol1minus8.dat};
				\addlegendentry{\tiny{$tol_{\rm eq}=10^{-8}$}};
				\addplot[ 
				black,
				mark=triangle,
				mark options={ color = black},
				]  %
				table{dat/max_errors/onlyEQ/avg_err_tol1minus10.dat};
				\addlegendentry{\tiny{$tol_{\rm eq}=10^{-10}$}};
				\addplot[ 
				black,
				mark=o,
				mark options={ color = black},
				]  %
				table{dat/max_errors/onlyEQ/avg_err_tol1minus14.dat};
				\addlegendentry{\tiny{$tol_{\rm eq}=10^{-14}$}};
			\end{loglogaxis}
		\end{tikzpicture}
	}
	\subfloat[]{
		\begin{tikzpicture}
			\begin{loglogaxis}[
				xmode=linear,
				ymode=log,
				grid=both,
				minor grid style={gray!25},
				major grid style={gray!25},
				title = {},
				xlabel={$m$},
				ylabel={$E_{\rm max}$},
				width=.35\textwidth,
				xtick={2,4,6,8,10,12,14,16},
				legend style={at={(1,0.9)},anchor=east},
				ymin=5*1e-5,
				ymax=0.25,
				]
				]
				\addplot[ 
				red,
				mark=triangle*,
				mark options={ color = red},
				]  %
				table{dat/max_errors/onlyEQ/errHFQ.dat};
				\addlegendentry{\tiny{HFQ}};
				\addplot[ 
				black,
				mark=square,
				mark options={ color = black},
				]  %
				table{dat/max_errors/onlyEQ/err_tol1minus8.dat};
				\addlegendentry{\tiny{$tol_{\rm eq}=10^{-8}$}};
				\addplot[ 
				black,
				mark=triangle,
				mark options={ color = black},
				]  %
				table{dat/max_errors/onlyEQ/err_tol1minus10.dat};
				\addlegendentry{\tiny{$tol_{\rm eq}=10^{-10}$}};
				\addplot[ 
				black,
				mark=o,
				mark options={ color = black},
				]  %
				table{dat/max_errors/onlyEQ/err_tol1minus14.dat};
				\addlegendentry{\tiny{$tol_{\rm eq}=10^{-14}$}};
			\end{loglogaxis}
		\end{tikzpicture}
	}
	\quad
	\subfloat[${\Omega}^{\rm a}_{\rm int}$]{
		\begin{tikzpicture}
			\begin{semilogyaxis}[
				xmode=linear,
				ymode=linear,
				grid=both,
				minor grid style={gray!25},
				major grid style={gray!25},
				title = {},
				ymin= 0, ymax = 23,
				xlabel={$m$},
				ylabel={$\%$ \text{elements}},
				yticklabel={\pgfmathparse{\tick}\pgfmathprintnumber{\pgfmathresult}\%},
				width=0.35\textwidth,
				xtick={2,4,6,8,10,12,14,16},
				legend style={at={(0.04,0.9)},anchor=west},
				]
				]	
				\addplot[ 
				black,
				mark=square,
				mark options={ color = black},
				]  %
				table{dat/perc_al3.dat};
				\addlegendentry{\tiny{$tol_{\rm eq}=10^{-8}$}};
				\addplot[ 
				black,
				mark=triangle,
				mark options={ color = black},
				]  %
				table{dat/perc_al1.dat};
				\addlegendentry{\tiny{$tol_{\rm eq}=10^{-10}$}};
				\addplot[ 
				black,
				mark=o,
				mark options={ color = black},
				]  %
				table{dat/perc_al2.dat};
				\addlegendentry{\tiny{$tol_{\rm eq}=10^{-14}$}};
			\end{semilogyaxis}
		\end{tikzpicture}
	}
	\subfloat[${\Omega}^{\rm a}_{\rm ext}$]{
		\begin{tikzpicture}
			\begin{semilogyaxis}[
				xmode=linear,
				ymode=linear,
				grid=both,
				minor grid style={gray!25},
				major grid style={gray!25},
				title = {},
				ymin= 0, ymax = 23,
				xlabel={$m$},
				yticklabel={\pgfmathparse{\tick}\pgfmathprintnumber{\pgfmathresult}\%},
				width=0.35\textwidth,
				xtick={2,4,6,8,10,12,14,16},
				legend style={at={(-0.001,0.9)},anchor=west},
				]
				]
				\addplot[ 
				black,
				mark=square,
				mark options={ color = black},
				]  %
				table{dat/perc_ext3.dat};
				\addlegendentry{\tiny{$tol_{\rm eq}=10^{-8}$}};
				\addplot[ 
				black,
				mark=triangle,
				mark options={ color = black},
				]  %
				table{dat/perc_ext1.dat};
				\addlegendentry{\tiny{$tol_{\rm eq}=10^{-10}$}};
				\addplot[ 
				black,
				mark=o,
				mark options={ color = black},
				]  %
				table{dat/perc_ext2.dat};
				\addlegendentry{\tiny{$tol_{\rm eq}=10^{-14}$}};
			\end{semilogyaxis}
		\end{tikzpicture}
		\label{fig:elements_ext}
	}
	\caption{hyper-reduction of the port-to-bubble maps 
		for several tolerances $tol_{\rm eq}$ and port space sizes $m$, with $n=m$.
		Behavior of the (a)
		average, (b) max out-of-sample prediction.
		(c)-(d)
		percentage of sampled elements in ${\Omega}^{\rm a}_{\rm int}$ and 
		${\Omega}^{\rm a}_{\rm ext}$.
	}
	\label{fig:hyper_port2bubble}
\end{figure}
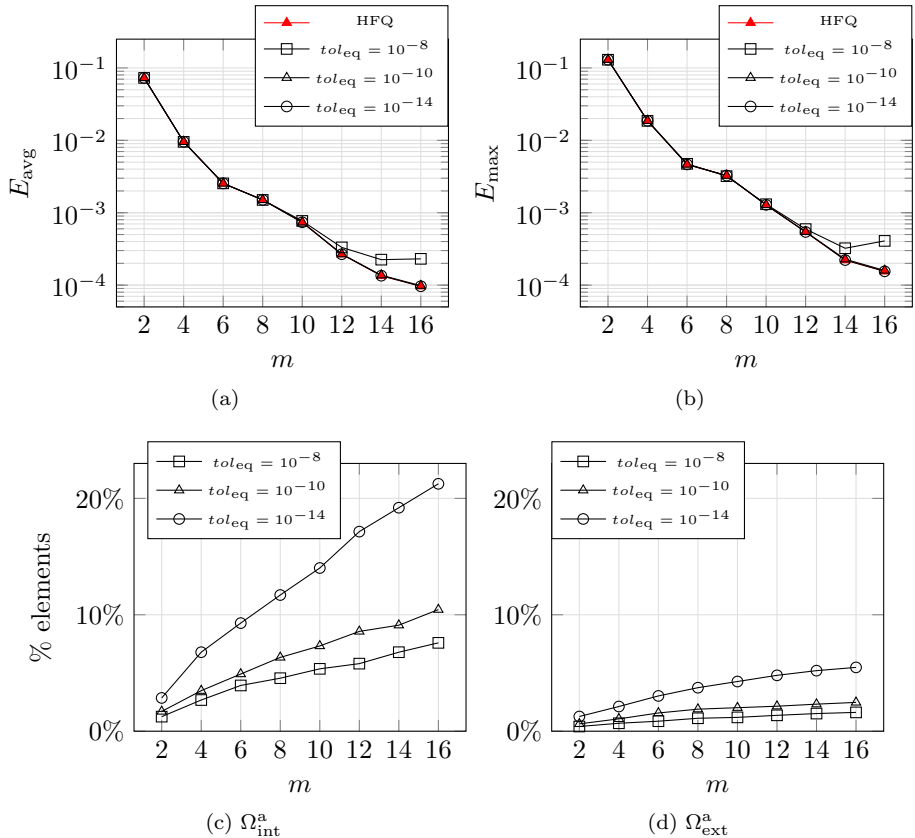

\Cref{fig:hyper_port2bubble}  investigates performance of the EQ rule for different tolerances $tol_{\rm eq}$ (cf. \ref{sec:app_hyper_port2bubble}):
Figures \ref{fig:hyper_port2bubble}(a) and  \ref{fig:hyper_port2bubble}(b) show the behavior of the out-of-sample relative error compared to the OS2 ROM with HF quadrature (dubbed HFQ); in \Cref{fig:hyper_port2bubble}(a) we depict $E_{\rm avg}$, in \Cref{fig:hyper_port2bubble}(b) we depict $E_{\rm max}$.
Figures \ref{fig:hyper_port2bubble}(c) and \ref{fig:hyper_port2bubble}(d) show the percentage of sampled elements as a function of $m$, for the two archetype components and for several tolerances.
We observe that for $tol_{\rm eq} \leq 10^{-10}$ the hyper-reduced OS2 ROM is   as accurate as the OS2 ROM with HF quadrature for all values of $m$ considered.
We further observe that the percentage of sampled elements is between three and  five times larger in the internal component --- since $N_{\rm ext}^{\rm e} \approx 3.5 N_{\rm int}^{\rm e}$, we have that the absolute number of sampled elements is nearly the same for the two components.

\subsection{Hyper-reduction of the objective function}

\begin{figure}[H]
	\centering
	\subfloat[$\Omega_{\rm int}^{\rm a}$]{
		\begin{tikzpicture}
			\begin{loglogaxis}[
				xmode=linear,
				ymode=log,
				grid=both,
				xmin=1,
				xmax=8,			
				ymin=0.0000001,
				ymax=0.001,
				minor grid style={gray!25},
				major grid style={gray!25},
				title = {},
				xlabel={$m$},
				ylabel={$E^{\infty}_{\rm avg,eim}$},
				width=0.4\textwidth,
				xtick={2,4,6,8,10,12,14,16},
				legend style={at={(0.99,0.5)},anchor=east},
				]
				]
				\addplot[ 
				black,
				mark=triangle*,
				mark options={ color = black},
				]  %
				table{dat/int_err1.dat};
			\end{loglogaxis}
		\end{tikzpicture}
	}
	\subfloat[$\Omega_{\rm ext}^{\rm a}$]{
		\begin{tikzpicture}
			\begin{loglogaxis}[
				xmode=linear,
				ymode=log,
				grid=both,
				xmin=1,
				xmax=8,			
				ymin=0.0000001,
				ymax=0.001,
				minor grid style={gray!25},
				major grid style={gray!25},
				title = {},
				xlabel={$m$},
				ylabel={$E^{\infty}_{\rm avg,eim}$},
				width=0.4\textwidth,
				xtick={2,4,6,8,10,12,14,16},
				legend style={at={(0.99,0.5)},anchor=east},
				]
				]
				\addplot[ 
				black,
				mark=triangle*,
				mark options={ color = black},
				]  %
				table{dat/int_err2.dat};
			\end{loglogaxis}
		\end{tikzpicture}
	}	
	\caption{
		application of the EIM procedure for vector-valued fields  (cf. \Cref{alg:EIM}).
		(a)-(b) 
		behavior of the in-sample $L^{\infty}$ approximation error $E_{\rm avg,eim}^{\infty}$  for the internal and the  external component.}
	\label{fig:int_err}
\end{figure}
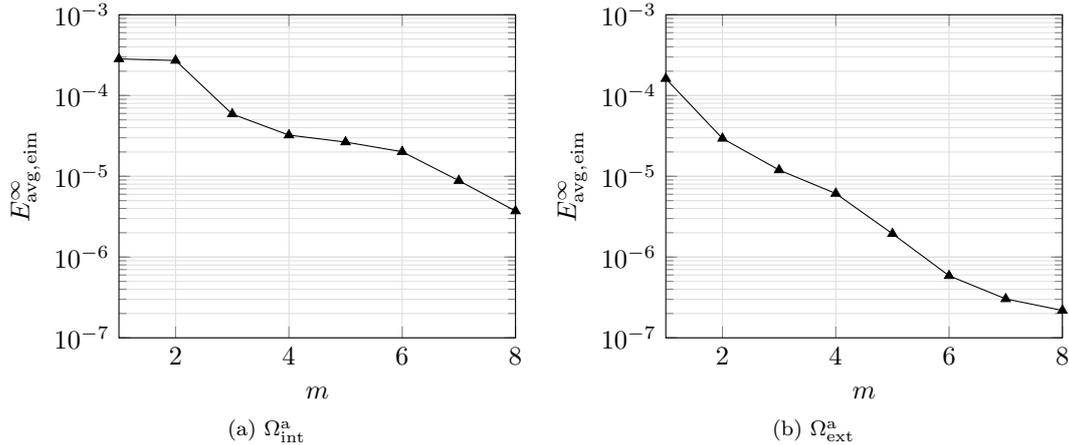

In  \Cref{fig:int_err}, we show the behavior of the $L^{\infty}$  error
$$
E_{\rm avg, eim}^{\infty}(\ell,m) :=\frac{1}{n_{\rm train,\ell}}
\sum_{k=1}^{n_{\rm train,\ell}}
\Big\|u_{\ell,k}^{\rm p}-\mathcal{I}_{\ell,m}[u^p_{k}]\Big\|_{\infty}
$$
where 
$\{ u_{\ell,k}^{\rm p} \}_{k=1}^{   n_{\rm train,\ell} }$  are the port fields associated with the $\ell$-th component and employed to generate the port basis (cf.   \Cref{alg:data_compression}). 
We observe near-exponential convergence of the $L^{\infty}$  error for both components; 
interestingly, the interpolation error for the internal component is one order of magnitude larger than the error for the external component.

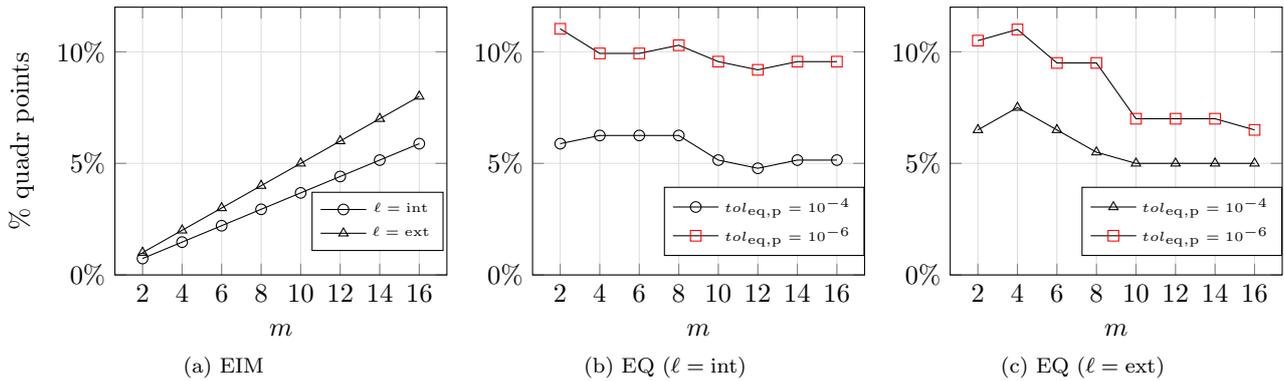
\begin{figure}[H]
	\centering
	\subfloat[EIM]{
		\begin{tikzpicture}
			\begin{axis}[
				xmode=linear,
				ymode=linear,
				grid=both,
				ymin=0,
				ymax=12,
				minor grid style={gray!25},
				major grid style={gray!25},
				title = {},
				xlabel={$m$},
				width=0.35\textwidth,
				xtick={2,4,6,8,10,12,14,16},
				legend style={at={(0.99,0.2)},anchor=east},
				ylabel={\% \text{quadr points}},
				yticklabel={\pgfmathparse{\tick}\pgfmathprintnumber{\pgfmathresult}\%},
				ymax=12,
				]
				]
				\addplot[ 
				black,
				mark=o,
				mark options={ color = black},
				]  %
				table{dat/EQvsEIM/percEIM_int.dat};
				\addlegendentry{\tiny{$\ell=\rm int$}};
				\addplot[ 
				black,
				mark=triangle,
				mark options={ color = black},
				]  %
				table{dat/EQvsEIM/percEIM_ext.dat};
				\addlegendentry{\tiny{$\ell=\rm ext$}};
			\end{axis}
		\end{tikzpicture}
	}
	\subfloat[EQ ($\ell=\rm int$)]{
		\begin{tikzpicture}
			\begin{axis}[
				xmode=linear,
				ymode=linear,
				grid=both,
				ymin=0,
				ymax=12,
				minor grid style={gray!25},
				major grid style={gray!25},
				title = {},
				xlabel={$m$},
				width=0.35\textwidth,
				xtick={2,4,6,8,10,12,14,16},
				legend style={at={(1.0,0.2)},anchor=east},
				ylabel={},
				ylabel style={
				},
				yticklabel={\pgfmathparse{\tick}\pgfmathprintnumber{\pgfmathresult}\%},
				]
				]
				\addplot[ 
				black,
				mark=o,
				mark options={ color = black},
				]  %
				table{dat/EQvsEIM/percEQobj_int_tolminus4.dat};
				\addlegendentry{\tiny{$tol_{\rm eq,p}=10^{-4}$}};
				\addplot[ 
				black,
				mark=square,
				mark options={ color = red},
				]  %
				table{dat/EQvsEIM/percEQobj_int_tolminus6.dat};
				\addlegendentry{\tiny{$tol_{\rm eq,p}=10^{-6}$}};
			\end{axis}
		\end{tikzpicture}
	}
	\subfloat[EQ ($\ell=\rm ext$)]{
		\begin{tikzpicture}
			\begin{axis}[
				xmode=linear,
				ymode=linear,
				grid=both,
				minor grid style={gray!25},
				major grid style={gray!25},
				ymin=0,
				ymax=12,
				title = {},
				xlabel={$m$},
				width=0.35\textwidth,
				xtick={2,4,6,8,10,12,14,16},
				legend style={at={(1.0,0.2)},anchor=east},
				ylabel={},
				ylabel style={
				},
				yticklabel={\pgfmathparse{\tick}\pgfmathprintnumber{\pgfmathresult}\%},
				]
				]
				\addplot[ 
				black,
				mark=triangle,
				mark options={ color = black},
				]  %
				table{dat/EQvsEIM/percEQobj_ext_tolminus4.dat};
				\addlegendentry{\tiny{$tol_{\rm  eq,p}=10^{-4}$}};
				\addplot[ 
				black,
				mark=square,
				mark options={ color = red},
				]  %
				table{dat/EQvsEIM/percEQobj_ext_tolminus6.dat};
				\addlegendentry{\tiny{$tol_{\rm eq,p}=10^{-6}$}};
			\end{axis}
		\end{tikzpicture}
	}
	\caption{
		hyper-reduction of the objective function for internal and external archetype components, with respect to $m$,  with $n=m$.
		(a) percentage of sampled quadrature points based on EIM.
		(b)-(c) 
		percentage of sampled quadrature points based on the EQ procedure, for two tolerances $tol_{\rm eq,p}$.}
	\label{fig:EQ_EIM_obj}
\end{figure}

In  \Cref{fig:EQ_EIM_obj}, we report the percentage of sampled quadrature points by the two hyper-reduction procedures. By construction, EIM selects $m_{\rm p,eq}=m$ points; on the other hand, the number of points selected by the EQ procedure of   \cref{sec:EQ_bnd} weakly depends on the size $m$ of the port basis.

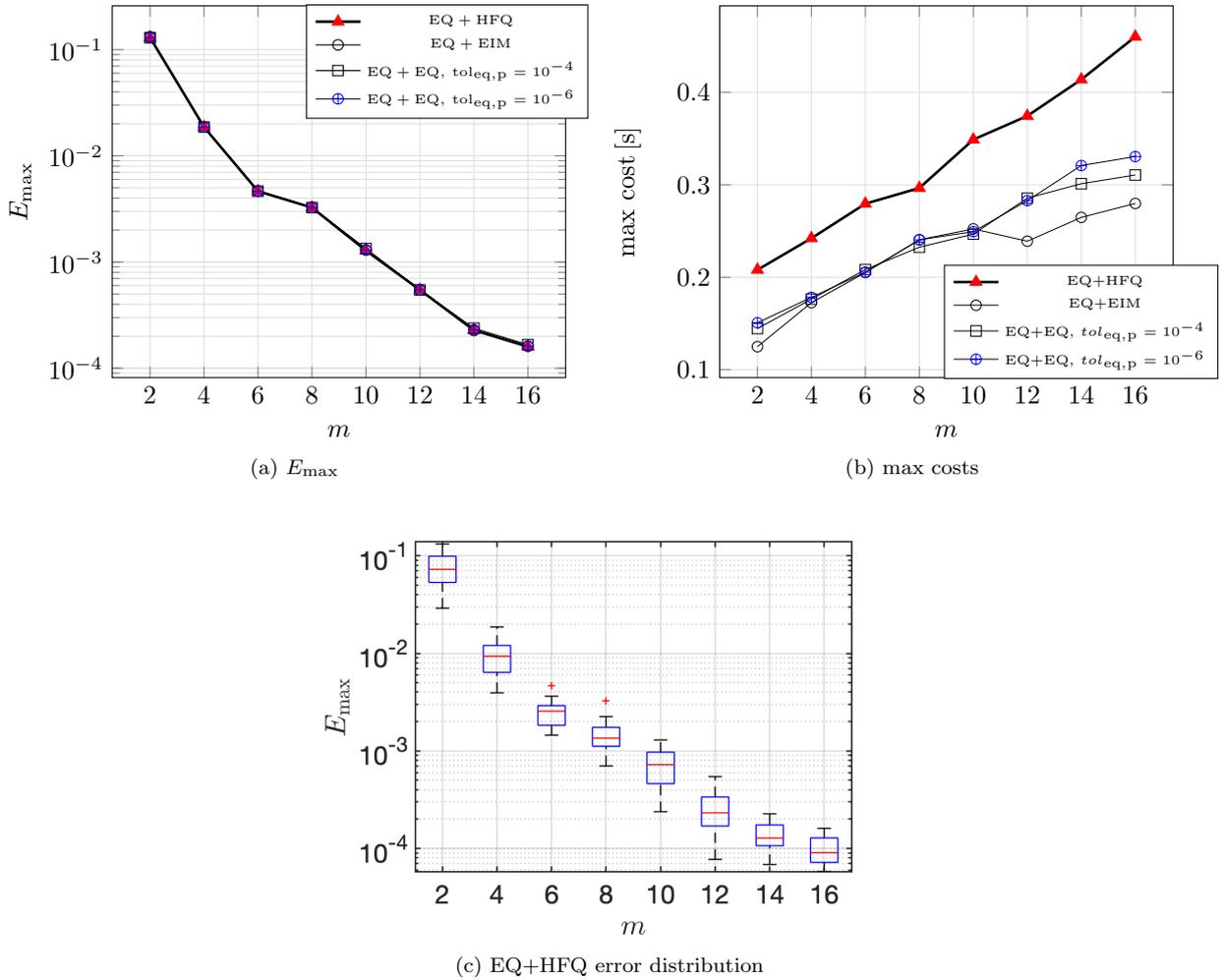
\begin{figure}[H]
	\centering
	\subfloat[$E_{\rm max}$]{
		\begin{tikzpicture}
			\begin{loglogaxis}[
				xmode=linear,
				ymode=log,
				grid=both,
				minor grid style={gray!25},
				major grid style={gray!25},
				title = {},
				xlabel={$m$},
				ylabel={$E_{\rm max}$},
				width=0.45\textwidth,
				xtick={2,4,6,8,10,12,14,16},
				legend style={at={(1.05,0.85)},anchor=east},
				]
				]
				
				\addplot[ 
				black,
				line width=0.35mm,
				mark=triangle*,
				mark options={ color = red},
				]  %
				table{dat/max_errors/EQ2_EQEIM/errEQ_HFQ.dat};
				\addlegendentry{\tiny{$\rm EQ+HFQ$}};

				\addplot[ 
				black,
				mark=o,
				mark options={ color = black},
				]  %
				table{dat/max_errors/EQ2_EQEIM/errEIM.dat};
				\addlegendentry{\tiny{$\rm EQ+EIM$}};
				
				\addplot[ 
				black,
				mark=square,
				mark options={ color = black},
				]  %
				table{dat/max_errors/EQ2_EQEIM/errGal1eminus4.dat};
				\addlegendentry{\tiny{$\rm EQ+EQ, \, tol_{\rm eq,p}=10^{-4}$}};

				
				\addplot[ 
				black,
				mark=oplus,
				mark options={ color = blue},
				]  %
				table{dat/max_errors/EQ2_EQEIM/errGal1eminus6.dat};
				\addlegendentry{\tiny{$\rm EQ+EQ, \, tol_{\rm eq,p}=10^{-6}$}};
			\end{loglogaxis}
		\end{tikzpicture}
	}
	\subfloat[max costs]{
		\begin{tikzpicture}
			\begin{loglogaxis}[
				xmode=linear,
				ymode=linear,
				grid=both,
				minor grid style={gray!25},
				major grid style={gray!25},
				title = {},
				xlabel={$m$},
				ylabel={$ \max \, \rm cost \, [s]$},
				width=0.45\textwidth,
				xtick={2,4,6,8,10,12,14,16},
				legend style={at={(1.1,0.15)},anchor=east},
				]
				]	
				\addplot[ 
				black,
				line width=0.35mm,
				mark=triangle*,
				mark options={ color = red},
				]  %
				table{dat/EQvsEIM/EQ_HFQ_cost.dat};
				\addlegendentry{\tiny{EQ+HFQ}};
				\addplot[ 
				black,
				line width=0.1mm,
				mark=o,
				mark options={ color = black},
				]  %
				table{dat/EQvsEIM/EIMobj_cost.dat};
				\addlegendentry{\tiny{EQ+EIM}};
				\addplot[ 
				black,
				line width=0.1mm,
				mark=square,
				mark options={ color = black},
				]  %
				table{dat/EQvsEIM/EQ2_tol1eminus4_cost.dat};
				\addlegendentry{\tiny{EQ+EQ, $tol_{\rm eq,p}=10^{-4}$}};
				\addplot[ 
				black,
				line width=0.1mm,
				mark=oplus,
				mark options={ color =blue},
				]  %
				table{dat/EQvsEIM/EQ2_tol1eminus6_cost.dat};
				\addlegendentry{\tiny{EQ+EQ, $tol_{\rm eq,p}=10^{-6}$}};
			\end{loglogaxis}
		\end{tikzpicture}
	}
	\quad 
	\subfloat[EQ+HFQ error distribution]{
		\includegraphics[width=0.46\textwidth]{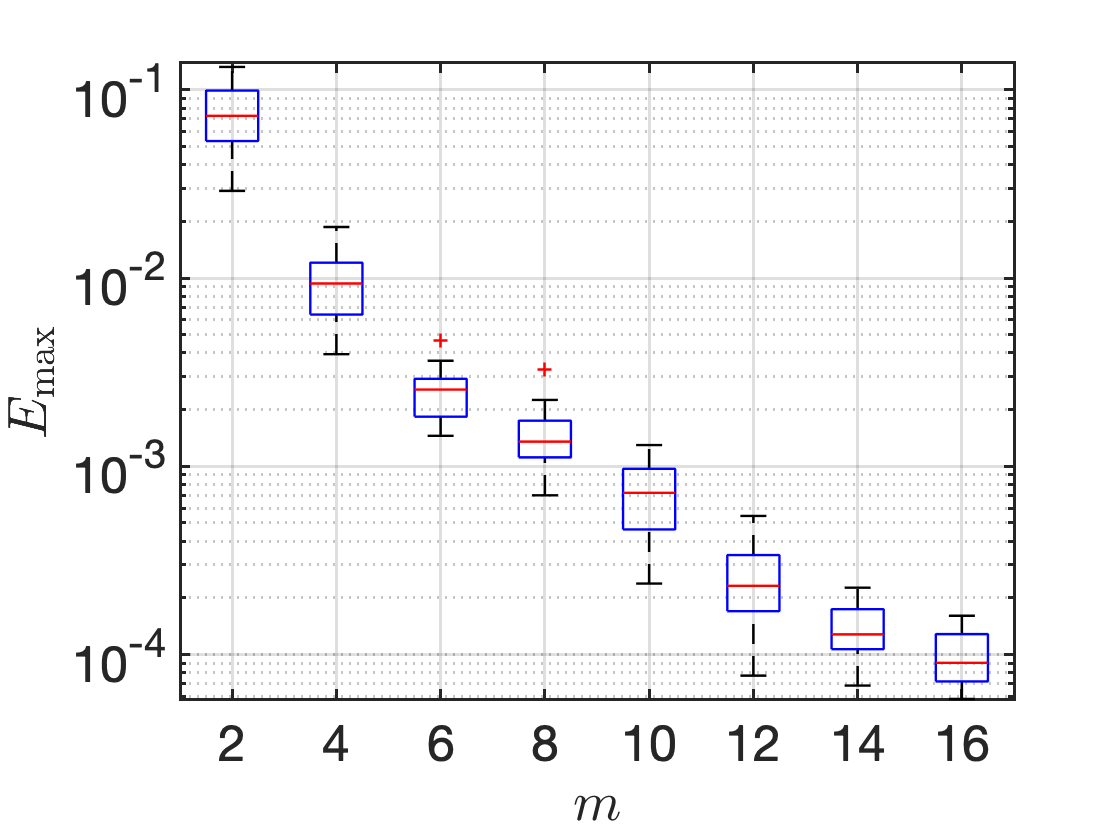}}
	\caption{
		hyper-reduction of the objective function based on EIM and EQ.
		(a) worst out-of-sample performance of the hyper-reduced OS2 ROM for several choices of $m$, with  $n=m$.
		(b) maximum computational cost over the test set.
		(c) Out-of-sample error distributions for the EQ+HFQ case.
		Results are based on the EQ tolerance 
		$tol_{\rm eq}=10^{-10}$ for the local problems and the tolerances $tol_{\rm eq,p}=10^{-4}$ and 
		$tol_{\rm eq,p}=10^{-6}$ for the objective function (for EQ+EQ).}
	\label{fig:hyperobj_err_cost}
\end{figure}

In  \Cref{fig:hyperobj_err_cost}, we investigate the performance of the fully hyper-reduced ROM:
\Cref{fig:hyperobj_err_cost}(a) shows the behavior of the prediction error \eqref{eq:Eavg}, while 
\Cref{fig:hyperobj_err_cost}(b) shows the behavior of the maximum wall-clock time over the test set. 
We observe that the speed-up due to hyper-reduction of the   objective function is of the order $1.5$ for all choices of $m$; on the other hand, performance of  the two considered hyper-reduction strategies is comparable for all tests. 
In \Cref{fig:hyperobj_err_cost}(c) the out-of-sample error distributions are depicted in the case without hyper-reduction on the objective function for different values of $m$.

\begin{figure}[H]
	\centering
	\subfloat[$m=8$]{
		\begin{tikzpicture}
			\begin{axis}[
				xmode=linear,
				ymode=linear,
				grid=both,
				minor grid style={gray!25},
				major grid style={gray!25},
				title = {},
				xlabel={$N_{\rm dd}$},
				ylabel={speed-up},
				width=0.45\textwidth,
				xtick={3,4,5,6,7,8},
				legend style={at={(0.93,0.3)},anchor=east},
				ymax=30,
				]
				]	
				\addplot[ 
				black,
				line width=0.35mm,
				mark=triangle*,
				mark options={ color = red},
				]  %
				table{dat/speedup/m8/speedup_HFQ2.dat};
				\addlegendentry{\tiny{HFQ+HFQ}};
				\addplot[ 
				black,
				mark=square,
				mark options={ color = black},
				]  %
				table{dat/speedup/m8/speedup_EQ2_tol1eminus4.dat};
				\addlegendentry{\tiny{EQ+EQ, $tol_{\rm eq,p}=10^{-4}$}};
				\addplot[ 
				black,
				mark=oplus,
				mark options={ color = blue},
				]  %
				table{dat/speedup/m8/speedup_EQ2_tol1eminus6.dat};
				\addlegendentry{\tiny{EQ+EQ, $tol_{\rm eq,p}=10^{-6}$}};
				\addplot[ 
				black,
				mark=o,
				mark options={ color = black},
				]  %
				table{dat/speedup/m8/speedup_EIM.dat};
				\addlegendentry{\tiny{EQ+EIM}};
			\end{axis}
		\end{tikzpicture}
		\label{fig:speedup_m8}
	}
	\subfloat[$m=16$]{
		\begin{tikzpicture}
			\begin{axis}[
				xmode=linear,
				ymode=linear,
				grid=both,
				minor grid style={gray!25},
				major grid style={gray!25},
				title = {},
				xlabel={$N_{\rm dd}$},
				ylabel={speed-up},
				width=0.45\textwidth,
				xtick={3,4,5,6,7,8},
				legend style={at={(0.93,0.3)},anchor=east},
				ymax=30,
				]
				]	
				\addplot[ 
				black,
				line width=0.35mm,
				mark=triangle*,
				mark options={ color = red},
				]  %
				table{dat/speedup/m16/speedup_HFQ2.dat};
				\addlegendentry{\tiny{HFQ+HFQ}};
				\addplot[ 
				black,
				mark=square,
				mark options={ color = black},
				]  %
				table{dat/speedup/m16/speedup_EQ2_tol1eminus4.dat};
				\addlegendentry{\tiny{EQ+EQ, $tol_{\rm eq,p}=10^{-4}$}};
				\addplot[ 
				black,
				mark=oplus,
				mark options={ color = blue},
				]  %
				table{dat/speedup/m16/speedup_EQ2_tol1eminus6.dat};
				\addlegendentry{\tiny{EQ+EQ, $tol_{\rm eq,p}=10^{-6}$}};
				\addplot[ 
				black,
				mark=o,
				mark options={ color = black},
				]  %
				table{dat/speedup/m16/speedup_EIM.dat};
				\addlegendentry{\tiny{EQ+EIM}};
			\end{axis}
		\end{tikzpicture}
		\label{fig:speedup_m16}
	}
	\caption{Speed-up  of the   OS2 ROMs  with respect to the HF  monolithic solver for several values of the number of subdomains. 
		(a)  performance for $m=8$; 
		(b) performance for $m=16$. 
		EQ tolerance for the port-to-bubble maps is set equal to $tol_{\rm eq}=10^{-10}$.}
	\label{fig:speedup}
\end{figure}
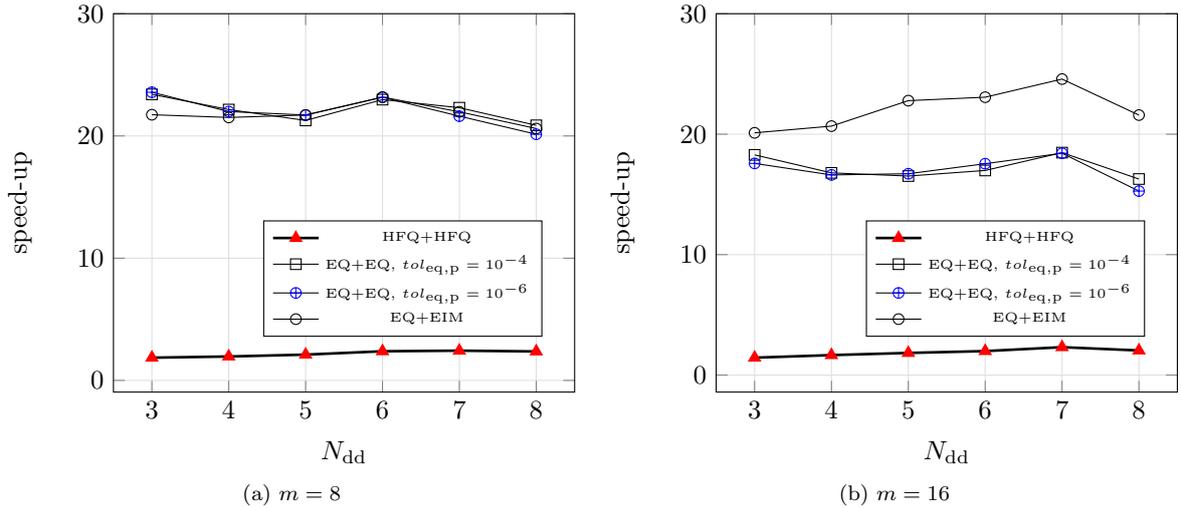

In Figure \ref{fig:speedup} we show the speed-up factor  of the hyper-reduced OS2 solvers   with respect to a representative  monolithic HF solver of comparable accuracy    for different numbers of subdomains.
The  monolithic P2 FE solver  runs in approximately\footnote{Computational times  are based 
	on an average over $5$ tests for each number of subdomains; the computational grid has  $17177$ FE nodes for $N_{\rm dd}=2$ and it has $38637$ nodes for $N_{\rm dd}=8$.}  $2.7806 \, [{\rm s}]$ for $N_{\rm dd}=2$ and in $9.9971 \, [{\rm s}]$ for $N_{\rm dd}=8$; the CB HF solver \eqref{eq:OS2_original} that is used to generate  training and test data is roughly a factor three slower than the corresponding monolithic solver.
We define the speed-up factor as:
\begin{equation*}
	\rm{speed}{\text -}\rm{up}(N_{\rm dd}):=\frac{t_{\rm hf}(N_{\rm dd})}{t_{\rm OS2}(N_{\rm dd})}
\end{equation*}
where $t_{\rm hf}$ is the estimated execution time of the  monolithic HF solver averaged over $5$ tests and  $t_{\rm OS2}$ is the execution time associated with the CB ROM, averaged over the same $5$ configurations, for  $N_{\rm dd} \in \{3,\ldots, 8\}$.
We perform  hyper-reduction 
of the port-to-bubble maps using the    tolerance $tol_{\rm eq}=10^{-10}$ and 
we consider 
the   tolerances $tol_{\rm eq,p}=10^{-4}$ and 
$tol_{\rm eq,p}=10^{-6}$ for the hyper-reduction of the objective function (for the EQ+EQ case).

We  observe that the speed-up factors depicted in Figure \ref{fig:speedup} depend weakly  on the number of subdomains. 
The EIM method leads to  slightly larger speed-ups than the EQ method for  $m=n=16$  (cf. Figure \ref{fig:speedup}(b)), 
while performance is comparable for the case 
$m=n=8$ (Figure \ref{fig:speedup}(a)).
We envision that more effective implementations of \Cref{alg:objective} 
---
which  rely on parallelization of the port-to-bubble loop at Lines 4-7 and on pointwise EQ hyper-reduction of the port-to-bubble maps,  as opposed to element-wise EQ 
---
will lead to significantly larger speed-ups.

\subsection{Optimization strategy: comparison between Gauss-Newton, quasi-Newton and  overlapping Schwarz}

We compare the performance of the Gauss-Newton method and  the quasi-Newton method  discussed in section \ref{sec:minimization} for various choices of $m$ and $n=m$; to provide a concrete reference, we also consider the  multiplicative overlapping Schwarz method with Dirichlet interface conditions.
More precisely, we implement the iterative procedure described in \Cref{alg:os_poor}:
note that the OS method simply requires the solution to a sequence of problems on the subdomains with information propagating through the boundary conditions;
since the discretization is not conforming across components, we should define the $i$-th port mode using projection (cf. Line 5, \Cref{alg:os_poor}). Note that at step $i$ of the for loop at Lines $4-7$ we use the values of $\widehat{\boldsymbol{\beta}}_{i},\ldots,\widehat{\boldsymbol{\beta}}_{N_{\rm dd}}$ at the previous  iteration and the values 
$\widehat{\boldsymbol{\beta}}_{1},\ldots,\widehat{\boldsymbol{\beta}}_{i-1}$ at the current iteration: the for loop is thus not parallelizable.
We set $tol=10^{-6}$ in \Cref{alg:objective}  (cf. Line 11) and we consider the same termination criterion for the quasi-Newton solver and the OS solver. 
In this test, we perform hyper-reduction at the local level (EQ tolerance $10^{-10}$), but we do not hyper-reduce the objective function.

\begin{algorithm}[H]                      
	\caption{Overlapping Schwarz method.}     
	\label{alg:os_poor}     
	
	\small
	\begin{flushleft}
		\emph{Inputs:}  
		$\boldsymbol{\alpha}^{(0)} = [\boldsymbol{\alpha}_1^{(0)} , \ldots,\boldsymbol{\alpha}_{N_{\rm dd}}^{(0)} ]$,
		$\boldsymbol{\beta}^{(0)}  = [\boldsymbol{\beta}_1^{(0)} , \ldots,\boldsymbol{\beta}_{N_{\rm dd}}^{(0)}]$ initial conditions (cf. Eq. \eqref{eq:IC_GNM}),
		$tol>0, \texttt{maxit}$.
		\smallskip
		
		\emph{Outputs:} 
		$\widehat{\boldsymbol{\beta}}$ port coefficients,
		$\widehat{\boldsymbol{\alpha}}= \widehat{\texttt{F}}^{\rm eq}(  \widehat{\boldsymbol{\beta}}  )$  bubble coefficients.
	\end{flushleft}                      
	
	\normalsize 
	
	\begin{algorithmic}[1]
		\State Set $\widehat{\boldsymbol{\beta}}^{(0)} = \boldsymbol{\beta}^{(0)}$ and 
		$\widehat{\boldsymbol{\alpha}}=\boldsymbol{\alpha}^{(0)}$.
		\medskip
		
		\For {$k=1, \ldots, \texttt{maxit}$ }
		
		\State
		Initialize $\widehat{\boldsymbol{\alpha}}^{(k)}
		=\widehat{\boldsymbol{\alpha}}^{(k-1)}$ and
		$\widehat{\boldsymbol{\beta}}^{(k)}
		=\widehat{\boldsymbol{\beta}}^{(k-1)}$.
		\medskip

		\For {$i=1, \ldots, N_{\rm dd}$ }
		
		\State
		Update $\widehat{\boldsymbol{\beta}}_i^{(k)}
		\in {\rm arg}
		\min_{    {\boldsymbol{\beta}}\in \mathbb{R}^m   }
		\sum_{j \in {\rm Neigh}_i} \;
		\|   W_i^{\rm p}  {\boldsymbol{\beta}} 
		\, - \, 
		Z_j^{\rm b} 
		\widehat{\texttt{F}}_j^{\rm eq}( \widehat{\boldsymbol{\beta}}_j^{(k)} )
		\, 	- \, 
		W_j^{\rm p}  \widehat{\boldsymbol{\beta}}_j^{(k)}
		\|_{L^2(\Gamma_{i,j})}^2$.
		
		\State
		Update $\widehat{\boldsymbol{\alpha}}_i^{(k)}
		=
		\widehat{\texttt{F}}_i^{\rm eq}( \widehat{\boldsymbol{\beta}}_i^{(k)} )$. 
		
		\EndFor
		\vspace{3pt}
		
		\If {$\| \widehat{\boldsymbol{\beta}}^{(k)} - \widehat{\boldsymbol{\beta}}^{(k-1)}   \|_2 < tol \| \widehat{\boldsymbol{\beta}}^{(k)} \|_2$}, \texttt{BREAK}
		
		\EndIf
		
		\EndFor
		
	\end{algorithmic}
\end{algorithm}

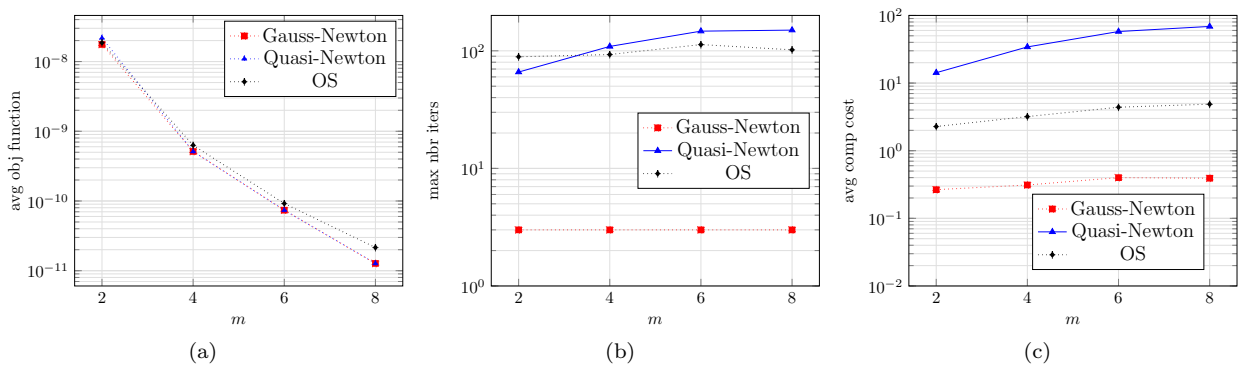
\begin{figure}[h]
	\centering
	\subfloat[]{
		\begin{tikzpicture}[scale=0.63]
			\begin{loglogaxis}[
				xmode=linear,
				ymode=log,
				grid=both,
				minor grid style={gray!25},
				major grid style={gray!25},
				title = {},
				xlabel={$m$},
				ylabel={avg obj function},
				xtick={2,4,6,8},
				]
				]
				
				\addplot[ 
				red,
				dotted,
				mark=square*,
				mark options={ color = red},
				]  %
				table{dat/QNtest/fstar_gn.dat};
				\addlegendentry{\large{Gauss-Newton}};		
				
				\addplot[ 
				blue,
				dotted,
				mark=triangle*,
				mark options={ color = blue},
				]  %
				table{dat/QNtest/fstar_qn.dat};
				\addlegendentry{\large{Quasi-Newton}};

				\addplot[ 
				black,
				dotted,
				mark=diamond*,
				mark options={ color = black},
				]  %
				table{dat/QNtest/fstar_OS.dat};
				\addlegendentry{\large{OS}};
				
			\end{loglogaxis}
		\end{tikzpicture}
	}
	\subfloat[]{
		\begin{tikzpicture}[scale=0.63]
			\begin{axis}[
				xmode=linear,
				ymode=log,
				grid=both,
				minor grid style={gray!25},
				major grid style={gray!25},
				title = {},
				ymin=1,
				ymax=200,
				xlabel={$m$},
				ylabel={max nbr iters},
				xtick={2,4,6,8},
				legend style={at={(0.97,0.5)},anchor=east},
				]
				]
				\addplot[ 
				red,
				dotted,
				mark=square*,
				mark options={ color = red},
				]  %
				table{dat/QNtest/iters_gn.dat};
				\addlegendentry{\large{Gauss-Newton}};
				\addplot[ 
				blue,
				mark=triangle*,
				mark options={ color = blue},
				]  %
				table{dat/QNtest/iters_qn.dat};
				\addlegendentry{\large{Quasi-Newton}};

				\addplot[ 
				black,
				dotted,
				mark=diamond*,
				mark options={ color = black},
				]  %
				table{dat/QNtest/iters_os.dat};
				\addlegendentry{\large{OS}};			
				
			\end{axis}
		\end{tikzpicture}
	}
	\subfloat[]{
		\begin{tikzpicture}[scale=0.63]
			\begin{axis}[
				xmode=linear,
				ymode=log,
				grid=both,
				minor grid style={gray!25},
				major grid style={gray!25},
				title = {},
				ymin=0.01,
				ymax=100,
				xlabel={$m$},
				ylabel={avg comp cost},
				xtick={2,4,6,8},
				legend style={at={(0.9,0.2)},anchor=east},
				]
				]
				\addplot[ 
				red,
				dotted,
				mark=square*,
				mark options={ color = red},
				]  %
				table{dat/QNtest/cost_gn.dat};
				\addlegendentry{\large{Gauss-Newton}};
				\addplot[ 
				blue,
				mark=triangle*,
				mark options={ color = blue},
				]  %
				table{dat/QNtest/cost_qn.dat};
				\addlegendentry{\large{Quasi-Newton}};

				\addplot[ 
				black,
				dotted,
				mark=diamond*,
				mark options={ color = black},
				]  %
				table{dat/QNtest/cost_os.dat};
				\addlegendentry{\large{OS}};			
				
			\end{axis}
		\end{tikzpicture}
	}
	\caption{comparison between OS2 with Gauss-Newton optimization and with quasi-Newton optimization, and multiplicative overlapping Schwarz methods.
		(a) average value of the objective function with respect to $m$ and for $n=m$.	
		(b) maximum number of iterations to meet the convergence criterion.
		(c) average  wall-clock cost with respect to $m$ and for $n=m$.	
	}
	\label{fig:gn_qn}
\end{figure}

\Cref{fig:gn_qn}(a) shows the behavior of the objective function  in  \eqref{eq:OS2_eq}  with respect to the ROB sizes over the test set, while 
\Cref{fig:gn_qn}(b) shows the number of  iterations required to meet the termination criterion: we observe that GNM requires many fewer iterations without any deterioration in accuracy. 
\Cref{fig:gn_qn}(c) shows the wall-clock average cost for the three methods:
even if GNM has a slightly larger per-iteration cost, we empirically find that OS2 with GNM is significantly more rapid than the other two approaches.
Furthermore, since the OS internal loop (cf. Lines 4-7 \Cref{alg:os_poor}) is not parallelizable as opposed to the corresponding loop of the OS2 solver (cf. Lines 4-7 \Cref{alg:objective}), we expect significantly larger computational gains if we resort to parallel computing.

\begin{figure}[h]
	\centering
	\subfloat[]{
		\begin{tikzpicture}[scale=0.63]
			\begin{loglogaxis}[
				xmode=linear,
				ymode=log,
				grid=both,
				minor grid style={gray!25},
				major grid style={gray!25},
				title = {},
				xlabel={$m$},
				ylabel={avg obj function},
				legend style={at={(0.04,0.2)},anchor=west},
				xtick={2,4,6,8},
				]
				]
				\addplot[ 
				red,
				dotted,
				mark=square*,
				mark options={ color = red},
				]  %
				table{dat/QNtest_IC0/fstar_gn.dat};
				\addlegendentry{\large{Gauss-Newton}};

				\addplot[ 
				blue,
				dotted,
				mark=triangle*,
				mark options={ color = blue},
				]  %
				table{dat/QNtest_IC0/fstar_qn.dat};
				\addlegendentry{\large{Quasi-Newton}};

				\addplot[ 
				black,
				dotted,
				mark=diamond*,
				mark options={ color = black},
				]  %
				table{dat/QNtest_IC0/fstar_OS.dat};
				\addlegendentry{\large{OS}};
				
			\end{loglogaxis}
		\end{tikzpicture}
	}
	\subfloat[]{
		\begin{tikzpicture}[scale=0.63]
			\begin{axis}[
				xmode=linear,
				ymode=log,
				grid=both,
				minor grid style={gray!25},
				major grid style={gray!25},
				title = {},
				ymin=1,
				ymax=200,
				xlabel={$m$},
				ylabel={max nbr iters},
				xtick={2,4,6,8},
				legend style={at={(0.97,0.5)},anchor=east},
				]
				]
				\addplot[ 
				red,
				dotted,
				mark=square*,
				mark options={ color = red},
				]  %
				table{dat/QNtest_IC0/iters_gn.dat};
				\addlegendentry{\large{Gauss-Newton}};
				\addplot[ 
				blue,
				mark=triangle*,
				mark options={ color = blue},
				]  %
				table{dat/QNtest_IC0/iters_qn.dat};
				\addlegendentry{\large{Quasi-Newton}};

				\addplot[ 
				black,
				dotted,
				mark=diamond*,
				mark options={ color = black},
				]  %
				table{dat/QNtest_IC0/iters_os.dat};
				\addlegendentry{\large{OS}};			
				
			\end{axis}
		\end{tikzpicture}
	}
	\subfloat[]{
		\begin{tikzpicture}[scale=0.63]
			\begin{axis}[
				xmode=linear,
				ymode=log,
				grid=both,
				minor grid style={gray!25},
				major grid style={gray!25},
				title = {},
				ymin=0.01,
				ymax=100,
				xlabel={$m$},
				ylabel={avg comp cost},
				xtick={2,4,6,8},
				legend style={at={(0.9,0.2)},anchor=east},
				]
				]
				\addplot[ 
				red,
				dotted,
				mark=square*,
				mark options={ color = red},
				]  %
				table{dat/QNtest_IC0/cost_gn.dat};
				\addlegendentry{\large{Gauss-Newton}};
				\addplot[ 
				blue,
				mark=triangle*,
				mark options={ color = blue},
				]  %
				table{dat/QNtest_IC0/cost_qn.dat};
				\addlegendentry{\large{Quasi-Newton}};

				\addplot[ 
				black,
				dotted,
				mark=diamond*,
				mark options={ color = black},
				]  %
				table{dat/QNtest_IC0/cost_os.dat};
				\addlegendentry{\large{OS}};			
				
			\end{axis}
		\end{tikzpicture}
	}
	\caption{comparison between Gauss-Newton, quasi-Newton methods, and multiplicative overlapping Schwarz methods with zero initial condition.
		(a) average value of the objective function with respect to $m$ and for $n=m$.	
		(b) maximum number of iterations to meet the convergence criterion.
		(c) average  wall-clock cost with respect to $m$ and for $n=m$.	
	}
	\label{fig:gn_qn_newic}
\end{figure}
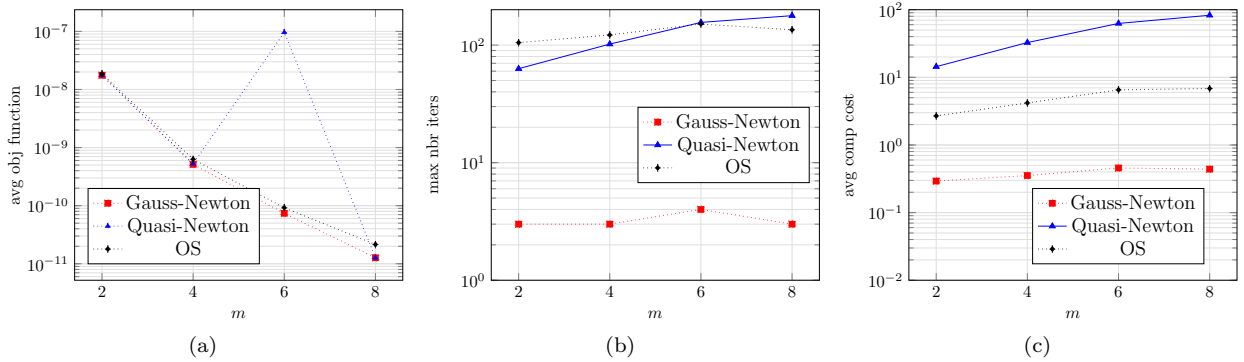

In \Cref{fig:gn_qn_newic}, we repeat the test of \Cref{fig:gn_qn} for the choice of the initial conditions
$\boldsymbol{\alpha}^{(0)} = \mathbf{0}$  and
$\boldsymbol{\beta}^{(0)} = \mathbf{0}$  in \Cref{alg:objective} and \Cref{alg:os_poor}.
We observe that OS and OS2 with GNM show similar performance with respect to all metrics, while OS2 with QN, instead of converging to the optimal solution, converges to a different local minimum for two configurations for $m=6$.

\section{Conclusions}
\label{sec:conclusions}
In this work we developed and numerically validated the one-shot overlapping Schwarz (OS2) approach to  component-based MOR of steady nonlinear PDEs.
The key features of the approach are 
(i)
a constrained  optimization statement that penalizes the jump at the components' interfaces  subject to the approximate satisfaction of the PDE  in each deployed (instantiated) component;
(ii)
the decomposition of the local solutions    into a \emph{port component} --- associated with the solution on interior boundaries (\emph{ports}) --- and a  \emph{bubble component}  that vanishes at \emph{ports}, to enable effective parallelization of the online solver.
Hyper-reduction of the local sub-problems and of the objective function is performed to reduce online assembly costs.
We illustrate the many elements of the formulation through the application to a two-dimensional nonlinear mechanics (Neo-Hookean) PDE model; for this problem, we are able to devise a CB-ROM  that reduces online costs by a factor $20$ compared to a standard monolithic FE model with less than $0.1\%$ prediction error, and without resorting to any parallelization of the online ROM solver.
We also observe that for the particular model problem considered in this paper the OS2 formulation provides acceptable results also for under-resolved ROBs.

We aim to extend our approach in several directions.
First, we wish to apply the OS2 method to more challenging problems in nonlinear mechanics, with particular emphasis on thermo-hydro-mechanical (THM) systems \cite{iollo2022adaptive}: towards this end, we should extend the OS2 formulation to unsteady PDEs and we should also devise specialized routines to deal with internal variables.
In this respect, we envision to combine our approach with the recently-developed OS method discussed in \cite{mota2022schwarz, de2022novel,barnett2022schwarz}.
Second, we wish to devise localized training techniques to avoid the solution to global HF problems at training stage:
in this regard, 
we aim to extend the approach in \cite{smetana2022localized} to unsteady PDEs with internal variables.
Third, we aim to 
combine data-fitted and projection-based ROMs in the OS2 framework: we envision that the 
successful combination of first-principle and data-fitted models might offer new  solutions for data assimilation (state estimation) applications, for a broad range of engineering tasks.
Fourth, we wish to investigate the possibility of generalizing our approach to non-overlapping decompositions: this requires to add in the objective function a term that penalizes the jump of the normal stresses; a deep investigation of the connection with discontinuous Galerkin methods should also be considered.

\section*{Acknowledgements}
The authors acknowledge the financial support of Andra (National Agency for Radioactive Waste Management) and thank Dr.
Marc Leconte and Dr. Antoine Pasteau (Andra) for fruitful discussions; 
The authors also thank Dr. Irina K. Tezaur (Sandia National Laboratories) for her pertinent and useful comments on this work.
The authors acknowledge the support by European Union’s Horizon 2020 research and innovation programme under the Marie Skłodowska-Curie Actions, grant agreement 872442 (ARIA).

\appendix 

\section{Hyper-reduction of port-to-bubble problems}
\label{sec:app_hyper_port2bubble}

We review the element-wise EQ hyper-reduction procedure that is employed here     to speed up the solution to the port-to-bubble problems.
The approach exploits the methods first proposed 
in \cite{farhat2015structure,yano2019lp}: we refer to \cite{taddei2021discretize} for further details.
We recall that, for any $\ell \in \mathcal{L}$,
$\{  x_{\ell, j}^{\rm a,v}  \}_{j=1}^{N_{\ell}^{\rm v}}$ are the nodes of the reference mesh of the $\ell$-th component, while 
$\texttt{T}_{\ell} \in \mathbb{N}^{N_{\ell}^{\rm e}  \times n_{\rm lp}}$ is the connectivity matrix, with 
$n_{\rm lp}$ equal to the number of elemental degrees of freedom.
We denote by $e_1,\ldots , e_D$ the vectors of the canonical basis in $\mathbb{R}^D$ and we denote by $\varphi_{\ell,k,i}$ the FE basis associated with the 
$i$-th degree of freedom of the $k$-th element of the $\ell$-th component.

Given $u\in \mathcal{X}_{\ell}^{\rm a}$,
we denote by $\mathbf{u}^{\rm un} \in \mathbb{R}^{n_{\rm lp} \times  N_{\ell}^{\rm e} \times D}$ the corresponding  third-order tensor such that
$$
u_{i,k,d}^{\rm un}
=
\left(
u\left(  x_{\ell, \texttt{T}_{\ell,k,i}}^{\rm a,v}  \right)
\right)_d,
\quad
i=1,\ldots,n_{\rm lp}, \;\;
k=1,\ldots,N_{\ell}^{\rm e}, \;\;
d=1,\ldots,D.
$$
Similarly, given the ROB basis 
$Z_{\ell}^{\rm a,b}:\mathbb{R}^n \to \mathcal{Z}_{\ell}^{\rm a,b}$, we denote by 
$\mathbf{Z}_{\ell}^{\rm a,b,un} \in \mathbb{R}^{n_{\rm lp} \times  N_{\ell}^{\rm e} \times D \times n}$ the corresponding fourth-order tensor. 
We further define the unassembled residual associated with the field $u$ and the parameter $\mu_{\ell}$,
$$
R_{\ell,i,k,d}^{\rm a,un}(u; \mu_{\ell}) \, : = \,
\int_{\texttt{D}_{\ell,k}} \;
\eta_{\ell}^{\rm a,e}  (u,  \varphi_{\ell,k,i} \;  e_d  ; \mu_{\ell}) \, dx
\; + \;
\int_{\partial \texttt{D}_{\ell,k}} \;
\eta_{\ell}^{\rm a,f} (u,  \varphi_{\ell,k,i} \; e_d  ; \mu_{\ell})
\, dx,
$$
for $
\ell\in \mathcal{L}, \; i=1,\ldots,n_{\rm lp}, \;
k=1,\ldots,N_{\ell}^{\rm e}, \;\; 
d=1,\ldots,D.$
Then, it is easy to verify that
\begin{subequations}
	\label{eq:residual_eq}
	\begin{equation}
		\label{eq:residual_eq_a}
		\left(
		\widehat{\mathbf{R}}_{\ell}^{\star}
		(\boldsymbol{\gamma}_{\ell}  )
		\right)_j
		\; = \;
		\sum_{i,k,d} \rho_{\ell,k}^{\star} \, 
		{Z}_{\ell,i,k,d,j}^{\rm a,b,un} \;\;
		R_{\ell,i,k,d}^{\rm a,un}(  \boldsymbol{\gamma}_{\ell}  )
		\; = \;
		\left( \mathbf{G}_{\ell}^{\rm a}(  \boldsymbol{\gamma}_{\ell}  ) \boldsymbol{\rho}_{\ell}^{\star}
		\right)_j,
		\;\;
		j=1,\ldots,n,
	\end{equation}
	where 
	$\star \in \{\rm hf, \rm eq\}$, 
	$\boldsymbol{\gamma}_{\ell} = (\boldsymbol{\alpha}_{\ell}, \boldsymbol{\beta}_{\ell}  , \mu_{\ell})$ denotes the triplet of bubble coefficients,  port coefficients and parameter, 
	$\mathbf{G}_{\ell}^{\rm a}: \mathbb{R}^n \times \mathbb{R}^m \times \mathcal{P}_{\ell} \to \mathbb{R}^{n \times N_{\ell}^{\rm e}}$ is the matrix-valued function that satisfies
	$
	\left(
	\mathbf{G}_{\ell}^{\rm a}(  \boldsymbol{\gamma}_{\ell}  )
	\right)_{j,k} 
	=
	\sum_{i,d}  
	{Z}_{\ell,i,k,d,j}^{\rm a,b,un} \;\;
	R_{\ell,i,k,d}^{\rm a,un}(  \boldsymbol{\gamma}_{\ell}  )
	$ for 
	$j=1,\ldots,n$ and $k=1,\ldots,N_{\ell}^{\rm e}$. The latter identity implies that
	\begin{equation}
		\widehat{\mathbf{R}}_{\ell}^{\rm hf}
		(\boldsymbol{\gamma}_{\ell}  )
		\, - \,
		\widehat{\mathbf{R}}_{\ell}^{\rm eq}
		(\boldsymbol{\gamma}_{\ell}  )
		\, = \,
		\mathbf{G}_{\ell}^{\rm a}(  \boldsymbol{\gamma}_{\ell}  ) \left( 
		\boldsymbol{\rho}_{\ell}^{\rm hf} - 
		\boldsymbol{\rho}_{\ell}^{\rm eq} 
		\right)
	\end{equation}
\end{subequations}

For any $\ell \in \mathcal{L}$, EQ procedures aim to find a vector $\boldsymbol{\rho}_{\ell}^{\rm eq} \in \mathbb{R}^{N_{\ell}^{\rm e}}$ such that
(i) $\boldsymbol{\rho}_{\ell}^{\rm eq}$ is as sparse as possible;
(ii) the constant function  is integrated accurately, that is
\begin{equation}
	\label{eq:constant_accuracy_EQ}
	\big|
	\sum_{k=1}^{N_{\ell}^{\rm e}}  {\rho}_{\ell,k}^{\rm eq}
	\, | \texttt{D}_{\ell, k}  |
	- | \Omega_{\ell}^{\rm a}   |
	\big| \ll 1;
\end{equation}
(iii) given the training set of  triplets
$ \Sigma_{\ell}^{\rm train,eq}  : =  \{  \boldsymbol{\gamma}_{\ell}^{(j)}  \}_{j=1}^{n_{\rm train,\ell}}$,  the  residual is adequately calculated for  all elements of the training set,
\begin{equation}
	\label{eq:manifold_accuracy_EQ}
	\big|
	\mathbf{J}_{\ell}^{\rm b} \left(  \boldsymbol{\gamma}_{\ell} \right)^{-1}
	\left(
	\widehat{\mathbf{R}}_{\ell}^{\rm hf}
	(\boldsymbol{\gamma}_{\ell} )
	\, - \,
	\widehat{\mathbf{R}}_{\ell}^{\rm eq}
	(\boldsymbol{\gamma}_{\ell}  )
	\right)
	\big|
	\ll 1,
	\quad
	{\rm where} \; 
	\mathbf{J}_{\ell}^{\rm b} : =
	\partial_{ \boldsymbol{\alpha}  }
	\widehat{\mathbf{R}}_{\ell}^{\rm hf},
	\quad
	\forall \, 
	\boldsymbol{\gamma}_{\ell} \in 
	\Sigma_{\ell}^{\rm train,eq}.
\end{equation}
As discussed in section \ref{sec:methods} (cf. 
\eqref{eq:importance_constant_accuracy}), the constant accuracy constraint
\eqref{eq:constant_accuracy_EQ}
is designed to control the $\ell^1$ norm of the weights that is related to the stability of the quadrature rule (see, e.g., \cite[section 2.3]{huybrechs2009stable}); the constraints
\eqref{eq:manifold_accuracy_EQ} are directly linked to the approximation error between the ROM estimate with HF quadrature and the hyper-reduced ROM estimate (cf. \cite[Proposition 3.2]{yano2019lp}).

We observe that the EQ problem can be recast  as a sparse representation problem of the form 
\begin{equation}
	\label{eq:sparse_representation_problem}
	\min_{\boldsymbol{\rho} \in \mathbb{R}^{N_{\ell}^{\rm e}}}
	\|   \boldsymbol{\rho}  \|_{\ell^0},
	\;\;
	{\rm s. \; t.} \;\;
	\|  \mathbf{C}_{\ell}^{\rm eq} \left(  \boldsymbol{\rho}_{\ell}^{\rm hf} - 
	\boldsymbol{\rho}_{\ell}^{\rm eq}  \right)
	\|_2\leq tol_{\rm eq},
\end{equation}
where  $\|   \boldsymbol{\rho}   \|_{\ell^0}$ is the $\ell^0$ norm that counts the number of non-zero entries in the vector $  \boldsymbol{\rho}$,
$\mathbf{C}_{\ell}^{\rm eq}$ is a suitable matrix that can be readily derived from \eqref{eq:constant_accuracy_EQ} and \eqref{eq:manifold_accuracy_EQ}, and  $tol_{\rm eq}$ is a suitable tolerance. Problem \eqref{eq:sparse_representation_problem} is NP hard; however, several effective approximate strategies have been proposed in the literature to determine 
parsimonious quadrature rules for MOR applications,
\cite{farhat2015structure,yano2019lp,chapman2017accelerated,manucci2022sparse}.
In this work, we resort to the non-negative least-square algorithm  implemented in the Matlab routine \texttt{lsqnonneg}, which takes as input the matrix $\mathbf{C}_{\ell}^{\rm eq}$, the vector $\mathbf{b}_{\ell}^{\rm eq}: = \mathbf{C}_{\ell}^{\rm eq} \boldsymbol{\rho}_{\ell}^{\rm hf} $ and the tolerance $tol_{\rm eq}$, and returns the sparse quadrature rule.

\section{Proofs}
\label{sec:tedious_proofs}

\subsection{Proof of   \Cref{th:link_port2full}}
\begin{proof}
	Let $u^{\star}$ be the solution to \eqref{eq:lin_coercive}.  Then, we find
	$$
	a(u^{\star}|_{\Omega_i}, v)
	\overset{\eqref{eq:technical_assumption}}{=}
	a(u^{\star} ,  v^{\rm ext})
	=
	f(  v^{\rm ext})
	=f(v) \;\; \forall \; v\in \mathcal{X}_{i,0};
	$$
	therefore, $u^{\star}|_{\Omega_i} = T_i u^{\star}|_{\Gamma_i} + G_i f$. The latter implies that $\lambda^{\star}  = \left( u^{\star}|_{\Gamma_1}, u^{\star}|_{\Gamma_2}   
	\right)$ satisfies
	$$
	\lambda^{\star} 
	- E \lambda^{\star}  - G f
	=
	\left(
	u^{\star}|_{\Gamma_1} - u^{\star}|_{\Gamma_1}, \;\;
	u^{\star}|_{\Gamma_2} - u^{\star}|_{\Gamma_2}
	\right) = 0,
	$$
	and thus that $\lambda^{\star}$ solves \eqref{eq:port_problem}.

	Let $\lambda^{\star}$ satisfy \eqref{eq:port_problem}. We define $u_i^{\star} = T_i  \lambda_i^{\star}  + G_i f$ for $i=1,2$. If we define the space
	$\mathcal{X}_{1,2}^0 = 
	\{ v|_{\Omega_1\cap \Omega_2} : v\in \mathcal{X}, \; v|_{\Gamma_1 \cup \Gamma_2} = 0 \}, 
	$ we observe that $u_1^{\star},u_2^{\star}$ satisfy 
	\begin{equation}
		\label{eq:intersection_problem}
		u_i^{\star}|_{\Gamma_1} = \lambda_1^{\star}, \;\;
		u_i^{\star}|_{\Gamma_2} = \lambda_2^{\star}, \;\;
		a(   u_i^{\star}, v   ) = f(v) \;\;\; \forall \; v \in \mathcal{X}_{1,2}^0,
		\qquad
		i=1,2.
	\end{equation}
	Since $a: \mathcal{X}_{1,2}^0 \times \mathcal{X}_{1,2}^0 \to \mathbb{R}$ is coercive, the solution to \eqref{eq:intersection_problem} exists and is unique: therefore, $u_1^{\star} = u_2^{\star}$ in $\Omega_1\cap \Omega_2$. In particular, if we define
	$u^{\star} = \sum_{i=1}^2 u^{\star}_i$, we have
	$u^{\star}|_{\Omega_i} = u_i^{\star}$ for $i=1,2$.
	
	Given $v\in \mathcal{X}$, we have 
	$v \phi_i \in \mathcal{X}_{i,0}$, since, by construction, $\rm{supp} \, \phi_i \subset \overline{\Omega}_i$. We thus have
	$$
	a(u^{\star},v)
	=
	\sum_{i=1}  a(u^{\star},\phi_i v)
	=
	\sum_{i=1}  a(u^{\star}\big|_{\Omega_i},\phi_i v)
	=
	\sum_{i=1}  a(u_i^{\star},\phi_i v)
	=
	\sum_{i=1}  f( \phi_i v)
	=
	f(v),
	$$
	which is the desired result.
\end{proof}

\subsection{Proof of Proposition  \ref{th:stability_port}}

\begin{proof}
	Continuity of $a_{\rm p}$ follows from the continuity of the trace operators, and the local operators $T_1,T_2$. We omit the details. To prove inf-sup  stability of the problem, we resort to the Fredholm's alternative: since $T$ is compact, provided that $\nu=1$ is not an eigenvalue of $T$, the equation $\lambda- T \lambda = f$ admits a unique solution for any $f\in \mathcal{U}$ and there exists a constant $C$ such that
	$\vertiii{\lambda} \leq C \vertiii{f}$ (see, e.g.,
	\cite[Theorem 6.6.8]{salsa2016partial}). It thus suffices to prove that $T \lambda = \lambda$ only holds for $\lambda=0$.

	Towards this end, we consider the problem:
	$$
	{\rm find} \; w\in \mathcal{X}_{1,2} : \,
	a(w,v) = 0 \; \forall \, v \in \mathcal{X}_{1,2}^0 ,
	\quad
	w|_{\Gamma_1} = \gamma_1, \;\;
	w|_{\Gamma_2} = \gamma_2,
	$$
	with 
	$\mathcal{X}_{1,2}^0 = 
	\{ v|_{\Omega_1\cap \Omega_2} : v\in \mathcal{X}  \}$, and 
	$\mathcal{X}_{1,2}^0 = 
	\{ v|_{\Omega_1\cap \Omega_2} : v\in \mathcal{X}, \; v|_{\Gamma_1 \cup \Gamma_2} = 0 \}$. 
	Since $T_1\lambda_1 = \lambda_1$ on $\Gamma_1$ by definition and 
	$T_1\lambda_1 = \lambda_2$ on $\Gamma_2$ since
	$T \lambda  = \lambda$, we have that 
	$T_1\lambda_1 |_{\Omega_{1,2}} = w$; similarly, we find   $T_2\lambda_2 |_{\Omega_{1,2}} = w$. 
	As observed in the proof of   \Cref{th:link_port2full}, there exists a unique solution to the problem 
	$w\in \mathcal{X}_{1,2}$: this implies that 
	$T_1\lambda_1 |_{\Omega_{1,2}}  = 
	T_2\lambda_2 |_{\Omega_{1,2}} $. Given the partition of unity $\phi_1,\phi_2$ associated with $\{ \Omega_i\}_{i=1}^2$, we define the field
	$u = \sum_{i=1}^2 \phi_i  T_i\lambda_i \in \mathcal{X}$, which satisfies
	$u|_{\Omega_i} = T_i\lambda_i$ for $i=1,2$.
	We observe that
	$$
	a(u, v)
	=
	\sum_{i=1}^2
	a\left(  u, \phi_i  v \right)
	\overset{\eqref{eq:technical_assumption}}{=}
	\sum_{i=1}^2
	a\left(  u|_{\Omega_i}, \phi_i  v\right)
	=
	\sum_{i=1}^2
	a\left(     T_i\lambda_i,  \phi_i   v\right)
	=
	0.
	$$
	Since $a$ is coercive, we must have $u\equiv 0$ and thus $\lambda\equiv 0$.
\end{proof}

\subsection{Proofs of the  estimate  \eqref{eq:perturbation_bound}}
\begin{proof}
	We first introduce the orthonormal basis $\{\psi_i\}_{i=1}^M$ of $\mathcal{Z}^{\rm p}$; given $\lambda\in \mathcal{Z}^{\rm p}$, we denote by
	$\boldsymbol{\lambda}\in \mathbb{R}^M$ the corresponding vector of coefficients such that
	$\lambda = \sum_{i=1}^M (\boldsymbol{\lambda} )_i \psi_i$. 
	By straightforward calculations, we find that 
	\begin{equation}
		\label{eq:algebraic_lincoerc}
		\widetilde{\mathbf{A}}  \, \widetilde{\boldsymbol{\lambda}} 
		\, = \, \widetilde{\mathbf{F}},
		\;\;
		\widehat{\mathbf{A}}  \, \widehat{\boldsymbol{\lambda}} 
		\, = \, \widehat{\mathbf{F}},
		\;\;
		{\rm with} 
		\left\{
		\begin{array}{l}
			\displaystyle{
				\left(
				\widetilde{\mathbf{A}} 
				\right)_{i,j} =
				\langle
				(Id - T) \psi_j, (Id - T) \psi_i
				\rangle,
				\;\;
				\left(
				\widetilde{\mathbf{F}} 
				\right)_{i} =
				\langle
				(Id - T) \psi_i,G f
				\rangle,
			}
			\\[3mm]
			\displaystyle{
				\left(
				\widehat{\mathbf{A}} 
				\right)_{i,j} =
				\langle
				(Id - \widehat{T}) \psi_j, (Id - \widehat{T}) \psi_i
				\rangle,
				\;\;
				\left(
				\widehat{\mathbf{F}} 
				\right)_{i} =
				\langle
				(Id - \widehat{T}) \psi_i,\widehat{G} f
				\rangle.
			}
			\\
		\end{array}
		\right.
	\end{equation}
	By straightforward calculations, we obtain
	$$
	\widetilde{\boldsymbol{\lambda}}  -
	\widehat{\boldsymbol{\lambda}} 
	\, = \,
	\widetilde{\mathbf{A}}^{-1} \left(
	\widetilde{\mathbf{F}}  - \widehat{\mathbf{F}}
	- 
	\left(
	\widetilde{\mathbf{A}}  - \widehat{\mathbf{A}}
	\right)
	\widehat{\boldsymbol{\lambda}} 
	\right)
	$$
	and thus
	\begin{equation}
		\label{eq:boring_evening_estimate}
		\|   \widetilde{\boldsymbol{\lambda}}  -
		\widehat{\boldsymbol{\lambda}}   \|_2
		\,  \leq  \, \underbrace{\|\widetilde{\mathbf{A}}^{-1}  \|_2}_{ \rm =:(I)}
		\Big(
		\underbrace{\|\widetilde{\mathbf{A}} - \widehat{\mathbf{A}}   \|_2}_{ \rm =:(II)}
		\underbrace{\|  \widehat{\boldsymbol{\lambda}}  \|_2}_{ \rm =:(III)}  
		\, + \,
		\underbrace{\| \widetilde{\mathbf{F}} - \widehat{\mathbf{F}} \|_2}_{ \rm =:(IV)}    
		\Big) .
	\end{equation}
	
	We estimate each term 
	of  \eqref{eq:boring_evening_estimate}
	independently: combination of the estimates for (I)-(IV) leads to \eqref{eq:perturbation_bound}.
	\begin{enumerate}
		\item[(I)] 
		Recalling the definition of $\alpha_{\rm p}$, we have $\vertiii{  \psi - T \psi    } \geq \alpha_{\rm p} \vertiii{\psi}$; therefore, we have
		$$
		\boldsymbol{\psi}^T \widetilde{\mathbf{A}} \boldsymbol{\psi}
		= 
		\vertiii{\psi - T \psi}^2
		\geq 
		\alpha_{\rm p}^2  \vertiii{\psi}^2
		=
		\alpha_{\rm p}^2
		\| \boldsymbol{\psi} \|_2^2,
		$$
		which implies (I).
		\item[(II)]
		By summing and subtracting 
		$\langle
		(Id - T) \psi_j, (Id - \widehat{T}) \psi_i
		\rangle$ to $
		\big|  \left(
		\widetilde{\mathbf{A}} 
		\right)_{i,j}  \, - \, 
		\left(
		\widehat{\mathbf{A}} 
		\right)_{i,j} 
		\big|$ and recalling the definitions of $\gamma_{\rm p}, \widehat{\gamma}_{\rm p}$ and $\varepsilon_{\rm T}$, we obtain
		$$
		\big|  \left(
		\widetilde{\mathbf{A}} 
		\right)_{i,j}  \, - \, 
		\left(
		\widehat{\mathbf{A}} 
		\right)_{i,j} 
		\big|
		\leq
		\left(
		\gamma_{\rm p} + \widehat{\gamma}_{\rm p}
		\right)  \varepsilon_{\rm T},
		\quad
		\forall \, i,j=1,\ldots,M.
		$$
		Estimate  (II) then follows  by exploiting the fact that 
		for any $M\times M$ matrix $\mathbf{A}$, we have
		$\|  \mathbf{A}\|_2 \leq M \max_{i,j} |   A_{i,j} |$.
		\item[(III)] Estimate (III) follows directly from the properties of minimum residual formulations of inf-sup stable problems. Indeed, since the bilinear form $a_{\rm p}$ is continuous and inf-sup stable, using the  Ne\v{c}as theorem
		(see, e.g., \cite[Thm 6.42]{salsa2016partial})	
		we have 
		$\|\widehat{\boldsymbol{\lambda}}\|_2=\vertiii{\widehat{{\lambda}}}\leq \frac{1}{\widehat{\alpha}_{\rm p}}  \vertiii{\widehat{G}f} $ for all $f \in \mathcal{X}'$.
		\item[(IV)] Proceeding as in (II), we find
		$$
		\big|  \left(
		\widetilde{\mathbf{F}} 
		\right)_{i}  \, - \, 
		\left(
		\widehat{\mathbf{F}} 
		\right)_{i}   \big|
		\leq
		\widehat{\gamma}_{\rm p} \varepsilon_{\rm G} +
		\vertiii{G f } \varepsilon_{\rm T},
		\quad
		\forall \, 
		i=1,\ldots,M,
		$$ 
		and thus $\|  \widetilde{\mathbf{F}} 
		\, - \, 
		\widehat{\mathbf{F}}  \|_2\leq\sqrt{M} 
		\|  \widetilde{\mathbf{F}} 
		\, - \, 
		\widehat{\mathbf{F}}  \|_{\infty}
		\leq
		\sqrt{M}  \left( \widehat{\gamma}_{\rm p} \varepsilon_{\rm G} +
		\vertiii{G f } \varepsilon_{\rm T}
		\right).
		$
	\end{enumerate}
\end{proof}

\subsection{Proof of \eqref{eq:Lions_variational_formulation_OS2}}
\begin{proof}
	For the two-subdomain problem,
	the OS2 statement \eqref{eq:OS2_lincoercive_mod} can be stated as:
	\begin{equation}
		\label{eq:one_more_OS2_lin}
		\min_{( \psi_1,\psi_2)\in \mathcal{Z}_1^{\rm p} \times \mathcal{Z}_2^{\rm p}}
		\,
		\|
		\widehat{u}_1(\psi_1) - \widehat{u}_2(\psi_2) 
		\|_{H^{1/2}(\Gamma_1\cup \Gamma_2)}
	\end{equation}
	where
	$\widehat{u}_i(\psi_i) = \widehat{u}_i^{\rm b}(\psi_i) + E_i \psi_i$ and 
	$\widehat{u}_i^{\rm b}(\psi_i)\in \mathcal{Z}_i^{\rm b}$ satisfies
	$
	a( \widehat{u}_i^{\rm b}(\psi_i) + E_i \psi_i, v ) = f(v)
	$ for all $v\in \mathcal{Z}_i^{\rm b}$
	and all $\psi_i\in \mathcal{Z}_i^{\rm p}$, for $i=1,2$. If we differentiate \eqref{eq:one_more_OS2_lin}, we obtain the optimality conditions
	$$
	\begin{array}{ll}
		\displaystyle{
			\left(
			\widehat{u}_1^{\rm p}
			-
			\chi_{\Gamma_1} \left(
			\widehat{u}_2^{\rm b} \left( \widehat{u}_2^{\rm p} \right)
			+
			E_2  \widehat{u}_2^{\rm p}
			\right), \; 
			\psi_1
			-
			\chi_{\Gamma_1} \left(
			\widehat{u}_2^{\rm b}(\psi_2)
			+
			E_2  \psi_2
			\right)
			\right)_{H^{1/2}(\Gamma_1)}
		} &
		\\[3mm]
		\displaystyle{
			+
			\left(
			\widehat{u}_2^{\rm p}
			-
			\chi_{\Gamma_2} \left(
			\widehat{u}_1^{\rm b} \left( \widehat{u}_1^{\rm p} \right)
			+
			E_1  \widehat{u}_1^{\rm p}
			\right), \; 
			\psi_2
			-
			\chi_{\Gamma_2} \left(
			\widehat{u}_1^{\rm b}(\psi_1)
			+
			E_1  \psi_1
			\right)
			\right)_{H^{1/2}(\Gamma_2)} }
		& 
		=0
		\quad
		\forall \, \psi=(\psi_1,\psi_2) \in {\mathcal{Z}}_1^{\rm p} \times {\mathcal{Z}}_2^{\rm p},
		\\
	\end{array} 
	$$
	which can rewritten as in  \eqref{eq:Lions_variational_formulation_OS2}.
\end{proof}

\subsection{Proofs of the estimates in \cref{sec:OSvsOS2_rates}}
In the following, we use the Taylor expansions:
\begin{equation}
	\label{eq:taylor_expansion}
	e^x \sim 1 + x + x^2, \quad
	\frac{1}{1-x} \sim 1 + x + x^2, \quad
	(1+x)^{1/2}  \sim 1 + \frac{1}{2} x - \frac{1}{8} x^2, \quad
	(1+x)^{2}  \sim 1 +  2 x,
\end{equation}
which are valid for $|x|\ll 1$.
We further employ the identiy:
\begin{equation}
	\label{eq:boring_identity}
	\max \left\{ |1 - \sigma \lambda_1  | , \;
	|1 - \sigma \lambda_2  |  \right\}
	=
	\left\{
	\begin{array}{ll}
		1 - \sigma \lambda_1 & \sigma < \frac{2}{\lambda_1 + \lambda_2}
		\\
		\sigma \lambda_2 - 1 & \sigma \geq  \frac{2}{\lambda_1 + \lambda_2}
		\\
	\end{array}
	\right.
\end{equation}
that is valid for any  $0 \leq \lambda_1 \leq \lambda_2$.

\subsubsection{Problem \eqref{eq:1Dpb_a}}
It is easy to verify that the local solutions 
$\widehat{u}_1,\widehat{u}_2$ satisfy
\begin{equation}
	\label{eq:local_solution_1a}
	\widehat{u}_1(x,\beta)
	=
	x^2 -
	\frac{ \delta^2}{1+\delta} (1+x) + 
	\frac{ \beta}{1+\delta} (1+x),
	\quad
	\widehat{u}_2(x,\beta)
	=
	x^2 -
	\frac{ \delta^2}{1+\delta} (1-x) + 
	\frac{ \beta}{1+\delta} (1+x).
\end{equation}
By imposing $\beta_1 = \widehat{u}_2(\delta,\beta_2)$ and
$\beta_2 = \widehat{u}_1(\delta,\beta_1)$ we obtain the system of equations:
$$
\mathbf{A}_{\delta} \, 
\boldsymbol{\beta}
\, = \, 
\mathbf{F}_{\delta},
\quad
{\rm with} \;
\mathbf{A}_{\delta} =\left[
\begin{array}{ll}
	1 &-c_{\delta} \\
	-c_{\delta} & 1 \\
\end{array}
\right],
\;\;
\mathbf{F}_{\delta} =\left[
\begin{array}{ll}
	d_{\delta} \\
	d_{\delta}   \\
\end{array}
\right],
$$
and $c_{\delta} = \frac{1-\delta}{1+\delta}$,
$d_{\delta} = \frac{2\delta^3}{1+\delta}$. The matrix $\mathbf{A}_{\delta}$ is symmetric with positive eigenvalues  $1-c_{\delta}$ and $1+c_{\delta}$; we thus have 
$$
{\rm cond} \left( 
\mathbf{A}_{\delta} \right)
= \frac{1+c_{\delta}}{1-c_{\delta}} = \frac{1}{\delta},
\;\;
\alpha_{\rm p}  = 
1-c_{\delta}
=
\frac{2 \delta}{1+\delta}, \;\;
\gamma_{\rm p}  = 
1+c_{\delta}
=
\frac{2 }{1+\delta}. 
$$
which are \eqref{eq:tedious_proof_b} and \eqref{eq:tedious_proof_c}.

Multiplicative OS corresponds to the application of the Gauss-Seidel iterative method to the linear system
$\mathbf{A}_{\delta} \, 
\boldsymbol{\beta}
\, = \, 
\mathbf{F}_{\delta}$. We thus find
$$
\boldsymbol{\beta}^{(k)}
=
\mathbf{P}_{\delta}^{\rm os} \boldsymbol{\beta}^{(k-1)}
+
\mathbf{F}_{\delta}^{\rm os},
\quad
{\rm with} \;\;
\mathbf{P}_{\delta}^{\rm os} = \left[
\begin{array}{ll}
	0 & c_{\delta} \\ 0 & c_{\delta}^2 \\
\end{array}
\right],
\;\;
\mathbf{F}_{\delta}^{\rm os} =\left[
\begin{array}{ll}
	d_{\delta} \\
	d_{\delta} + c_{\delta} d_{\delta}   \\
\end{array}
\right].
$$
We can then verify that the spectral radius of $\mathbf{P}_{\delta}^{\rm os}$ is equal to 
$$
\rho_{\delta}^{\rm os} =  c_{\delta}^2 \sim 1 - 4\delta.
$$

The OS2  method for \eqref{eq:1Dpb_a} reads as 
\begin{equation}
	\label{eq:tmp_OS2}
	\min_{\boldsymbol{\beta}  \in \mathbb{R}^2 }
	\frac{1}{2}
	\sum_{x \in \{ -\delta, \delta \}  }
	\left(
	\widehat{u}_1(x,\beta_1) - \widehat{u}_2(x,\beta_2)
	\right)^2
	=
	\frac{1}{2} \|  \mathbf{A}_{\delta} \, 
	\boldsymbol{\beta}
	- 
	\mathbf{F}_{\delta}  \|_2^2.
\end{equation}
If we apply the gradient descent method to \eqref{eq:tmp_OS2}, we obtain
$$
\boldsymbol{\beta}^{(k)}
=
\left( \mathbbm{1} - \sigma 
\mathbf{A}_{\delta}^T \mathbf{A}_{\delta} \right)  \boldsymbol{\beta}^{(k-1)}
+
\sigma 
\mathbf{A}_{\delta}^T \mathbf{F}_{\delta}.
$$
By tedious calculations, we can verify that the eigenvalues of the transition matrix
$\mathbbm{1} - \sigma 
\mathbf{A}_{\delta}^T \mathbf{A}_{\delta} $ are equal to 
$1 - \sigma ( c_{\delta} + 1  )^2$ and
$1 - \sigma ( c_{\delta} - 1  )^2$: recalling \eqref{eq:boring_identity}, we find that the spectral radius of the transition matrix is minimized by $\sigma = \frac{1}{c_{\delta}^2 + 1}$ and is equal to 
$$
\rho_{\delta}^{\rm os2} = \frac{2 c_{\delta}}{c_{\delta}^2 + 1}
\sim
1-4 \delta^2.
$$

\subsubsection{Problem \eqref{eq:1Dpb_b}}

The local solutions 
$\widehat{u}_1,\widehat{u}_2$ satisfy
\begin{equation}
	\label{eq:local_solution_1b}
	\widehat{u}_1(x,\beta)
	=
	\beta
	\frac{e^{\gamma x} - e^{-\gamma}}{e^{\gamma \delta} - e^{-\gamma}},
	\quad
	\widehat{u}_2(x,\beta)
	=
	\frac{e^{\gamma x} - e^{-\gamma \delta}  }{ e^{\gamma} - e^{-\gamma \delta} }
	\, + \,
	\beta
	\frac{e^{\gamma} - e^{ \gamma x}  }{ e^{\gamma} - e^{-\gamma \delta} }.
\end{equation}
Exploiting the Taylor expansions in 
\eqref{eq:taylor_expansion}, we obtain
$$
\widehat{u}_1(-\delta,\beta)
\sim
\beta \left(1 - 2 c_{\gamma} \delta  
+ 2 c_{\gamma}^2 \delta^2   \right),
\quad
\widehat{u}_2( \delta,\beta)
\sim
2  d_{\gamma} \delta - 
2  d_{\gamma}^2 \delta^2
\, + \,
\beta \left(
1 - 2  d_{\gamma} \delta
+ 2  d_{\gamma}^2 \delta^2
\right)
$$
where $c_{\gamma} := \frac{\gamma}{1 - e^{-\gamma}}$ and
$d_{\gamma} := \frac{c_{\gamma} }{ e^{\gamma}}$.
We thus find the (approximate) system of equations
$$
\mathbf{A}_{\delta} \, 
\boldsymbol{\beta}
\, = \, 
\mathbf{F}_{\delta},
\quad
{\rm with} \;
\mathbf{A}_{\delta} =\left[
\begin{array}{ll}
	1 &\left(-1 + 2 d_{\gamma} \delta - 2 d_{\gamma}^2 \delta^2 \right)   \\
	\left(-1 + 2 c_{\gamma} \delta - 2 c_{\gamma}^2 \delta^2 \right) 
	& 1 \\
\end{array}
\right],
\;\;
\mathbf{F}_{\delta} =\left[
\begin{array}{ll}
	2  d_{\gamma} \delta - 
	2  d_{\gamma}^2 \delta^2\\
	0   \\
\end{array}
\right].
$$
Therefore, the Gauss-Seidel transition matrix is approximately equal to
$$
\mathbf{P}_{\delta}^{\rm os} \sim \left[
\begin{array}{ll}
	0 & - 1 + 2 d_{\gamma} \delta \\ 0 & 
	-\left(1 - 2 d_{\gamma} \delta \right)
	\left(1 - 2 c_{\gamma} \delta \right) \\
\end{array}
\right]
$$
and thus 
$$
\rho_{\delta}^{\rm os} \sim 
1 - 2 \left(
c_{\gamma} + d_{\gamma} 
\right) \delta
=
1 - 2 \frac{e^{\gamma}+1}{e^{\gamma}-1} \gamma \delta .
$$
On the other hand, the eigenvalues of 
$\mathbf{A}_{\delta}^T \mathbf{A}_{\delta}$ are approximately equal to 
$$
\lambda_1 \sim
\frac{( c_{\gamma} + 2 d_{\gamma}   )^2}{4} \delta^2, \quad
\lambda_2 \sim
4 \, - \,
( 2 c_{\gamma} + 4 d_{\gamma}   ) \delta,
$$
and thus
$$
\alpha_{\rm p} = \sqrt{\lambda_1} 
\sim
\frac{4 (e^{\gamma}+2) \gamma\delta}{ 2 (e^{\gamma}-1)   },
\quad
\gamma_{\rm p} = \sqrt{\lambda_2} 
\sim 2.
$$

Exploiting \eqref{eq:boring_identity}, we find that the approximately optimal choice of the step size $\sigma$ is equal to $\sigma= \frac{1}{2} \left(  1 + \left( \frac{c_{\gamma}}{2} + d_{\gamma}\right) \delta \right)$ and thus
$$
\rho_{\delta}^{\rm os2}
\sim 1 - \sigma 
\lambda_1 \sim
1 - \frac{1}{8} \left(    c_{\gamma} + 2 d_{\gamma}      \right)^2 \delta^2.
$$
On the other hand, we obtain that the condition number of $\mathbf{A}_{\delta}$ is given by
$$
{\rm cond} ( \mathbf{A}_{\delta}  )
=
\sqrt{
	\frac{\lambda_{\rm max}(  \mathbf{A}_{\delta}^T \mathbf{A}_{\delta}  )}{
		\lambda_{\rm min}(  \mathbf{A}_{\delta}^T \mathbf{A}_{\delta}  )}
}
\sim
\frac{
	\sqrt{4 \, - \,
		( 2 c_{\gamma} + 4 d_{\gamma}   ) \delta}
}{ 
	\frac{( c_{\gamma} + 2 d_{\gamma}   )}{2} \delta     }
\sim
\frac{4}{ (c_{\gamma} + 2 d_{\gamma} )\delta  }
=
\frac{ 4 (e^{\gamma}-1)   }{4 (e^{\gamma}+2) \gamma\delta}
$$
\footnotesize 
 
\bibliographystyle{abbrv}
\bibliography{biblio}
 
\end{document}